
\documentstyle{amsppt}
\magnification=1200
\topmatter
\title Mixed Hodge Complexes on Algebraic Varieties
\endtitle
\author Morihiko Saito\endauthor
\affil RIMS Kyoto University, Kyoto 606-8502 Japan \endaffil
\endtopmatter
\tolerance=1200
\baselineskip=12pt
\hsize=6.5truein
\vsize=8.9truein

\document
\NoBlackBoxes
\def\ssbull{\raise.3ex\hbox{${\scriptscriptstyle\bullet}$}}
\def\Lmd{{\Lambda}}
\def\init{\text{{\rm init}}}
\def\ad{\text{{\rm ad}}\,}
\def\Perv{\hbox{{\rm Perv}}}
\def\Ext{\hbox{{\rm Ext}}}
\def\Gr{\text{{\rm Gr}}}
\def\Im{\hbox{{\rm Im}}}

\def\Coim{\hbox{{\rm Coim}}}
\def\Coker{\hbox{{\rm Coker}}}
\def\Ker{\hbox{{\rm Ker}}}
\def\Hom{\hbox{{\rm Hom}}}
\def\End{\hbox{{\rm End}}}

\def\DR{\hbox{{\rm DR}}}
\def\IC{\hbox{{\rm IC}}}

\def\Diff{\text{{\rm Diff}}}

\def\MHM{\text{{\rm MHM}}}
\def\MHW{\text{{\rm MHW}}}
\def\MH{\text{{\rm MH}}}
\def\MHS{\text{{\rm MHS}}} 
\def\Dec{\text{{\rm Dec}}\,}
\def\Decdual{{\text{{\rm Dec}}}^{*}}
\def\Supp{\hbox{{\rm Supp}}\,}
\def\supp{\hbox{{\rm supp}}\,}
\def\Sing{\hbox{{\rm Sing}}\,}

\def\inj{\text{\rm inj}}
\def\coh{\text{\rm coh}}
\def\hol{\text{\rm hol}}
\def\red{\text{\rm red}}

\def\an{\text{\rm an}}

\def\simto{\buildrel\sim\over\to}

\def\SameAuthor{\vrule height3pt depth-2.5pt width1cm}

\noindent
The notion of cohomological mixed Hodge complex was introduced by 
Deligne [4] as a tool to construct the mixed Hodge structure on the 
cohomology of complex algebraic varieties.
This is defined by using logarithmic forms and simplicial resolutions of 
varieties.
It is called cohomological, because its conditions are imposed only after 
taking the global section functor.
Forgetting the rational (or integral) structure and also the weight 
filtration, a refinement of this notion was obtained by Du Bois [6]
(following Deligne's idea).
For a complex algebraic variety
$ X $ and a closed subvariety
$ D $ of
$ X $, he introduced a filtered complex 
$ (\tilde{\Omega}^{\ssbull} _{X}\langle D \rangle , F) $ on
$ X $, whose differential is given by differential  operators of order 
at most one.
It is well defined in a certain triangulated category, and gives the Hodge 
filtration of Deligne's mixed Hodge structure on the cohomology of
$ X \backslash D $ by taking the hypercohomology if
$ X $ is proper.

On the other hand, the notion of mixed Hodge Module is introduced in [14], 
[15].
This gives also a mixed Hodge structure of the cohomology of 
a complex algebraic variety (without  using a simplicial resolution).
It is generally considered that the theory of mixed Hodge Modules is 
a generalization of Deligne's mixed Hodge theory.
There were, however, some gaps between the two theories.

One is that the theory of mixed Hodge Modules does not work on
simplicial schemes, because the direct images and the pull-backs are
defined only in the derived categories, and not in the level of complexes.
In particular, it was not clear whether the two mixed structures obtained
coincide in general (except when the variety is embeddable 
into a smooth variety).
The difficulty comes from the fact that the double complex construction 
associated with a cosimplicial complex on a simplicial scheme does not
work in the  derived category, because
$ d^{2} $ is not zero, but only homotopic to zero.
Note that this difficulty is not solved in this paper, nor seems
to be solved soon.
In fact, the problem is avoided in this paper, 
because we do not have to
construct a complex of mixed Hodge Modules on the simplicial scheme in
the  proof of (0.2).

Another problem is the difference in the systems of weight.
In Deligne's theory [4], the weight filtration 
$ W $ on a mixed Hodge complex
$ K $ is defined so that the weight of
$ H^{i}{\Gr}^{W}_{k} K $ is
$ i + k $.
This weight filtration 
$ W $ is called the {\it standard} weight.
On the other hand, the weight of
$ H^{i}{\Gr}^{W'}_{k}M $ is
$ k $ for a complex of mixed Hodge structures
$ M $ (or a mixed Hodge complex in the sense of [1, 3.2])
with the weight filtration
$ W' $.
This weight filtration
$ W' $ is called the {\it absolute} weight.
For a complex of mixed Hodge structures,
the passage from
$ W' $ to
$ W $ is done by taking the convolution of
$ W' $ and
$ \sigma $ (the filtration ``b\^ete'' in [4]).
But this does not work for a mixed Hodge complex in the sense of
[1, 3.2].
For the converse we have
$ W' = \Dec W $ as in [4].
This works for mixed Hodge complexes in the sense of loc.~cit.
These transformations of weights give certain difficulties for the 
unification of the two theories
(e.g.  in the compatibility with the direct images).
The main difficulty comes from the fact that the standard weight
filtration is not strict, but cannot be avoided in order that the
mixed Hodge complexes be stable by the direct images.

In this paper we construct a triangulated category of mixed 
$ {\Bbb Q} $-Hodge complexes (resp. mixed 
$ {\Bbb Q} $-Hodge
$ {\Cal D} $-complexes) on a complex algebraic variety
$ X $,
which will be denoted by
$ {D}_{\Cal H}^{b}(X,{\Bbb Q}) $ (resp.
$ {D}_{\Cal H}^{b}(X,{\Bbb Q})_{{\Cal D}}) $, together with a 
natural functor
$ \DR_{X}^{-1} : {D}_{\Cal H}^{b}(X,{\Bbb Q}) \rightarrow  {D}
_{\Cal H}^{b}(X,{\Bbb Q})_{{\Cal D}} $.

A mixed
$ {\Bbb Q} $-Hodge complex
$ K $ consists of a filtered differential complex 
$ (K_{F}, F) $ with a filtration
$ W $, and two filtered complexes
$ (K_{{\Bbb Q}},W) $,
$ (K_{\Bbb C},W) $ with
$ {\Bbb Q} $ and
$ {\Bbb C} $-coefficients.
They are endowed with two
quasi-isomorphisms
$ \alpha_{F} : (K_{F},W) \to (K_{\Bbb C},W) $,
$ \alpha_{{\Bbb Q}} : (K_{{\Bbb Q}},W)\otimes_{{\Bbb Q}}{\Bbb C} 
\to (K_{\Bbb C},W) $,
and satisfy certain good properties.
See (2.2).
A mixed
$ {\Bbb Q} $-Hodge 
$ \Cal D $-complex is defined similarly 
by replacing the filtered differential complex
$ (K_{F},F) $ with a complex of filtered
$ \Cal D $-Module
$( K_{\Cal D},F) $.
A morphism of mixed
$ {\Bbb Q} $-Hodge complex
$ u : K \to L $ consists not only of morphisms
$ u_{F} $,
$ u_{{\Bbb Q}} $,
$ u_{\Bbb C} $ between the components of
$ K $,
$ L $, but also of {\it homotopies}
$ u'_{F} $,
$ u'_{{\Bbb Q}} $, because
$ u_{*} $ and
$ \alpha_{*} $ commute up to a homotopy, which is given by
$ u'_{*} $.
Note that this homotopy does not appear in [1],
but is essential for Theorem (0.1) below.
We can show that the mixed
$ {\Bbb Q} $-Hodge complexes (or
$ {\Cal D} $-complexes) are stable by the direct image under a proper 
morphism of algebraic varieties in a compatible way with
$ \DR^{-1} $.

\bigskip\noindent
{\bf 0.1.~Theorem.} {\sl We have a natural functor
$ \varepsilon  $ from
$ D^{b}\MHM(X,{\Bbb Q}) $ (the derived category of 
bounded complexes of mixed
$ {\Bbb Q} $-Hodge Modules on
$ X) $ to
$ {D}_{\Cal H}^{b}(X,{\Bbb Q})_{{\Cal D}} $.
It commutes with the direct image under a proper morphism.
See (2.7-8).
}
\bigskip

For 
$ D $ a closed subvariety of
$ X $, and
$ j : X \backslash D \to X $ the inclusion morphism,
we construct a mixed
$ {\Bbb Q} $-Hodge complex
$ {\Cal C}^{H}_{X\langle D \rangle}{\Bbb Q} $ whose underlying
$ {\Bbb Q} $-complex is
$ \bold{R}j_{*}{\Bbb Q}_{X\backslash D} $ by using 
a simplicial resolution.
See (3.3).
Forgetting the weight filtration and the rational structure, it coincides
essentially with the Du Bois complex
$ (\tilde{\Omega}^{\ssbull} _{X}\langle D \rangle , F) $.

\bigskip\noindent
{\bf 0.2.~Theorem.} {\sl We have a natural isomorphism
$$
\DR_{X}^{-1}({\Cal C}^{H}_{X\langle D \rangle}{\Bbb Q}) = \varepsilon (j_{*}
{\Bbb Q}_{X\backslash D}^{H})\quad  \text{in }{D}_{\Cal H}^{b}(X,{\Bbb Q})
_{{\Cal D}},
$$
where
$ j_{*}{\Bbb Q}_{X\backslash D}^{H} \in  D^{b}\MHM(X) $ is as in [15].
See (2.8) below.
}
\bigskip

A naive idea for the proof of (0.2) would be to construct a complex of 
mixed Hodge Modules on a simplicial resolution of  X.
But this is very difficult, and is not enough to show (0.2), 
because there still remains the comparison between this complex with 
$ {\Cal C}^{H}_{X\langle D \rangle}{\Bbb Q} $, which is not trivial at all.
Instead, we use a smooth affine stratification of
$ X $ to get a good 
representative of
$ j_{*}{\Bbb Q}_{X\backslash D}^{H} $.
The corresponding object can be constructed in
$ {D}_{\Cal H}^{b}(X,{\Bbb Q})_{{\Cal D}} $,
and it is isomorphic to
$ \DR_{X}^{-1}({\Cal C}^{H}_{X\langle D \rangle}{\Bbb Q}) $.

As a corollary of (0.1) and (0.2),
we get the coincidence of the two mixed Hodge structures in [4] and [15].
But what is more important is the relation between
$ {\Bbb Q}_{X}^{H} $ and the Du Bois complex
$ (\tilde{\Omega}^{\ssbull} _{X}, F) $.
Let
$ (M,F) $ denote the underlying complex of filtered
$ {\Cal D} $-Module of
$ {\Bbb Q}_{X}^{H} $.
Then
$ {\Gr}_{p}^{F} \DR_{X}(M) $ are well-defined in
$ D^{b}_{\coh}(X, {\Cal O}_{X}) $.
See (1.3.4).
Let 
$ n = \dim X $, and
$ X' $ be the union of the
$ n $-dimensional irreducible components of
$ X $.
Let
$ \tilde{X}' \to X' $ be a resolution of singularities, and
$ \pi :  \tilde{X}' \to X $  its composition with the natural morphism to
$ X $.
Let
$ a_{\ssbull} : X_{\ssbull} \to X $ be a simplicial resolution.
Then (0.2) implies that
$ {\Gr}_{p}^{F} \DR_{X}(M) = 0 $ for
$ p > 0 $ or
$ p < -n $, and

\bigskip\noindent
{\bf 0.3.~Corollary.} {\sl
$ {\Gr}_{0}^{F} \DR_{X}(M) = 
\bold{R}(a_{\ssbull})_{*}{\Cal O}_{X_{\ssbull}} $, \,\,\,
$  {\Gr}_{-n}^{F} \DR_{X}(M) = \pi_{*}\omega_{\tilde{X}'}[-n] $.
}
\bigskip

We have
$ R^{i} \pi_{*}\omega_{\tilde{X}'} = 0 $ for
$ i > 0 $ by [8], and
$ \pi_{*}\omega_{\tilde{X}'} $ is independent of the choice of the
desingularization.
(See also [10].)
Let
$ (M',F) $  denote the underlying filtered
$ {\Cal D} $-Module of the intersection complex
$ \IC_{X'} {\Bbb Q}^{H} $.
Then we have also
$  {\Gr}_{-n}^{F} \DR_{X}(M') =
\pi_{*}\omega_{\tilde{X}'} $ by
[16].
Note that we cannot prove (0.3) without using (0.2)
even if
$ X $ is embeddable into a smooth variety.

We say that
$ X $ has (at most) Du Bois singularities, or 
$ X $ is Du Bois for short, if a natural morphism
$ {\Cal O}_{X} \to \bold{R} (a_{\ssbull})_{*}{\Cal O}_{X_{\ssbull}} $
is an isomorphism in the derived category.
See [20].
Of course,
$ {\Bbb C}_{X} \to \bold{R} (a_{\ssbull})_{*}{\Bbb C}_{X_{\ssbull}} $
 is always an isomorphism, 
but the same does not necessarily hold for 
$ {\Cal O} $.
This isomorphism gives a strong restriction to
$ {\Cal O}_{X} $.
(For example, 
$ X $ must be weakly normal. 
See Remark (i) after (5.2).
If
$ {\Bbb C}_{X}[\dim X] $ is a perverse sheaf,
$ X $ must be Cohen-Macaulay.
See Remark (ii) after (5.3).)

Let
$ K_{X} $ denote the dualizing complex for
$ {\Cal O}_{X} $-Modules, and
$ ({\Bbb D}M,F) $ the underlying filtered
$ {\Cal D} $-Module of
$ {\Bbb D} {\Bbb Q}_{X}^{H} $ (which is the dual of
$ {\Bbb Q}_{X}^{H} $) so that
$ {\Gr}_{0}^{F} \DR_{X}({\Bbb D}M) $ is the dual of
$ {\Gr}_{0}^{F} \DR_{X}(M) $ at least locally.
Then we have uniquely a morphism in
$ D_{\coh}^{b}(X,{\Cal O}_{X}) $
$$
{\Gr}_{0}^{F} \DR_{X}({\Bbb D}M) \to K_{X},
$$
whose restriction to the smooth part of
$ X $ is a natural isomorphism.
See (5.3).
By (0.3), this morphism is an isomorphism if and only if
$ X $ is Du Bois.
Using this, we can prove that a rational singularity is Du Bois,
as announced in Remark (ii) after (2.4) of [18].
This was conjectured, and proved in the isolated singularity case,
by Steenbrink [20].
After writing the first version of this paper,
I am informed that the assertion is proved by Koll\'ar [11,12.9]  
in the case where the  singularity admits a
projective compactification having only rational  singularities, 
and then by Kov\'acs [12] in the general case.

I would like to thank the referee for useful comments.

In \S1 we review the theory of filtered differential complexes and filtered
$ {\Cal D} $-Modules which are needed in this paper.
Then we introduce the notion of mixed Hodge 
complexes on algebraic varieties, and prove (0.1) in \S 2. 
The appendix to \S 2 gives a review on the theories of compatible
filtrations [14] and the realization functor [2], which are necessary for
the proof of (0.1).
In \S 3 the notion of geometric complexes is used to define a refinement of
the  filtered complex of Du Bois.
Then we prove (0.2) using a smooth affine stratification in \S 4.
The application to the Du Bois singularity is given in \S 5.

In this paper a variety means a separated scheme of finite type over a 
field, and an analytic space means a separated complex analytic space.

 \bigskip\bigskip\centerline{{\bf 1. Filtered Differential Complexes}}

\bigskip
\noindent
{\bf 1.1.} Let
$ X $ be an algebraic variety over a field
$ k $ of characteristic zero, or a complex analytic space.
(In the analytic case, variety will mean analytic space in this section.)
We consider filtered
$ {\Cal O}_{X} $-Modules
$ (L,F) $ such that the filtrations
$ F $ are increasing and exhaustive, and
$ F_{p}L = 0 $ for
$ p \ll 0 $.
(If the reader prefers, he can assume that
$ {\Gr}_{p}^{F}L $ are annihilated by the nilpotent ideal of
$ {\Cal O}_{X} $ corresponding to
$ X_{\red} \subset X $, so that
$ \DR_{X}^{-1}(L,F) $, which will be defined in (1.3),
satisfies the corresponding condition for filtered
$ {\Cal D} $-Modules.
This condition is satisfied by mixed Hodge Modules due to [14, 3.2.6].)

For
$ {\Cal O}_{X} $-Modules
$ L $,
$ L' $, we have the notion of differential operator of order
$ \le n $, which is invariant by the direct image of
$ L $,
$ L' $ under a closed embedding of varieties.
See Remark (i) below.
The group of differential operators of order
$ \le n $ are denoted by
$ \Diff_{X}^{n}(L,L') $.
If
$ X $ is smooth, this can be defined by the image of the composition
of injective morphisms (see [14, 2.2.2]) :
$$
\aligned
\Hom_{\Cal O_{X}}(L,L'\otimes_{\Cal O_{X}}F_{n}{\Cal D}_{X}) &\to
\Hom_{\Cal D}(L\otimes_{\Cal O_{X}}{\Cal D}_{X},
L'\otimes_{\Cal O_{X}}{\Cal D}_{X}) \\
&\to \Hom_{\Bbb C}(L,L'),
\endaligned
\tag 1.1.1
$$
where the last morphism is given by the tensor with
$ {\Cal O}_{X} $ over
$ {\Cal D}_{X} $.
(Here
$ {\Cal D}_{X} $ is the ring of differential operators on
$ {\Cal O}_{X} $, and has the Hodge filtration
$ F $ such that
$ F_{n}{\Cal D}_{X} = {\Diff}_{X}^{n}({\Cal O}_{X},{\Cal O}_{X}) $.)

For filtered
$ {\Cal O}_{X} $-Modules
$ (L,F) $,
$ (L',F) $ as above, we say that a 
$ {\Bbb C} $-linear morphism of filtered sheaves
$ u : (L,F) \to (L',F) $ is a filtered differential morphism if the composition
$$
F_{p}L \rightarrow  L \overset u \to \rightarrow L'
\rightarrow  L'/F_{p-q-1}L'
\tag 1.1.2
$$
is a differential operator of order  $ \le  q $ for any
$ p, q $.
(In particular, 
$ {\Gr}^{F}u : {\Gr}_{p}^{F}L \to {\Gr}_{p}^{F}L' $ is
$ {\Cal O}_{X} $-linear.)
The condition implies that if we put
$ \init(L,F) = \min\{p \in {\Bbb Z} : F_{p}L \ne 0 \} $, then
the restriction of
$ u $ to
$ F_{p}L $ is a differential operator of order
$ \le p - \init(L',F) $.
We will denote the group of filtered differential morphisms by
$ \Hom_{\Diff}((L,F),(L',F)) $.
This is invariant by the direct image of
$ (L,F) $,
$ (L',F) $ under a closed embedding of
varieties.

If
$ X $ is smooth, we define
$ \DR_{X}^{-1} (L,F) $ to be the right 
$ {\Cal D} $-Module 
$ L \otimes_{\Cal O_{X}}{\Cal D}_{X} $ with the tensor product filtration
$ F $, i.e.,
$ F_{p} = \sum _{i} F_{i}L\otimes _{{\Cal O}_{X}}F_{p-i}{\Cal D}_{X} $.
Then we see that the tensor with
$ {\Cal O}_{X} $ over
$ {\Cal D}_{X} $ induces a bijection
$$
\Hom_{\Cal D}(\DR_{X}^{-1}(L,F), \DR_{X}^{-1}(L',F)) \simto
\Hom_{\Diff}((L,F), (L',F)),
\tag 1.1.3
$$
where
$ \Hom $ on the left-hand side is taken in the category of filtered
right 
$ {\Cal D}_{X} $-Modules.
In fact, the morphism is injective by the injectivity of (1.1.1), 
and the condition on (1.1.2) corresponds to that the morphism of
$ {\Cal O}_{X} $-Modules 
$ L \to L'\otimes_{{\Cal O}_{X}}{\Cal D}_{X} $
preserves the filtration
$ F $, because
$$
F_{p}(L'\otimes_{\Cal O_{X}} {\Cal D}_{X}) =
\bigcap_{q} (F_{p-q-1}L'\otimes_{\Cal O_{X}} {\Cal D}_{X} +
L'\otimes_{\Cal O_{X}} F_{q}{\Cal D}_{X}).
$$
See [2, (3.1.2.8)].
(As remarked by the referee, we have the inverse morphism of (1.1.3) 
by identifying 
$ F_{p}\DR_{X}^{-1}(L) $ with
$ \Hom_{\Diff}(({\Cal O}_{X},F), (L,F[-p])) $ and taking the composition with
$ \Hom_{\Diff}((L,F),(L',F)) $.)
Note that (1.1.3) implies in general that the filtered differential morphisms 
are stable by composition.

We define
$ MF({\Cal O}_{X},\Diff) $ to be the category whose objects are
filtered
$ {\Cal O}_{X} $-Modules
$ (L,F) $ as above, and whose morphisms are filtered differential
morphisms.
In the smooth case, let
$ MF({\Cal D}_{X}) $ denote the category of filtered right
$ {\Cal D}_{X} $-Modules 
$ (M,F) $ such that the filtration
$ F $ is exhaustive and
$ F_{p}M = 0 $ for
$ p \ll 0 $.
Then, by the bijectivity of (1.1.3), we get a fully faithful functor
$$
\DR_{X}^{-1} : MF({\Cal O}_{X}, \Diff) \to MF({\Cal D}_{X}).
\tag 1.1.4
$$

Let
$ CF({\Cal O}_{X}, \Diff) $ be the category of complexes 
$ (L,F) $ of
$ MF({\Cal O}_{X}, \Diff) $, where we assume
$ F_{p}L = 0 $ for
$ p \ll 0 $ independently of the complex degree.
Then
$ KF({\Cal O}_{X}, \Diff) $ is defined by considering morphisms up to
homotopy as in [21].
Since the filtered acyclic objects form a thick subcategory in the sense
of loc.~cit, we get
$ DF({\Cal O}_{X}, \Diff) $ by inverting filtered quasi-isomorphisms.
Similarly we define
$ C^{*}F({\Cal O}_{X}, \Diff) $,  $ K^{*}F({\Cal O}_{X}, \Diff) $,  $ D^{*}F
({\Cal O}_{X}, \Diff) $ for
$ * = +, -, b $ by assuming the corresponding boundedness conditions 
on the complexes as in loc.~cit.

We consider also the category
$ MF({\Cal O}_{X},  \Diff; W) $ whose objects are
$ {\Cal O}_{X} $-Modules having a filtration
$ F $  as above together with an increasing finite filtration
$ W $, and whose morphism are filtered differential morphisms
with respect to $ F $, 
and preserve the filtration
$ W $.
Then for
$ * = +, -, b $ or empty, we can define the categories
$ C^{*}F({\Cal O}_{X}, \Diff; W) $,
$ K^{*}F({\Cal O}_{X}, \Diff; W) $, 
$ D^{*}F({\Cal O}_{X}, \Diff; W) $ in a similar way.
Here 
$ W $ is assumed to be finite uniformly in the complex degree.

We denote by
$ D_{\coh}^{*}F({\Cal O}_{X},\Diff;W) $ the subcategory of
$ D^{*}F({\Cal O}_{X},\Diff;W) $ consisting of objects
$ (L,F;W) $ such that
$ {\Cal H}^{i}{\Gr}_{p}^{F}{\Gr}_{k}^{W}L $ are coherent
$ {\Cal O}_{X} $-Modules.
Forgetting the filtration
$ W $, we have similarly
$ D_{\coh}^{*}F({\Cal O}_{X},\Diff) $.

\bigskip
\noindent
{\it Remarks.} (i) 
We have three equivalent definitions of differential operator of order
$ \le n $.
(The third definition is in the algebraic case.)

The first one takes (1.1.1) for definition in the smooth case,
and the general case is reduced to this by using locally defined closed 
embeddings into smooth varieties.
The well-definedness follows from the invariance of morphisms
of (filtered)
$ {\Cal D} $-Modules by the
direct image under a closed embedding of smooth varieties.
(In fact, if
$ i : X \to Y $ is a closed embedding of smooth varieties
such that 
$ X $ is defined locally by
$ y_{1} = 0 $ with
$ (y_{1}, \dots , y_{m}) $ a local coordinate system of
$ Y $, then
the direct image of a filtered right
$ {\Cal D}_{X} $-Module
$ (M,F) $ by
$ i $ is decomposed by the eigenvalue of the semisimple action of
$ y_{1}\partial/\partial y_{1} $ so that
$ M = \Ker\,y_{1}\partial/\partial y_{1} $.)

The second definition is due to Grothendieck [7, 16.8.1].
Let
$ I_{X} $ be the ideal of the diagonal of
$ X \times  X $,
and
$$
{P}_{X}^{n} = {\Cal O}_{X\times X }/ {I}_{X}^{n+1}.
$$
Let
$ L $,
$ L' $ be
$ {\Cal O}_{X} $-Modules, and
$ u : L \rightarrow  L' $ a
$ {\Bbb C} $-linear morphism of sheaves.
We say that
$ u $ is a differential operator of order
$ \le  n $ in the sense of Grothendieck if there exists an
$ {\Cal O}_{X} $-linear morphism
$$
u' : {P}_{X}^{n}\otimes _{{\Cal O}_{X}}L \rightarrow  L' 
$$
such that
$ u $ coincides with the composition of the natural inclusion
$$
L \rightarrow  {P}_{X}^{n}\otimes _{{\Cal O}_{X}}L 
\quad\text{(defined by }
 m \rightarrow  1\otimes m) 
$$
with
$ u' $.
Note that
$ u' $ is uniquely determined by
$ u $,
using a natural isomorphism
$$
\Hom_{{\Bbb C}}(L, L') = \Hom_{{\Cal O}_{X}}(({\Cal O}_{X}\otimes 
_{{\Bbb C}}{\Cal O}_{X})\otimes _{{\Cal O}_{X}}L, L'),
$$
because a natural morphism
$ {\Cal O}_{X}\otimes _{{\Bbb C}}{\Cal O}_{X} \rightarrow  {P}_{X}^{n} $ is 
surjective.
The last surjectivity can be reduced to the case
$ X $ smooth, because
$ P_{Y}^{n} \rightarrow P_{X}^{n} $ is surjective for a closed embedding
$ X \rightarrow Y $.  See also Remark after (1.2).
The equivalence of the first and the second definitions is reduced to
the smooth case (see Remark after (1.2)) and follows 
for example from [19, (1.20.2)].

The third definition in the algebraic case is as follows.
For
$ {\Cal O}_{X} $-Modules
$ L $,
$ L' $,
we say that  a 
$ {\Bbb C} $-linear morphism of sheaves
$ u : L \to L' $ is a differential morphism of order
$ \le n $, if for any local sections
$ g_{0}, \dots, g_{n} $ of
$ {\Cal O}_{X} $, we have
$$
\ad g_{0} \cdots \ad g_{n}(u) = 0,
\tag 1.1.5
$$
where
$ \ad g_{i}(v) = g_{i}v - vg_{i} $ for a morphisms of sheaves
$ v : L \to L' $.
(This was used for example in J. Bernstein's talk
on algebraic 
$ D $-Modules at the Luminy conference 1983.)
This definition is clearly invariant by the direct image under a
closed embedding of varieties, and we can verify the
equivalence with the second definition in the smooth case.
See for example [19, (1.20)].
(The key point is that
$ \ad g_{i} $ corresponds to the multiplication by
$ g_{i}\otimes 1 - 1\otimes g_{i} \in I_{X} $.)
By this definition, the condition on (1.1.2) means that 
for any local sections
$ g_{1}, \dots, g_{n} $ of
$ {\Cal O}_{X} $, 
$ \ad g_{1} \cdots \ad g_{n}(u) $
pushes down the filtration
$ F $ by
$ n $ (as remarked by the referee).
This  is also good for the direct image in (1.2),
(However the third definition is not so simple in the analytic case, 
because (1.1.5) is not strong enough, 
and we have to assume further a condition on the continuity of  
$ u $ in some topology.)

(ii) If
$ u : L \rightarrow  L' $ is a differential operator of order
$ \le  n $ and
$ L, L' $ have a filtration
$ W $ such that
$ u(W_{i}L) \subset  W_{i}L' $,
then
$$
u : W_{i}L/W_{j}L \rightarrow  W_{i}L'/W_{j}L'
$$
are differential operators of order
$ \le  n $ for
$ i > j $.
So if
$ u : (L,F) \rightarrow (L',F) $ is a filtered differential 
morphism preserving a filtration 
$ W $ of 
$ {\Cal O }_{X} $-Modules
$ L $ and
$ L' $, then
$$
u : (W_{i}L/W_{j}L, F) \rightarrow  (W_{i}L'/W_{j}L', F)
$$
is a filtered differential morphism for
$ i > j $.

(iii) It is not clear whether
$ D^{*}F({\Cal O}_{X},\Diff) $ is a full subcategory of
$ DF({\Cal O}_{X},\Diff) $ for
$ * = b, +, - $, because the usual truncation
$ \tau_{\le n} $ does not exist for filtered differential complexes.
However, 
$ D^{*}F^{b}({\Cal O}_{X},\Diff) $ (which is defined by assuming that the filtration
$ F $ is finite in the definition of
$ D^{*}F({\Cal O}_{X},\Diff) $) is a full subcategory of
$ DF^{b}({\Cal O}_{X},\Diff) $ and the latter is a full subcategory of
$ DF({\Cal O}_{X},\Diff) $ due to the functors which associate to
$ (L,F) $ respectively
$ ((\Dec F)_{p}L,F) $ and
$ (F_{p}L, F) $.
(Note that the
$ (\Dec F)_{p}L^{i} $ are 
$ {\Cal O}_{X} $-Modules because
$ {\Gr}^{F}d $ is
$ {\Cal O}_{X} $-linear.
We may consider
$ \Dec F $ as the canonical truncation relative to the filtration
$ F $.)

\bigskip
\noindent
{\bf 1.2.} Let
$ f : X \rightarrow  Y $ be a morphism of algebraic varieties or analytic 
spaces.
Then we have a natural morphism of
$ {\Cal O}_{Y} $-Modules
$$
{P}_{Y}^{n}\otimes _{{\Cal O}_{Y}}f_{*}L \rightarrow  
f_{*}{P}_{X}^{n}\otimes _{{\Cal O}_{Y}}f_{*}L \rightarrow  
f_{*}({P}_{X}^{n}\otimes _{{\Cal O}_{X}}L).
$$
So using the second definition of differential operators in Remark (i) after
(1.3),  we see that the direct image functor  
$ f_{*} $ for
$ {\Cal O} $-Modules induces
$$
\aligned
f_{*} : {\Diff}_{X}^{n}(L,L') &\rightarrow  {\Diff}_{Y}^{n}(f_{*}L,f_{*}L'), \\
f_{*} : \Hom_{\Diff}((L,F), (L',F)) &\rightarrow  \Hom_{\Diff}(f_{*}(L,F), 
f_{*}(L',F)),
\endaligned
$$
i.e., differential operators and filtered differential morphisms are stable 
by direct images.
(Here we can also verify this easily using the third definition of differential
operators in the algebraic case.)
So we get the direct image functor
$$
\bold{R}f_{*} : D^{+}F({\Cal O}_{X}, \Diff)  \rightarrow  D^{+}F
({\Cal O}_{Y}, \Diff)  .
\tag 1.2.1
$$
by using for example the canonical flasque resolution of Godement.
This direct image is compatible with the direct image of filtered
$ {\Cal D} $-Modules via the functor
$ \DR^{-1} $ in (1.3.2) below.
See [14, 2.3.6].

\bigskip

\noindent
{\it Remark.} Let
$ i : X \rightarrow  Y $ be a closed embedding.
Then the direct image induces natural isomorphisms
$$
\aligned
i_{*} : {\Diff}_{X}^{n}(L,L') 
&\simto {\Diff}_{Y}^{n}(i_{*}L,i_{*}L'),\\
i_{*} : \Hom_{\Diff}((L,F), (L',F)) 
&\simto \Hom_{\Diff}(i_{*}(L,F), 
i_{*}(L',F)).
\endaligned
\tag 1.2.2
$$
This is clear if we take the first definition of differential
operators (or the third in the algebraic case) in Remark (i) after (1.1).
For the second definition, the argument is as follows.
For the first isomorphism, it is enough to show that a natural morphism
$$
{\Cal O}_{X}\otimes _{{\Cal O}_{Y}}{P}_{Y}^{n}\otimes _{{\Cal O}_{Y}}
{\Cal O}_{X} \rightarrow  i_{*}{P}_{X}^{n}
$$
is an isomorphism.
Let
$ J, J' $ denote the defining ideals of
$ X $ and
$ X\times X $ in
$ Y $ and
$ Y\times Y $ respectively.
Since the source of the morphism is isomorphic to the quotient of
$ {\Cal O}_{Y\times Y} $ by the ideal generated by
$ J\otimes 1, 1\otimes J $ and
$ {I}_{Y}^{n} $,
it is enough to show the surjectivity of
$ I_{Y} \rightarrow  I_{X} $.
But this follows from the surjectivity of
$ J' \rightarrow  J $ using the snake lemma (applied to a morphism of the 
exact sequence associated with
$ J', J) $.
The second isomorphism follows from the first.
(Note that (1.2.2) is used in the proof of the equivalence of the three
definitions of differential operators 
when it is reduced to the smooth case.)

\bigskip
\noindent
{\bf 1.3.} Let
$ X $ be as in (1.1).
The notion of filtered
$ {\Cal D} $-Module is generalized to the case of
singular varieties as follows. 
See also [14, \S 2].
Let
$ LE(X) $ denote the category whose objects are closed embeddings
$ U \rightarrow  V $ where
$ U $ are open subsets of
$ X $ and
$ V $ are smooth.
A morphism of
$ \{ U \rightarrow  V\} $ to
$ \{ U' \rightarrow  V'\} $ in
$ LE(X) $ is a pair of morphisms
$ U \rightarrow  U' $ and
$ V $  $ \rightarrow  V' $ which gives a commutative diagram, where
$ U \rightarrow  U' $ is compatible with  natural inclusions
$ U \rightarrow  X $ and
$ U' \rightarrow  X $.

For
$ \{U \to V \} \in LE(X) $, let
$ MF({\Cal D}_{V})_{U} $ be the category of filtered right
$ {\Cal D}_{V} $-Modules
$ (M,F) $ such that the filtration
$ F $ is exhaustive,
$ F_{p}M = 0 $ for
$ p \ll 0 $, and the
$ {\Gr}_{p}^{F}M $ are annihilated by the ideal of
$ U $ in 
$ V $.
(If the reader prefers, he can assume that the
$ {\Gr}_{p}^{F}M $ are annihilated by the ideal of
$ U_{\red} $ rather than that of
$ U $.
This condition is satisfied by mixed Hodge Modules.)
By definition, a filtered (right)
$ {\Cal D} $-Module
$ (M,F) $ on
$ X $ consists of
$ (M_{U\rightarrow V},F) \in MF({\Cal D}_{V})_{U} $ for
$ \{ U \rightarrow  V\} \in  LE(X) $ together with isomorphisms
$$
v_{*}(M_{U\rightarrow V},F)|_{V'\backslash (U'\backslash U)} = 
(M_{U'\rightarrow V'},F)|_{V'\backslash (U'\backslash U)}
\tag 1.3.1
$$
for any morphisms
$ (u,v) $ of
$ \{ U \rightarrow  V\} $ to
$ \{ U' \rightarrow  V'\} $,
where
$ v_{*} $ denotes the direct image as filtered
$ {\Cal D} $-Modules (see [14, 2.3.5] and also Remark (iv) below), 
and (1.3.1) satisfies the natural
compatibility  condition for the composition of morphisms of
$ LE(X) $.
(Note that the left-hand side of (1.3.1) is a filtered
$ {\Cal D} $-Module.
See Remark (i) below.)
We assume also
$ F_{p}M_{U\rightarrow V} = 0 $ for
$ p \ll 0 $ independently of
$ \{ U \to V \} $.

Let
$ \init \,(= \init(M,F)) = \min\{p \in {\Bbb Z} : F_{p}M_{U\to V}\ne 0\} $.
Then we see that the
$ {\Cal O}_{U} $-Modules
$ F_{\init}M_{U\to V} $ are independent of the embedding of
$ U $, and defines a well-defined
$ {\Cal O}_{X} $-Module, which is denoted by
$ F_{\init}M $.
(This does not hold if we use left 
$ {\Cal D} $-Modules.)

The category of filtered
$ {\Cal D} $-Modules on
$ X $ is denoted by
$ MF(X,{\Cal D}) $, and
the category of filtered
$ {\Cal D} $-Modules with a finite increasing filtration
$ W $ by
$ MF(X,{\Cal D};W) $.
For 
$ * = +, -, b $ or empty, we denote by
$ C^{*}F(X,{\Cal D}) $, 
$ C^{*}F(X,{\Cal D};W) $
the categories of complexes 
of these objects satisfying the corresponding boundedness condition.
Their derived categories are denoted by
$ D^{*}F(X,{\Cal D}) $,
$ D^{*}F(X,{\Cal D};W) $.
Here we assume
$ \min\{\init(M^{i},F)\} \ne -\infty $ for a filtered complex
$ (M,F) $, and
$ W $ is finite uniformly in the complex degree.
We have similarly
$ M(X,{\Cal D}) $,
$ D^{*}(X,{\Cal D}) $,
$ D^{*}(X,{\Cal D};W) $,
forgetting the filtration
$ F $, but assuming that any local section of
$ M_{U\to V} $ is annihilated by a power of
the ideal of
$ U $ (this condition can be replaced with
$ \Supp M_{U\to V} \subset U $ if
$ M_{U\to V} $ is quasi-coherent).

We will denote by
$ D^{*}_{\coh} F(X,{\Cal D};W) $ the full subcategory of
$ D^{*} F(X,{\Cal D};W) $
consisting of coherent complexes  (i.e.
$ \bigoplus_{p} H^{i} F_{p} {\Gr}_{k}^{W} M_{U \to V} $
are coherent over
$ \bigoplus_{p} F_{p} {\Cal D}_{V} $), and by
$ D^{*}_{\hol} F(X,{\Cal D};W) $ its full subcategory consisting of
holonomic complexes (i.e.,  furthermore
$ H^{i} {\Gr}_{k}^{W} M_{U \to V} $ are holonomic
$ {\Cal D}_{V} $-Modules).
Forgetting the filtration
$ W $, we define similarly
$ D^{*}_{\coh} F(X,{\Cal D}) $,
$ D^{*}_{\hol} F(X,{\Cal D}) $.
We have also
$ M_{\hol} (X,{\Cal D}) $,
$ D^{*}_{\hol} (X,{\Cal D}) $,
$ D^{*}_{\hol} (X,{\Cal D};W) $, etc. as usual,
forgetting the filtration
$ F $.

We have a functor
$$
{\DR}_{X}^{-1} :MF({\Cal O}_{X}, \Diff)  \rightarrow  MF({\Cal D}_{X})
$$
defined by
$ (L,F) \to \{ (M_{U \to V},F) \} $ with
$ (M_{U \to V},F) = (L,F)|_{U}\otimes _{{\Cal O}_{V}}({\Cal D}_{V},F)  $.
This induces
$$
{\DR}_{X}^{-1} :DF({\Cal O}_{X}, \Diff)  \rightarrow  DF({\Cal D}_{X}).
\tag 1.3.2
$$

Let
$ \{ U_{i} \rightarrow  V_{i}\}_{i\in \Lambda} \in  LE(X) $ 
(indexed by an ordered set
$ \Lambda  $) such that  
$ \{ U_{i}\}_{i\in \Lambda } $ is a locally finite covering of
$ X $.
For
$ I \subset  \Lambda  $,
let
$ U_{I} = \bigcap _{i\in I} U_{i}, V_{I} = \prod _{i\in I} V_{i} $ so that
$ \{ U_{I} \rightarrow  V_{I}\} \in  LE(X) $.
Let
$ M \in  M(X,{\Cal D}) $,
and put
$$
K_{I} = (j_{I})_{!}\DR_{V_{I}}(M_{U_{I}\rightarrow V_{I}}) 
\quad\text{(resp. }
(j_{I})_{!}\DR_{V_{I}}(M_{U_{I}\rightarrow V_{I}})^{\an})
$$
if
$ X $ is analytic (resp. algebraic).
(For
$ \DR $ see Remark (iii) below.)
Here
$ j_{I} : U_{I} \rightarrow  X $ (resp.
$ j_{I} : {U}_{I}^{\an} \rightarrow  X^{\an}) $ denote natural inclusions.

For
$ I \supset J $,
$ \DR_{V_{I}}(M_{U_{I} \to V_{I}}) $ has a double complex structure
associated with the decomposition
$ V_{I} = V_{J}\times V_{I\backslash J} $.
Let
$ v : V_{I} \to V_{J} $ denote the natural projection, and
$ v_{\ssbull} $ the sheaf-theoretic direct image.
Then
$$
v_{*}(M_{U_{I}\to V_{I}},F)|_{V_{J}\backslash(U_{J}\backslash U_{I})} = 
v_{\ssbull} \DR_{V_{I}/V_{J}}(M_{U_{I}\to V_{I}},F)
|_{V_{J}\backslash(U_{J}\backslash U_{I})} 
$$
as filtered right
$ {\Cal D} $-Modules, and
$$
\DR_{V_{J}}(v_{\ssbull} \DR_{V_{I}/V_{J}}(M_{U_{I}\to V_{I}},F))=
v_{\ssbull} \DR_{V_{I}}(M_{U_{I}\to V_{I}},F) .
$$
Since we have a filtered quasi-isomorphism
$$
v_{\ssbull} \DR_{V_{I}/V_{J}}(M_{U_{I}\to V_{I}},F)
|_{V_{J}\backslash(U_{J}\backslash U_{I})} \to (M_{V_{J}\to V_{J}},F)
|_{V_{J}\backslash(U_{J}\backslash U_{I})} 
$$
by (1.3.1),
we get a natural morphism
$ K_{I} \rightarrow  K_{J} $ which induces a quasi-isomorphism on
$ U_{I} $.
So we get a co-Cech complex whose
$ i $-th component is
$ \bigoplus _{|I| = 1-i} K_{I} $,
and it is denoted by
$ \DR_{X}(M) $.
This induces functors
$$
\aligned
\DR_{X} :D^{*}_{\hol} (X,{\Cal D}) &\rightarrow 
 D^{*}_{c} (X,{\Bbb C}), \\
\DR_{X} : M_{\hol} (X,{\Cal D}) &\rightarrow 
\Perv (X,{\Bbb C}),
\endaligned
\tag 1.3.3
$$
where
$ D_{c}^{*}(X,{\Bbb C}) $ in the algebraic case denotes the full
subcategory of 
$ D_{c}^{*} (X^{\an} ,{\Bbb C}) $ consisting of objects with
algebraic stratifications, and similarly for
$ \Perv (X,{\Bbb C}) $
(see [2] for
$  D^{*}_{c} (X,{\Bbb C}) $,
$ \Perv (X,{\Bbb C}) $ in the analytic case).

We have similarly functors
$$
{\Gr}_{p}^{F}\DR_{X} :D^{*}_{\coh} F(X,{\Cal D}) \rightarrow 
 D^{*}_{\coh} (X,{\Cal O}_{X}),
\tag 1.3.4
$$
and we can verify
$ {\Gr}_{\init}^{F} \DR_{X}(M) = F_{\init}M $.
But we do not have
$ \DR_{X} :D^{*}_{\coh} F(X,{\Cal D}) \rightarrow 
 D^{*}_{\coh} F({\Cal O}_{X},\Diff) $ unless
$ X $ is smooth.

We can verily that (1.3.3-4) are independent of the choice of the covering
$ \{ U_{i} \rightarrow  V_{i}\}_{i\in \Lambda } $ by using
$ \{ U_{i} \rightarrow  V_{i}\}_{i\in \Lambda ''} $ with
$ \Lambda '' $ the disjoint union of
$ \Lambda  $ and
$ \Lambda ' $ if
$ \{ U_{i} \rightarrow  V_{i}\}_{i\in \Lambda '} $ is another covering of
$ X $.

\bigskip

\noindent
{\it Remarks.} (i) The left-hand side of (1.3.1) is a filtered
$ {\Cal D} $-Module.
In fact, for
$ x \in U $, we take a minimal embedding of
$ (U,x) $ into a smooth variety
$ (Z,x) $ (i.e. 
$ \dim Z $ coincides with the dimension of the
Zariski tangent space of 
$ U $ at 
$ x $).
Then
$ (U,x) \to (V,x) $ is the composition of
$ (U,x) \to (Z,x) $ with a closed embedding
$ i_{V} : (Z,x) \to (V,x) $, and 
$ (M_{U\to V},F)_{x}  = (i_{V})_{*}(M_{Z},F)_{x} $ for a filtered
$ {\Cal D}_{Z} $-Module
$ (M_{Z}, F) $.
(In fact, if 
$ Z $ is defined by
$ y_{1} = 0 $ with
$ (y_{1}, \dots, y_{m}) $ a local coordinate system of
$ (V,x) $, then the filtration
$ F $ on
$ M_{U\to V} $ is stable by the action of
$ y_{1}\partial/\partial y_{1} $, and is compatible with
the decomposition of
$ M_{U\to V} $ by the eigenvalue of the action of 
$ y_{1}\partial/\partial y_{1} $.)
Since the composition of 
$ i_{V} $ with
$ v $ is a closed embedding
$ i_{V'} : (Z,x) \to (V',x) $, we get the assertion.

(ii) Let 
$ \{U \to V_{i}\} \in LE(X) $ for
$ i = 1, 2 $, and
$ v : V_{1} \to V_{2} $ a morphism giving a morphism of
$ LE(X) $ (i.e., it is compatible with 
$ \{U \to V_{i}\} $).
Then the direct image of filtered
$ {\Cal D} $-Modules induces an equivalence of categories
$$
v_{*} : MF({\Cal D}_{V})_{U} \to MF({\Cal D}_{V'})_{U}.
$$
This follows by the same argument as 
the above Remark (i) by reducing to the case where
$ (U,x) \to (V_{1},x) $ is a minimal embedding.

If there is no morphisms between
$ V_{1} $ and
$ V_{2} $, we can consider
$ V_{3} := V_{1} \times V_{2} $ with the diagonal morphism
$ U \to V_{3} $ so that we get morphisms
$ V_{3} \to V_{i} $ by the projection.
So we see for example that if 
$ X $ is embedded into a smooth variety
$ Z $, then
$ MF(X, {\Cal D}) $ is naturally equivalent to
$ MF({\Cal D}_{Z})_{X}  $.

(iii) For a right
$ {\Cal D}_{X} $-Module
$ M $ on a complex manifold
$ X $ or a smooth algebraic variety
$ X $,
the  de Rham complex
$ \DR_{X}(M) $ is defined so that its
$ i $-th component is
$ M\otimes _{{\Cal O}_{X}}\wedge ^{-i}\Theta _{X} $.
If
$ M $ has a filtration
$ F $,
$ \DR_{X}(M) $ has the filtration
$ F $ such that
$ F_{p} $ on the
$ i $-th component is 
$ F_{p+i}M\otimes _{{\Cal O}_{X}}\wedge ^{-i}\Theta _{X} $.

(iv)  If 
$ Y $ is a complex manifold, or a smooth algebraic variety,
$ \DR_{X\times Y/Y} $ is defined similarly.
We can verify that for a complex of filtered
$ {\Cal D}_{X\times Y} $-Modules
$ (M,F) $, the direct image 
$ p_{*}(M,F) $ by the projection
$ p : X\times Y \to Y $ is isomorphic to
$ \bold{R}p_{\ssbull}\DR_{X\times Y/Y}(M,F) $.

If
$ i : X \to Y $ is a closed embedding of complex manifolds or
smooth varieties such that locally
$ X = \{y_{1} = \dots = y_{r} = 0 \} $ by using a local coordinate system
$ (y_{1}, \dots, y_{m}) $, let
$ \partial_{i} = \partial/\partial y_{i} $ for 
$ 1 \le i \le r $,
$ {\Bbb C}[\partial] = {\Bbb C}[\partial_{1}, \dots, \partial_{r}] $ and
$ \partial^{\nu} = \prod_{1 \le i \le r} \partial_{i}^{\nu_{i}} $ for
$ \nu = (\nu_{1}, \dots , \nu_{r}) \in {\Bbb N}^{r} $.
Then we can show the isomorphism
$$
i_{*}(M,F) = i_{\ssbull}M\otimes_{\Bbb C}{\Bbb C}[\partial]
\quad \text{with }
F_{j}i_{*}M = \sum_{\nu} i_{\ssbull}F_{j-|\nu|}M\otimes \partial^{\nu}.
$$

\noindent
{\bf 1.4.} We review here some properties of  Hodge Modules which
are needed in this paper.
Let
$ A $ be a subfield of
$ {\Bbb R} $.
We denote by
$ MF_{\hol}(X,{\Cal D},A;W) $
the category whose objects are
$ {\Cal M} = ((M,F;W),(K,W);\alpha) $, where
$ (M,F;W) $ is a holonomic filtered
$ {\Cal D} $-Module on
$ X $ with a finite increasing filtration
$ W $,
$ (K,W) $ is an
$ A $-perverse sheaf on
$ X $ with a finite increasing filtration
$ W $, and
$$
\alpha : \DR_{X}(M,W) \simeq (K,W) \otimes_{A}{\Bbb C}
$$
is an isomorphism of filtered 
$ {\Bbb C} $-perverse sheaves.
Forgetting the filtration
$ W $, the category
$ MF_{\hol}(X,{\Cal D},A) $ is similarly defined, and we have
functors
$$
\Gr_{k}^{W} : MF_{\hol}(X,{\Cal D},A;W) 
\to MF_{\hol}(X,{\Cal D},A).
$$
Sometimes an object of
$ MF_{\hol}(X,{\Cal D},A;W) $ will be denoted by
$ (M,F,K;W) $ (and similarly for
$ (M,F,K) $) to simplify the notation.

We denote by
$ \MH(X,A,n) $  the full subcategory of
$ MF_{\hol}(X,{\Cal D},A) $
consisting of polarizable
$ A $-Hodge Modules of weight
$ n $ in the sense of [14, \S 5],
and
$ \MHW(X,A) $  the full subcategory of
$ MF_{\hol}(X,{\Cal D},A;W) $
consisting of 
$ {\Cal M} $ such that
$ \Gr_{k}^{W} {\Cal M} \in \MH(X,A,k) $ for any
$ k $.
This means that the objects of
$ \MHW(X,A) $ are obtained by extensions of polarizable Hodge Modules, 
but there is no condition on the extensions.
The category
$ \MHM(X,A) $ of mixed 
$ A $-Hodge Modules on
$ X $ is a full subcategory of
$ \MHW(X,A) $,
and its objects satisfy certain good properties.
See [15].
We say that 
$ {\Cal M} \in \MHW(X,A) $ is pure of weight
$ n $ if
$ \Gr_{k}^{W}{\Cal M} = 0 $ for
$ k \ne n $.

We state here some properties of Hodge Modules used in this paper.

(i) A polarizable 
$ A $-Hodge Module admits a strict support decomposition
$$
(M,F,K) = \bigoplus _Z (M_{Z},F,K_{Z} ),
\tag 1.4.1
$$
with
$ Z $ running over the irreducible closed subvarieties of
$ X $, where
$ M_{Z} $ is a holonomic 
$ {\Cal D} $-Module on
$ X $ with strict support
$ Z $ (i.e. 
$ \supp M_{Z} = Z $ or 
$ \emptyset $, and 
$ M_{Z} $ has no nontrivial sub or quotient object whose
support  has dimension 
$ < \dim Z $), and
$ K_{Z} $ is an intersection complex with support
$ Z $ or
$ \emptyset $ (in the sense of [2]).
See [14, 5.1.6].

(ii) If a polarizable
$ A $-Hodge Module of weight
$ n $ has strict support
$ Z $, its restriction to a smooth open subvariety of
$ Z $ is a polarizable variation of 
$ A $-Hodge structure of weight
$ n - \dim Z $.
See [14, 5.1.10 and 5.2.12]
Conversely, any polarizable variation of
$ A $-Hodge structure of weight
$ n - \dim Z $ on a smooth open subvariety of
$ Z $ extends uniquely to a polarizable 
$ A $-Hodge Module of weight
$ n $ with strict support
$ Z $.
See [15, 3.21].

(iii) $ \MHW(X,A) $ and
$ \MHM(X,A) $ are abelian categories whose morphisms are 
strictly compatible with the two filtrations
$ F $,
$ W $ in the  sense of [14, 1.1.15] (see also (A.2) below).
For polarizable 
$ A $-Hodge Modules
$ {\Cal M} $,
$ {\Cal M}' $ of weights
$ n $,
$ n' $ such that
$ n > n' $, we have
$$
\Hom({\Cal M},{\Cal M}') = 0,
\tag 1.4.2
$$
where
$ \Hom $ is taken in 
$ MF_{\hol}(X,{\Cal D},A) $.

(iv) Let
$ f : X \to Y $ be a proper morphism of algebraic varieties, or
a projective morphism of analytic spaces, and
$ (M,F,K) $ a polarizable 
$ A $-Hodge Module of weight
$ n $ on
$ X $.
Then the direct image
$ f_{*}(M,F) $ as filtered
$ {\Cal D} $-Module is strict, and
$$
({\Cal H}^{i}f_{*}(M,F), {}^{p}{\Cal H}^{i}f_{*}K) \in \MH(Y,A,n+i).
\tag 1.4.3
$$
Furthermore there exist (noncanonical) isomorphisms 
$$
\aligned
f_{*} (M,F) &\simeq \bigoplus _{i} 
({\Cal H}^{i}f_{*} (M,F))[-i]\quad  \text{in }D^{+}F(Y, {\Cal D}), \\
f_{*}K &\simeq \bigoplus _{i} ({}^{p}
{\Cal H}^{i}f_{*}K)[-i]\quad  \text{in }D^{+}(Y,A)
\endaligned
\tag 1.4.4
$$
compatible with
$ \alpha  $.

In fact, (1.4.3) follows from the stability of polarizable Hodge Modules 
under a proper morphism in the algebraic case [15] or a projective 
morphism in the analytic case [14].
Then (1.4.4) follows from the hard Lefschetz theorem [14] together 
with [5] in the analytic case.
In the algebraic case, we can verify easily the compatibility 
with the direct 
image of mixed Hodge Modules in [15], forgetting the filtration
$ W $.
Then the desired decomposition follows from the decomposition of the 
direct image of a pure Hodge Module.
(It is also possible to reduce the assertion to the projective case
by using Chow's lemma together with the above property (ii).)

(v)  Let
$ X $ be a smooth algebraic variety or complex manifold, and
$ D $ a divisor with normal crossings on
$ X $.
Put
$ U = X \backslash D $.
Let
$ (j^{\Cal D}_{*}{\Cal O}_{U},F;W') $ denote the underlying filtered
$ {\Cal D}_{X} $-Module
 of
$ j_{*}A_{U}^{H} $, where
$ j : U \to X $ denotes a natural inclusion.
(In the algebraic case, we have
$ j^{\Cal D}_{*}{\Cal O}_{U} = j_{*}{\Cal O}_{U} $.)
Then we have a canonical isomorphism in
$ D_{\hol}^{b}F(X,{\Cal D};W) $ :
$$
(j^{\Cal D}_{*}{\Cal O}_{U},F;W) = 
\DR_{X}^{-1}(\Omega_{X}^{\ssbull}(\log D),F;W)[\dim X],
\tag 1.4.5
$$
where
$ W = W' $ on
$ j^{\Cal D}_{*}{\Cal O}_{U} $ (because
$ j^{\Cal D}_{*}{\Cal O}_{U} $ is a
$ {\Cal D}_{X} $-Module), and the filtrations
$ F $,
$ W $ on the right-hand side is as in [4].
In fact, using [15, 3.11], we can show the bifiltered quasi-isomorphism
$$ 
(\Omega_{X}^{\ssbull}(\log D),F;W)[\dim X] \simto
\DR_{X} (j^{\Cal D}_{*}{\Cal O}_{U},F;W).
$$

Similarly, let
$ (j^{\Cal D}_{!}{\Cal O}_{U},F;W') $ denote the underlying filtered
$ {\Cal D}_{X} $-Module
 of
$ j_{!}A_{U}^{H} $.
Then we have a canonical  isomorphism in
$ D_{\hol}^{b}F(X,{\Cal D};W) $ :
$$ 
(j^{\Cal D}_{!}{\Cal O}_{U},F) = 
\DR_{X}^{-1}(\Omega_{X}^{\ssbull}(\log D)(-D),F)[\dim X],
\tag 1.4.6
$$
where
$ (-D) $ means the tensor with the sheaf of reduced ideals of
$ D $.

\bigskip\bigskip
\centerline{{\bf 2. Mixed Hodge Complexes}}

\bigskip
\noindent
In this section, we construct triangulated categories of mixed Hodge
complexes on an algebraic variety or an analytic space, together with
a natural functor from the bounded derived category of mixed Hodge 
Modules to this category.

\bigskip
\noindent
{\bf 2.1.} Let
$ X $ be a separated complex analytic space, and
$ A $ a subfield of
$ {\Bbb R} $.
Let
$ {C}_{\hol}^{b}F({\Cal O}_{X},\Diff,A;W) $ denote the category whose objects are
$$
K = ((K_{F}; F,W), (K_{A},W), (K_{{\Bbb C}},W); \alpha _{F}, \alpha _{A}),
$$
where
$ (K_{F}; F,W) $ 
is a bounded filtered differential complex with a finite 
filtration
$ W $ such that
$ \DR_{X}^{-1} (K_{F};F,W) $ is holonomic in the sense of (1.3),
$ (K_{\Lmd},W) $ for
$ \Lmd = A, {\Bbb C} $ is a bounded complex of
$ \Lmd $-Modules with a finite increasing filtration
$ W $, and
$$
\alpha _{\Lmd} : (K_{\Lmd,{\Bbb C}},W) 
\rightarrow  (K_{{\Bbb C}},W)
\quad\text{for } \Lmd = F, A
$$
is a filtered quasi-isomorphism of complexes of
$ {\Bbb C} $-Modules.
Here
$ (K_{F,{\Bbb C}},W) $ is the underlying filtered complex of
$ {\Bbb C} $-Modules of
$ (K_{F}; F,W) $ and
$$
(K_{A,{\Bbb C}},W) = (K_{A},W)\otimes _{A}{\Bbb C}.
$$

A morphism
$ u : K \rightarrow  L $ of
$ {C}_{\hol}^{b}F({\Cal O}_{X},\Diff,A;W) $ consists of
$ (u_{F}, u_{A}, u_{{\Bbb C}}, u'_{F}, u'_{A}) $ where
$$
u_{F} : (K_{F}; F,W) \rightarrow  (L_{F}; F,W)
$$
is a morphism of filtered differential complex preserving
$ W $,
$$
u_{\Lmd} : (K_{\Lmd},W) \rightarrow 
 (L_{\Lmd},W)
\quad\text{for }\Lmd = A, {\Bbb C}
$$
is a morphism of filtered complexes of
$ \Lmd $-Modules, and
$$
u'_{\Lmd} : (K_{\Lmd,{\Bbb C}}^{i},W) 
\rightarrow  (L_{{\Bbb C}}^{i-1}, W)
\,\, (i \in {\Bbb Z})
$$
are morphisms of filtered
$ {\Bbb C} $-Modules for
$ \Lmd = F, A,  $ such that
$$
\alpha _{\Lmd}\circ u_{\Lmd,{\Bbb C}} - 
u_{{\Bbb C}}\circ \alpha _{\Lmd} = du'_{\Lmd} + 
u'_{\Lmd}d\quad  \text{for }\Lmd = F, A,
\tag 2.1.1
$$
where
$ u_{\Lmd,{\Bbb C}} : (K_{\Lmd,{\Bbb C}},W) \rightarrow  (L_{\Lmd,{\Bbb C}}, W) $ 
is induced by
$ u_{\Lmd} $.

A homotopy
$ h : K \rightarrow  L[-1] $ of
$ {C}_{\hol}^{b}F({\Cal O}_{X},\Diff,A;W) $ consists of
$ (h_{F}, h_{A}, h_{{\Bbb C}}, h'_{F}, h'_{A}) $,
where
$$
h_{F} : ({K}_{F}^{i}; F,W) \rightarrow  ({L}_{F}^{i-1}; F,W[-1]) \,\,\,
(i \in {\Bbb Z})
$$
are filtered differential morphisms preserving
$ W $, and
$$
\aligned
h_{\Lmd} : (K_{\Lmd}^{i},W) 
&\rightarrow  (L_{\Lmd}^{i-1},W[-1])
\,\, (i \in {\Bbb Z}) \\
h'_{\Lmd} : (K_{\Lmd,{\Bbb C}}^{i},W) 
&\rightarrow  (L_{{\Bbb C}}^{i-2}, W[-1])
\,\, (i \in {\Bbb Z})
\endaligned
$$
are morphisms of filtered
$ \Lmd $-Modules for
$ \Lmd = A, {\Bbb C},  $ such that they satisfy the following condition:

\smallskip
If we put
$ u_{\Lmd} = dh_{\Lmd} + h_{\Lmd}d $ for
$ \Lmd = F, A, {\Bbb C},  $ and
$$
u'_{\Lmd} = dh'_{\Lmd} - h'_{\Lmd}d + 
\alpha _{\Lmd}\circ h_{\Lmd,{\Bbb C}} - h
_{{\Bbb C}}\circ \alpha _{\Lmd}
\quad \text{for }\Lmd = F, A,
\tag 2.1.2
$$
then
$ u = (u_{F}, u_{A}, u_{{\Bbb C}}, u'_{F}, u'_{A}) : K \rightarrow  L $ is
a morphism of $ {C}_{\hol}^{b}F({\Cal O}_{X},\Diff,A;W) $.
(Here
$ (W[m])_{k} = W_{k-m} $ in general.)

\smallskip
The above
$ u $ is denoted by
$ dh + hd $.
A morphism
$ u $ of
$ {C}_{\hol}^{b}F({\Cal O}_{X},\Diff,A;W)  $ is called homotopic to zero if
$ u = dh + hd $ for some homotopy
$ h $.

For a morphism
$ u : K \rightarrow  L $ of
$ {C}_{\hol}^{b}F({\Cal O}_{X},\Diff,A;W)  $,
the mapping cone
$ M = C(K \rightarrow  L) $ consists of
$ M_{\Lmd} = C(K_{\Lmd} \rightarrow  L_{\Lmd}) $ for
$ \Lmd = F, A , {\Bbb C} $ with
$$
\aligned
F_{p}M_{F} &= C(F_{p}K_{F} \rightarrow  F_{p}L_{F}),
\\
W_{k}M_{\Lmd} &= C(W_{k-1}K_{\Lmd} \rightarrow  W_{k}L_{\Lmd})\,\, 
(\Lmd = F, A , {\Bbb C}),
\endaligned
\tag 2.1.3
$$
where
$ \alpha _{\Lmd} : C(K_{\Lmd,{\Bbb C}} \rightarrow  L_{\Lmd,{\Bbb C}}) 
\rightarrow  C(K_{{\Bbb C}} \rightarrow  L_{{\Bbb C}}) $ for
$ \Lmd = F, A $ is given by
$$
\alpha _{\Lmd}(x,y) = (\alpha _{\Lmd}x, u'_{\Lmd}x + \alpha _{\Lmd}y)\quad  
\text{for }x \in  {K}_{\Lmd,{\Bbb C}}^{i+1}, y \in  {L}_{\Lmd,{\Bbb C}}^{i}.
$$

We define
$ K[1] = C(K\to 0) $, and more generally
$ K[m] $ for
$ m \in  {\Bbb Z} $ to be
$$
((K_{F}[m]; F,W[m]), (K[m]_{A},W[m]), (K_{{\Bbb C}}[m],W[m]); \alpha 
_{F}, \alpha _{A}).
$$
Then, as in [21], we have natural morphisms
$ K \rightarrow  L \rightarrow  C(K \rightarrow  L) \rightarrow  K[1] $.

We get similarly
$ {C}_{\hol}^{b}F(X,{\Cal D},A;W)  $ replacing
$ (K_{F};F,W) $ by
$ (K_{D};F,W) $ a holonomic complex of filtered
$ {\Cal D} $-Modules on
$ X $ with a finite increasing filtration
$ W $, and
$ (K_{F,{\Bbb C}}, W) $ by
$ (K_{{\Cal D},{\Bbb C}}, W) := \DR_{X} (K_{\Cal D},W) $.
Here
$ \DR_{X} $ is defined by choosing a covering
$ \{ U_{i} \rightarrow  V_{i}\}_{i} $ of
$ X $ as in (1.3).

We have also the notions of homotopy and mapping cone.
A natural functor
$$
\DR_{X}^{-1} : {C}_{\hol}^{b}F({\Cal O}_{X},\Diff,A;W)  
\rightarrow  {C}_{\hol}^{b}F(X,{\Cal D},A;W) ,
\tag 2.1.4
$$
is obtained by using the filtered quasi-isomorphism
$ \DR_{X}{\DR}_{X}^{-1}(K_{F},W) \rightarrow  (K_{F,{\Bbb C}},W) $.

Forgetting the filtration
$ W $, we have similarly
$ {C}_{\hol}^{b}F({\Cal O}_{X},\Diff,A) $,
$ {C}_{\hol}^{b}F(X,{\Cal D},A) $ together with natural functors
$$
\Gr_{k}^{W} : {C}_{\hol}^{b}F({\Cal O}_{X},\Diff,A;W) \to
{C}_{\hol}^{b}F({\Cal O}_{X},\Diff,A).
$$

For
$ K = ((K_{{\Cal D}}; F,W), (K_{A},W), (K_{{\Bbb C}},W); \alpha _{{\Cal D}}, 
\alpha _{A}) \in  {C}_{\hol}^{b}F(X,{\Cal D},A;W)  $, let
$ {\Cal H}^{i}K $ denote
$$
({\Cal H}^{i}(K_{{\Cal D}}, F;W[i]), {}^{p} {\Cal H}^{i}(K_{A},W[i]);
\alpha _{A}^{-1} \circ  \alpha_{{\Cal D}})
\in MF_{\hol}(X,{\Cal D},A;W),
\tag 2.1.5
$$
where
$$
\alpha _{A}^{-1} \circ  \alpha_{{\Cal D}} : 
\DR_{X} {\Cal H}^{i} (K_{\Cal D},W[i]) \simeq 
{}^{p} {\Cal H}^{i} (K_{A,{\Bbb C}},W[i])
$$
is an isomorphism in
the category of filtered perverse sheaves
$ \Perv(X,{\Bbb C};W) $. 
(See [2] for the perverse cohomology functor
$ {}^{p}{\Cal H}^{i} $, and (1.4) for
$ MF_{\hol}(X,{\Cal D},A;W) $.)
Here
$ F $, $ W $ denote also the induced filtrations on
$ {\Cal H}^{i} K_{{\Cal D}} $,
$ {}^{p} {\Cal H}^{i}K_{A} $, etc., and they are defined separately
using the exact sequence
$  {\Cal H}^{i} F_{p}K_{{\Cal D}} \to  {\Cal H}^{i} K_{{\Cal D}}
\to  {\Cal H}^{i} (K_{{\Cal D}}/ F_{p}K_{{\Cal D}}) $, etc.

We have similarly
$ {\Cal H}^{i}K \in MF_{\hol}(X,{\Cal D},A) $ for 
$ K \in  {C}_{\hol}^{b}F(X,{\Cal D},A) $.

\bigskip
\noindent
{\it Remarks.} (i) The morphisms of 
$ {C}_{\hol}^{b}F({\Cal O}_{X},\Diff,A;W)  $,
$ {C}_{\hol}^{b}F(X,{\Cal D},A;W) $, etc. are stable by composition.
In fact,
$ w = vu $ is given by
$ w_{\Lmd} = v_{\Lmd}u_{\Lmd} $,
$ w'_{\Lmd} = v'_{\Lmd}u_{\Lmd,{\Bbb C}} + v_{\Bbb C}u'_{\Lmd} $.
Furthermore, if
$ u $ or
$ v $ is homotopic to zero, so is
$ w $.
For example, if
$ v = dh + hd $, then
$ w = dg + gd $ with
$ g_{\Lmd} = h_{\Lmd}u_{\Lmd} $,
$ g'_{\Lmd} = h'_{\Lmd}u_{\Lmd,{\Bbb C}} + h_{\Bbb C}u'_{\Lmd} $.

(ii) For
$ K \in  {C}_{\hol}^{b}F({\Cal O}_{X},\Diff,A;W)  $,
we define a complex
$ \tilde{K}_{\Lmd} $ for
$ \Lmd = F, A $ by
$$
{\tilde{K}}_{\Lmd}^{i} = {K}_{\Lmd,{\Bbb C}}^{i}\oplus {K}_{\Bbb C}^{i-1}\quad  
\text{with }d(x,y) = (dx, \alpha _{\Lmd}x - dy).
$$
For a morphism
$ u : K \rightarrow  L $,
let
$$
\tilde{u}_{\Lmd} : \tilde{K}_{\Lmd} \rightarrow  \tilde{L}_{\Lmd}\quad  
\text{with }
\tilde{u}_{\Lmd}(x,y) = (u_{\Lmd}x, u'_{\Lmd}x + u_{{\Bbb C}}y).
$$
Then the condition (2.1.1) corresponds to
$ d\tilde{u}_{\Lmd} = \tilde{u}_{\Lmd}d $.
For a homotopy
$ h : K \rightarrow  L[-1] $,
let
$$
\tilde{h}_{\Lmd} : \tilde{K}_{\Lmd} \rightarrow  \tilde{L}_{\Lmd}[-1]\quad 
\text{with }
\tilde{h}_{\Lmd}(x,y) = (h_{\Lmd}x, - h'_{\Lmd}x - h_{{\Bbb C}}y).
$$
Then (2.1.2) corresponds to
$ \tilde{u}_{\Lmd} = d\tilde{h}_{\Lmd} + \tilde{h}_{\Lmd}d $.

\bigskip
\noindent
{\bf 2.2.~Definition.}
We say that
$ K = ((K_{{\Cal D}}; F), (K_{A}), (K_{{\Bbb C}}); \alpha _{{\Cal D}}, 
\alpha _{A}) \in  {C}_{\hol}^{b}F(X,{\Cal D},A)  $
splits in the derived categories in a compatible way
if we have (noncanonical) isomorphisms 
$$
\aligned
(K_{\Cal D},F) &\simeq \bigoplus _{i} 
({\Cal H}^{i} (K_{\Cal D},F))[-i]\quad  \text{in }D^{b}F({\Cal D}_{X}), \\
K_{\Lmd} &\simeq \bigoplus _{i} ({}^{p}
{\Cal H}^{i}K_{\Lmd})[-i]\quad  \text{in }D^{b}(X,\Lmd)
\endaligned
\tag 2.2.1
$$
for
$ \Lmd = A, {\Bbb C} $ in a compatible way with
$ \alpha _{\Cal D}, \alpha _{A} $.

The category
$ C_{\Cal H}^{b} (X,A,n)_{\Cal D} $ of pure
$ A $-Hodge
$ {\Cal D} $-complexes of weight
$ n $ on
$ X $ is the full subcategory of
$  {C}_{\hol}^{b}F(X,{\Cal D},A) $
consisting of
$ K $ such that
$ K $ splits in the derived categories 
in a compatible way, and the
$ {\Cal H}^{i}K $ are polarizable 
$ A $-Hodge Modules of weight
$ i + n $ for any
$ i $.
The category
$ C_{\Cal H}^{b} (X,A)_{\Cal D} $ of mixed
$ A $-Hodge
$ {\Cal D} $-complexes on
$ X $ is the full subcategory of
$  {C}_{\hol}^{b}F(X,{\Cal D},A;W) $
consisting of
$ K $ such that
$ \Gr_{k}^{W} K \in C_{\Cal H}^{b} (X,A,k)_{\Cal D} $ for any
$ k $.
The category
$ C_{\Cal H}^{b} (X,A) $ of mixed
$ A $-Hodge complexes on
$ X $ is the full subcategory of
$  {C}_{\hol}^{b}F({\Cal O}_{X},\Diff,A;W) $
consisting of
$ K $ such that
$ \DR_{X}^{-1} K \in C_{\Cal H}^{b} (X,A)_{\Cal D} $.
Then considering morphisms of
$ C_{\Cal H}^{b} (X,A,n)_{\Cal D} $,
$ C_{\Cal H}^{b} (X,A)_{\Cal D} $,
$ C_{\Cal H}^{b} (X,A) $ 
up to homotopy, we get respectively
$ K_{\Cal H}^{b} (X,A,n)_{\Cal D} $,
$ K_{\Cal H}^{b} (X,A)_{\Cal D} $,
$ K_{\Cal H}^{b} (X,A) $.

If
$ X $ is a complex algebraic variety, the categories
$ {\Cal K}_{\Cal H}^{b}(X,A) $,
$ {\Cal K}_{\Cal H}^{b}(X,A)_{\Cal D} $, etc. 
are defined in the same way as above, where
$ (K_{F}; F,W) $ is a bounded holonomic filtered (algebraic)
differential complex on
$ X $ with a finite increasing filtration
$ W $,
$ (K_{{\Cal D}}; F,W) $ is a holonomic filtered (algebraic)
$ {\Cal D} $-Module on
$ X $ with a finite increasing filtration
$ W $,
$ (K_{A},W), (K_{{\Bbb C}},W) $ are filtered complexes on
$ X^{\an} $,
 $ (K_{F,{\Bbb C}}, W) $ is the underlying filtered complex of
$ {\Bbb C} $-Modules of
$ (K_{F},W)^{\an} $,
and
$ (K_{{\Cal D},{\Bbb C}}, W) = \DR_{X}(K_{{\Cal D}},W) $.

If
$ X $ is an algebraic variety defined over a subfield
$ k $ of
$ {\Bbb C} $,
the categories
$ {\Cal K}_{\Cal H}^{b}(X,A) $,
$ {\Cal K}_{\Cal H}^{b}(X,A)_{\Cal D} $, etc. 
are defined similarly, 
by assuming
$ (K_{F}; F,W), (K_{{\Cal D}}; F,W) $ are defined over
$ k $.

An object
$ K $ of
$ K_{\Cal H}^{b} (X,A)  $ or
$ K_{\Cal H}^{b} (X,A)_{\Cal D} $ is weakly acyclic if
$ K_{A} $ is acyclic, and a morphism
$ u : K \rightarrow  L $ of
$ K_{\Cal H}^{b} (X,A)  $ or
$ K_{\Cal H}^{b} (X,A)_{\Cal D} $ is  a weak quasi-isomorphism if
$ u_{A} : K_{A} \rightarrow  L_{A} $ is a quasi-isomorphism forgetting
$ W $.
(A weak quasi-isomorphism is not necessarily 
a filtered quasi-isomorphism for
$ W $.)

\bigskip
\noindent
{\it Remarks.}  (i) The decomposition (2.2.1) implies that
$ F $ on
$ K_{\Cal D} $ is strict.

If we define
$ {\Cal D}_{\Cal H}^{b}(X,A,n)_{\Cal D} $ by inverting filtered
quasi-isomorphisms (or, equivalently, weak quasi-isomorphisms) of
$ {\Cal K}_{\Cal H}^{b}(X,A,n)_{\Cal D} $, 
we see that the decomposition (2.2.1) is equivalent to 
(the existence of) a noncanonical isomorphism
$$
K \simeq \bigoplus_{i}{\Cal H}^{i}K[-i] \quad\text{in }
{\Cal D}_{\Cal H}^{b}(X,A,n)_{\Cal D},
$$
because the uniqueness of the isomorphism is not in question.
Here 
$ {\Cal H}^{i}K = ({\Cal H}^{i}K_{\Cal D} $,
$ {}^{p}{\Cal H}^{i}K_{A} $,
$ {}^{p}{\Cal H}^{i}K_{\Bbb C}) $ can be defined in
$ {\Cal D}_{\Cal H}^{b}(X,A,n)_{\Cal D} $ by using (2.7) together with 
the forgetful functor which forgets the filtration
$ W $.
In fact, if we define
$ K_{\hol}^{b}F(X,{\Cal D},A) $ by considering morphisms of
$ C_{\hol}^{b}F(X,{\Cal D},A) $ up to homotopy, and then
$ D_{\hol}^{b}F(X,{\Cal D},A) $ by inverting filtered quasi-isomorphisms,
we see that a natural functor
$ {\Cal D}_{\Cal H}^{b}(X,A,n)_{\Cal D} \to D_{\hol}^{b}F(X,{\Cal D},A) $
is fully faithful.
(Note however that we do not have well-defined functors
$ \Gr_{k}^{W} : {\Cal D}_{\Cal H}^{b}(X,A)_{\Cal D}
\to {\Cal D}_{\Cal H}^{b}(X,A,n)_{\Cal D} $.)

More generally, let
$ K, L \in C_{\hol}^{b}F(X,{\Cal D},A) $, and assume that we have
isomorphisms
$ (K_{\Cal D},F) \simeq (L_{\Cal D},F) $ in
$ D^{b}F(X,{\Cal D}) $ and
$ K_{\Lmd} \simeq L_{\Lmd} $ in
$ D^{b}(X,{\Lmd}) \, (\Lmd = A, {\Bbb C}) $ 
compatible with 
$ \alpha_{\Cal D} $, 
$ \alpha_{A} $.
Then we have
$ (M_{\Cal D},F) \in C^{b}F(X,{\Cal D}) $,
$ M_{\Lmd} \in C^{b}(X,{\Lmd}) \,  (\Lmd = A, {\Bbb C}) $
with (filtered) quasi-isomorphisms
$ (K_{\Cal D},F) \to (M_{\Cal D},F) \leftarrow  (L_{\Cal D},F) $ in
$ C^{b}F(X,{\Cal D}) $,
$ K_{\Lmd} \to M_{\Lmd} \leftarrow  L_{\Lmd} $ in
$ C^{b}(X,{\Lmd}) $, which represent the given isomorphisms in the
derived categories.
Then replacing 
$ M_{\Bbb C} $ if necessary, we have quasi-isomorphisms
$ \alpha_{\Lmd} : M_{\Lmd,{\Bbb C}} \to M_{\Bbb C} $ for
$ \Lmd = {\Cal D}, A $,
which give a commutative diagram in
$ K^{b}(X,{\Bbb C}) $.
So we get filtered quasi-isomorphisms
$ K \to M \leftarrow L $ in
$ C_{\hol}^{b}F(X,{\Cal D},A) $.

(ii) If 
$ X $ is algebraic, we have  natural functors
$$
{K}_{\Cal H}^{b}(X,A) \to
{K}_{\Cal H}^{b}(X^{\an},A),\,\,\,
{K}_{\Cal H}^{b}(X,A)_{{\Cal D}} \to
{K}_{\Cal H}^{b}(X^{\an},A)_{{\Cal D}}.
$$

(iii) In the definition of morphisms of
$ {\Cal C}_{\Cal H}^{b}(X,A) $,
$ {\Cal C}_{\Cal H}^{b}(X,A)_{\Cal D}$, etc., we need the homotopy
$ u'_{\Lmd} $ (which was not used in [1] where
$ X = pt) $ in order to construct a natural functor from the derived category 
of bounded complexes of mixed Hodge Modules to
$ {D}_{\Cal H}^{b}(X,A)_{\Cal D} $.
See (2.7) below.
The decomposition (2.2.1) of
$ {\Gr}_{k}^{W}K $ is necessary to prove the stability by the direct image 
under a proper morphism (see (2.8) below), and has nothing to do with
$ p $-structure in [1] when
$ X = pt $.
In fact, the latter was used in loc.~cit.
to identify the bounded derived category
of polarizable mixed Hodge structures with a full subcategory of
the bounded derived category of mixed Hodge structures.
In the case
$ X = pt $, note that (2.2.1) is not necessary in the definition of
$ {D}_{\Cal H}^{b}(X,A) $, because it is not used in 
(2.10) below so that
$ {D}_{\Cal H}^{b}(pt,A) $ is equivalent to the category
which is defined without assuming (2.2.1).

\bigskip\noindent
{\bf 2.3.~Proposition.} {\sl
For
$ K = ((K_{{\Cal D}}; F,W), (K_{A},W), (K_{{\Bbb C}},W); \alpha _{{\Cal D}}, 
\alpha _{A}) \in  {C}_{\Cal H}^{b}(X,A)_{{\Cal D}},  $ we have the
weight spectral  sequence
$$
{E}_{1}^{p,q} = {\Cal H}^{p+q}{\Gr}_{-p}^{W}K \Rightarrow  
{\Cal H}^{p+q}K\quad  \text{in }\MHW(X,A),
\tag 2.3.1
$$
which degenerate at
$ E_{2} $,
and whose
$ d_{1} $ is strictly compatible with
$ F $.
(See (1.4) for
$ \MHW(X,A) .) $
Furthermore,
$ (K_{{\Cal D}}; F, \Dec W) $ is strict in the sense of [14, 1.2] 
(i.e., bistrict, see also (A.2) below), 
and gives the Hodge and weight filtrations on
$ {\Cal H}^{i}K $.
Consequently, a weak quasi-isomorphism
$  K \rightarrow  L $ of
$ {C}_{\Cal H}^{b}(X,A)_{{\Cal D}}) $ induces a bifiltered quasi-isomorphism
$$
(K_{{\Cal D}}; F, \Dec W) \rightarrow  (L_{{\Cal D}}; F, \Dec W).
\tag 2.3.2
$$
}

\demo\nofrills {Proof.\usualspace}
This follows from (1.4) and (A.3).
\enddemo

\medskip
\noindent
{\it Remarks.}  (i) In the case
$ X = pt $, (2.3) follows from [4].

(ii) The last assertion of (2.3) means
$$
(K_{{\Cal D}}; F, \Dec W)
\text{ is well-defined in }
 D_{\hol}^{b}F(X,{\Cal D};W')
\text{ for }
K \in  {D}_{\Cal H}^{b}(X,A)_{{\Cal D}}.
\tag 2.3.3
$$
This justifies the definition of
$ {D}_{\Cal H}^{b}(X,A) $,
$ {D}_{\Cal H}^{b}(X,A)_{{\Cal D}} $ by inverting the weak 
quasi-isomorphisms.

The well-definedness of
$ \Dec W $ in the exact category implies for example that
$ \Dec $ commutes with
$ {\Gr}_{p}^{F} $.
See (A.3).

\bigskip\noindent
{\bf 2.4.~Proposition-Definition.} {\sl The categories
$ {K}_{\Cal H}^{b}(X,A) $ and
$ {K}_{\Cal H}^{b}(X,A)_{{\Cal D}} $ are triangulated categories whose 
distinguished triangles are isomorphic to the images in
$ {K}_{\Cal H}^{b}(X,A) $, 
$ {K}_{\Cal H}^{b}(X,A)_{{\Cal D}} $ of
$$
\rightarrow  K \rightarrow  L \rightarrow  C(K \rightarrow  L) 
\overset +1 \to \rightarrow 
$$
given in (2.1).
The full subcategories of
$ {K}_{\Cal H}^{b}(X,A) $,
${K}_{\Cal H}^{b}(X,A)_{{\Cal D}} $ consisting of
weakly  acyclic objects are thick subcategories in the sense of [21], and the 
associated quotient categories, which will be denoted by
$ {D}_{\Cal H}^{b}(X,A) $,
$ {D}_{\Cal H}^{b}(X,A)_{{\Cal D}} $ are also 
triangulated categories whose distinguished triangles are as above.
The functor
$ {\DR}_{X}^{-1} $ in (2.1) induces a natural functor
$$
{\DR}_{X}^{-1} : {D}_{\Cal H}^{b}(X,A) \rightarrow  {D}
_{\Cal H}^{b}(X,A)_{{\Cal D}}.
$$
}

\demo\nofrills {Proof.\usualspace}
We prove the first assertion for
$ {K}_{\Cal H}^{b}(X,A), {D}_{\Cal H}^{b}(X,A) $.
The argument is similar for
$ {K}_{\Cal H}^{b}(X,A)_{{\Cal D}}, {D}_{\Cal H}^{b}(X,A)_{{\Cal D}} $.

For a morphism
$ u : K \rightarrow  L $ in
$ {C}_{\Cal H}^{b}(X,A) $,
its mapping cone
$ M = C(K \rightarrow  L) $ belongs to
$ {C}_{\Cal H}^{b}(X,A) $, because
$ {\Gr}_{k}^{W}M = ({\Gr}_{k-1}^{W}K)[1] \oplus {\Gr}_{k}^{W}L $
(since
$ u_{\Lmd} $ and
$ u'_{\Lmd} $ preserve the filtration
$ W $).
Then we can verify the axioms of triangulated category in [21] for  
$ {K}_{\Cal H}^{b}(X,A) $.

We see first that compositions of successive two morphisms of the
triangle are homotopic to zero.
For the composition of
$ K \to L \to C(K\to L) $, we have a homotopy
$ h $ defined by
$ h_{\Lmd}(x) = (x,0) $ and
$ h'_{\Lmd} = 0 $.
We verify similarly
$ C(id : K \to K) = 0 $ (i.e., 
$ id : C(id) \to C(id) $ is homotopic to zero), using a homotopy defined by
$ h_{\Lmd}(x,y) = (y,0) $ and
$ h'_{\Lmd} = 0 $.

For the invariance of triangles by shift, we can show for example that
the projection
$ C(u + id : K\oplus L \to L) \to K[1] $ is an isomorphism whose inverse 
$ v $ is
given by
$ v_{\Lmd}(x) = (x, -u_{\Lmd}(x),0) $,
$ v'_{\Lmd}(x) = (0, u'_{\Lmd}(x),0) $,
using the homotopy
$ h : C(u + id) \to C(u + id) $ defined by
$ h_{\Lmd}(x,y,z) = (0,z,0) $,
$ h'_{\Lmd} = 0 $.

If we have morphisms
$$
u : K \rightarrow  L,\,\,\,v : M \rightarrow  N,\,\,\,f : K \rightarrow  
M,\,\,\,g : L \rightarrow  N\quad  \text{in }{C}_{\Cal H}^{b}(X,A)
$$
and a homotopy
$ h : K \rightarrow  N[-1] $ such that
$ gu - vf = dh +hd $,
then we have
$ w : C(K \rightarrow  L) \rightarrow  C(M \rightarrow  N) $ by
$$
\aligned
w_{\Lmd}(x,y) &= (f_{\Lmd}x, h_{\Lmd}x + g_{\Lmd}y)\quad
\text{for }\Lmd = F, A, {\Bbb C}, 
\\
w'_{\Lmd}(x,y) &= (-f'_{\Lmd}x, -h'_{\Lmd}x + g'_{\Lmd}y)\quad
\text{for }\Lmd = F, A.
\endaligned
$$
So the mapping cone is unique up to a noncanonical isomorphism
due to [21, 1.2].
Then we verify the octahedral axiom by showing the isomorphism
$$
C(C(K \rightarrow  L) \rightarrow  C(M \rightarrow  N)) = C(C(K 
\rightarrow  M) \rightarrow  C(L \rightarrow  N)).
$$

Then it is clear that the subcategory consisting of weakly acyclic 
objects is thick,  and the remaining assertion follows from [21].
\enddemo

\bigskip\noindent
{\bf 2.5.} We are now going to construct a natural functor
from the bounded derived category of mixed Hodge Modules
$ D^{b}\MHM(X,A) $
to the triangulated category of mixed Hodge 
$ {\Cal D} $-complexes
$ D_{\Cal H}^{b}(X,A)_{\Cal D} $.
We first construct a category
$ C_{\MHM}^{b}(X,A) $ which is naturally equivalent to
the category of bounded complexes of mixed Hodge Modules 
$ C^{b}\MHM(X,A) $.

For
$ \Lmd = A, {\Bbb C} $,
let
$ C_{c}^{*}(X,\Lmd;G,W') $ be the category of bifiltered complexes of
$ \Lmd $-Modules 
$ (K_{\Lmd};G,W') $ on
$ X $ such that
$ G $,
$ W' $ are finite filtrations on
$ K_{\Lmd} $,
$ {\Cal H}^{i}{\Gr}_{G}^{j}{\Gr}_{W'}^{k}K_{\Lmd} $ are 
$ \Lmd $-constructible sheaves on
$ X $, and
$ K_{\Lmd} $ satisfy the boundedness condition corresponding to
$ * = +,-,b $.
Then
$ K_{c}^{*}(X,\Lmd;G,W') $,
$ D_{c}^{*}(X,\Lmd;G,W') $ are defined as in [21].
Let
$ K_{c}^{+,b}(X,\Lmd;G,W')_{\inj} $ be the full subcategory of
$ K_{c}^{+}(X,\Lmd;G,W') $ consisting of 
$ (K_{\Lmd};G,W') $ such that
$ {\Gr}_{G}^{j}{\Gr}_{W'}^{k}K_{\Lmd}^{i} $ are injective
$ \Lmd $-Modules and
$ {\Cal H}^{i}{\Gr}_{G}^{j}{\Gr}_{W'}^{k}K_{\Lmd} = 0 $ for
$ i \gg 0 $.
Then we have an equivalence of categories
$$
K_{c}^{+,b}(X,\Lmd;G,W')_{\inj} = D_{c}^{b}(X,\Lmd;G,W').
\tag 2.5.1
$$
 See [2] (and also (A.7) below).

Let
$ MF_{\hol}(X,{\Cal D};W') $ denote the category of filtered
$ {\Cal D} $-Modules 
$ (M,F) $ with a finite increasing filtration
$ W' $ on
$ X $ such that the
$ \Gr_{k}^{W'}(M,F) $ are holonomic.
We consider the category
$ {K}_{\hol}^{+,b}(X,{\Cal D},A;G,W') $ whose objects are
$$
K = ((K_{{\Cal D}}; F,W'), (K_{A}; G,W'), (K_{{\Bbb C}}; G,W'); \alpha 
_{{\Cal D}}, \alpha _{A})
$$
where
$ (K_{{\Cal D}}; F,W') $ is a bounded complex of
$ MF_{\hol}(X,{\Cal D};W') $,
$ (K_{\Lmd}; G,W') $ is an object of
$ K_{c}^{+,b}(X,\Lmd;G,W')_{\inj} $  for
$ \Lmd = A, {\Bbb C} $,
and
$$
\alpha _{\Lmd} : (K_{\Lmd,{\Bbb C}}; G,W') \rightarrow  (K_{{\Bbb C}}; G,W')
$$
is a bifiltered quasi-isomorphism in
$ K^{+}(X,{\Bbb C};G,W') $ for
$ \Lmd = {\Cal D}, A $.
Here
$$
(K_{{\Cal D},{\Bbb C}}; G,W') = \DR_{X}(K_{{\Cal D}}; \sigma ,W'),\quad 
(K_{A,{\Bbb C}}; G,W') = (K_{A}; G,W')\otimes _{A}{\Bbb C},
$$
and
$ \sigma  $ is the filtration ``b\^ete'' in [4].

A morphism
$ \tilde{u} : K \rightarrow  L $ of 
$ K_{\hol}^{+,b}F(X,{\Cal D},A;G,W') $
consists of
$ (\tilde{u}_{{\Cal D}}, \tilde{u}_{A}, \tilde{u}_{{\Bbb C}}) $ which are 
morphisms of
$ C^{b}F(X,{\Cal D};W') $, 
$ K_{c}^{+}(X,A;G,W') $,
$ K_{c}^{+}(X,{\Bbb C};G,W') $ respectively such that
$$
\alpha _{\Lmd}\circ \tilde{u}_{\Lmd,{\Bbb C}} = 
\tilde{u}_{{\Bbb C}}\circ \alpha _{\Lmd}\quad  \text{in }
K_{c}^{+}(X,{\Bbb C};G,W')
$$
for
$ \Lmd = {\Cal D}, A $,
where
$ \tilde{u}_{\Lmd,{\Bbb C}} : (K_{\Lmd,{\Bbb C}}; G,W') \rightarrow 
 (L_{\Lmd,{\Bbb C}}; G,W') $ denotes the induced morphism by
$ \tilde{u}_{\Lmd} $.
(Here
$ \tilde{u}_{{\Cal D}} $ is defined in
$ C_{\hol}^{b}F(X,{\Cal D};G,W') $ because the homotopy must preserve the
filtration $ \sigma  $.)

We define
$ C^{b}_{\MHM}(X,A) $ to be the full subcategory of
$ K_{\hol}^{+,b}F(X,{\Cal D},A;G,W') $ consisting of
$ K $ such that
$$ 
(({K}_{\Cal D}^{j}; F,W'), ({\Gr}_{G}^{j}K_{A}[j],W'), 
\alpha _{A}^{-1} \circ  \alpha_{{\Cal D}}) \in
\MHM(X,A)
$$
where
$ \alpha _{A}^{-1} \circ  \alpha_{{\Cal D}} : 
\DR_{X}({K}_{\Cal D}^{j},W') \simeq
({\Gr}_{G}^{j}K_{A}[j],W')\otimes _{A}{\Bbb C} $
is an isomorphism in the category of filtered perverse sheaves
$ \Perv(X, {\Bbb C};W') $.

\bigskip\noindent
{\it Remark.}
For
$ \Lmd = A, {\Bbb C} $, let
$ D_{c}^{b}(X,\Lmd;G_{\text{b\^ete}},W') $ denote the full subcategory of
$ D_{c}^{b}(X,\Lmd;G,W') $ whose objects satisfy the condition
$$
{\Gr}_{G}^{j}{\Gr}_{W'}^{k}K[j] \in  \Perv(X,\Lmd)
\quad\text{for any } j,k.
$$
Then we have a filtration
$ W' $ on
$ \overline{K}^{i} := {\Gr}_{G}^{i}K[i] \in  {\Cal C}  $ such that
$$
(W^{\prime a}/W^{\prime b})\overline{K}^{i} = 
{\Gr}_{G}^{i}(W^{\prime a}/W^{\prime b})K[i]\quad  \text{in } \Perv(X,\Lmd).
\tag 2.5.2
$$
So we get a functor
$$
{\Gr}_{G} : D_{c}^{b}(X,\Lmd;G_{\text{b\^ete}},W')
\rightarrow  C^{b}\Perv(X,\Lmd;W'),
\tag 2.5.3
$$
where the right-hand side is the category of bounded complexes of
filtered perverse sheaves.
See (A.4) below.
Furthermore, (2.5.3) is an equivalence of categories.
See (A.5) below.
For the full faithfulness of (2.5.3),
we use the spectral sequences
$$
\aligned
{E}_{1}^{p,q} = \prod _{j-i=p} {\Ext}_{D^{+}(X,\Lmd;W')}^{q}
&({\Gr}_{G}^{i}(K,W')[i], {\Gr}_{G}^{j}(L,W')[j]) \\
\Rightarrow  \,
&{\Ext}_{D^{+}(X,\Lmd;W')}^{p+q}((K,W'), (L,W')), \\
{E}_{1}^{p,q} =  \text{(same as above if }
&p \ge  k,\,\,\, \text{and }0\,\,\, \text{otherwise)  } \\
\Rightarrow  \,
& {\Ext}_{D^{+}(X,\Lmd;G,W')}^{p+q}((K; G,W'), (L; G[k],W')),
\endaligned
\tag 2.5.4
$$
such that  $ {E}_{1}^{p,q} = 0 $ for
$ q < 0 $.
See (A.5.3).
Note that the edge morphism to
$ E_{1}^{k,q} $ of the second spectral sequence is induced by the 
passage to
$ {\Gr}_{G} $ (i.e. by the functor (2.5.3)).

\bigskip\noindent
{\bf 2.6.~Proposition.} {\sl
We have naturally an equivalence of categories
$$
{C}_{\MHM}^{b}(X,A) = C^{b}\MHM(X,A) 
$$
where the right-hand side is the category of bounded complexes of mixed
$ A $-Hodge Modules.
}

\demo\nofrills {Proof.\usualspace}
Since (2.5.3) is an equivalence of categories by (A.5), we see 
(using (2.5.1)) that
$ K_{\hol}^{+,b}F(X,{\Cal D},A;G,W') $ is equivalent to the fiber product of
$ C^{b}MF_{\hol}(X,{\Cal D};W') $ and
$ C^{b}\Perv(X,A;W') $ over
$ C^{b}\Perv(X,{\Bbb C};W') $.
This is the category consisting of pairs
$$
(M,F;W) \in C^{b}MH_{\hol}(X,{\Cal D};W'), \quad
(\overline{K}_{A},W') \in C^{b}\Perv(X,A;W') 
$$
with a given isomorphism
$$
\alpha : \DR_{X}(M,W') \simeq (\overline{K}_{A},W')\otimes_{A}{\Bbb C} 
\quad\text{in }
C^{b}\Perv(X,{\Bbb C};W').
$$
So the assertion is clear by definition of
$ {C}_{\MHM}^{b}(X,A) $.
\enddemo

\bigskip\noindent
{\bf 2.7.~Theorem.} {\sl We have a natural functor
$$
\varepsilon  : D^{b}\MHM(X,A) \rightarrow  {D}_{\Cal H}^{b}(X,A)_{{\Cal D}},
\tag 2.7.1
$$
where the left-hand side is the derived category of bounded complexes of 
mixed
$ A $-Hodge Modules.
}
\bigskip

\demo\nofrills {Proof.\usualspace}
 We prove the assertion in the case
$ X $ is an analytic space.
The argument is similar in the case
$ X $ is algebraic.

Using the truncation
$ \tau _{\le i} $ for
$ i $ large enough, we can replace the bounded condition for
$ (K_{\Lmd}, W) \,\, (\Lmd = A, {\Bbb C}) $
in the definition of
$ {K}_{\Cal H}^{b}(X,A)_{{\Cal D}} $ by the cohomological one (i.e.,
$ {\Cal H}^{j}{\Gr}_{k}^{W}K_{\Lmd} = 0 $ for
$ \Lmd = A, {\Bbb C} $ if
$ j $ is large enough).
Then we have naturally a well-defined functor
$$
{C}_{\MHM}^{b}(X,A) \rightarrow  {K}_{\Cal H}^{b}(X,A)_{{\Cal D}}.
\tag 2.7.2
$$

In fact, we may change the above definition of
$ K_{\hol}^{+,b}F(X,{\Cal D},A;G,W') $
so that
$$
(K_{\Lmd}; G,W') \in  C_{c}^{+,b}(X,\Lmd;G,W')_{\inj}
\quad  \text{for }\Lmd = A, {\Bbb C},
$$
 $ \alpha _{\Lmd} $ is a morphism of
$ C^{+}(X,{\Bbb C};G,W') $ for
$ \Lmd = {\Cal D}, A $,
and
$ \tilde{u}_{\Lmd} $ is represented by a morphism
$ u_{\Lmd} $ of 
$ C^{+}(X,\Lmd;G,W') $ for
$ \Lmd = A, {\Bbb C} $.
So
$ \alpha _{\Lmd}\circ u_{\Lmd,{\Bbb C}} = 
u_{{\Bbb C}}\circ \alpha _{\Lmd} $ holds up to 
homotopy in
$ C^{+}(X,{\Bbb C};G,W') $,
and we consider morphisms
$ {u} = (u_{\Cal D}, u_{A}, u_{\Bbb C}) $ up to homotopy 
in the sense that 
$ {u} $ and
$ {v} $ are homotopic if 
$ u_{\Lmd} $ and
$ v_{\Lmd} $ are homotopic in
$ C^{+}(X,\Lmd;G,W') $ for
$ \Lmd = A, {\Bbb C} $ and
$ u_{{\Cal D}} = v_{{\Cal D}} $.
(Note that these homotopies are not given,
but the homotopies
$ u'_{\Lmd} $ in the definition of
$ C_{\Cal H}^{b}(X,A)_{\Cal D} $ are given.
For the well-definedness of (2.7.2), the difference must be absorbed by
the homotopy
$ h $.)

Let
$ W $ be the convolution of
$ G $ and
$ W' $,
i.e.,
$ W_{k} = \sum _{j} G^{j}W'_{k+j} $,
where
$ G = \sigma  $ on
$ K_{{\Cal D}} $.
Then
$$
((K_{{\Cal D}}; F,W), (K_{A},W), (K
_{{\Bbb C}},W); \alpha _{{\Cal D}}, \alpha _{A})
$$
belongs to
$ {C}_{\Cal H}^{b}(X,A)_{{\Cal D}} $,
because the decomposition (2.2.1) is satisfied by
$$
{\Gr}_{k}^{W} = \bigoplus _{j} {\Gr}_{G}^{j}{\Gr}_{k+j}^{W'}.
\tag 2.7.3
$$
Here
$ {\Cal H}^{i}{\Gr}_{G}^{j}{\Gr}_{k}^{W'}K_{\Lmd} = 0 $ and
$ {\Cal H}^{i}{\Gr}_{k}^{W}K_{\Lmd} = 0 $ for
$ i \gg 0 $ with
$ \Lmd = A, {\Bbb C} $ by assumption.
Furthermore, we have 
$$
\alpha _{\Lmd}\circ u_{\Lmd,{\Bbb C}} - 
u_{{\Bbb C}}\circ \alpha _{\Lmd} = du'_{\Lmd} + 
u'_{\Lmd}d\quad  \text{with }u'_{\Lmd} \in  G^{0}W'_{0}\Hom^{-1}(K_{\Lmd,
{\Bbb C}}, L_{{\Bbb C}})
$$
for
$ \Lmd = {\Cal D}, A $,
and they give a morphism
$ u $ in
$ {C}_{\Cal H}^{b}(X,A)_{{\Cal D}} $.
(In fact, an element of
$ G^{0}W'_{0}\Hom^{-1}(K_{\Lmd,{\Bbb C}}, L_{{\Bbb C}}) $ preserves
$ G, W' $ and hence
$ W $.)
Here
$ \Hom^{\ssbull} (K_{\Lmd,{\Bbb C}}, L_{{\Bbb C}}) $ denotes a complex
as in (A.5).
We have to show that
$ u $ is unique up to homotopy.

Assume we have
$ v_{\Lmd} $ which is homotopic to
$ u_{\Lmd} $ for
$ \Lmd = A, {\Bbb C} $.
Put
$ v_{{\Cal D}} = u_{{\Cal D}} $.
We define
$ v'_{\Lmd} $ as above for
$ \Lmd = {\Cal D}, A $ so that we get a morphism
$ v $ in
$ {C}_{\Cal H}^{b}(X,A)_{{\Cal D}} $.
Then
$$
u_{\Lmd} - v_{\Lmd} = dh_{\Lmd} + h_{\Lmd}d\quad  \text{with }h_{\Lmd} \in  
G^{0}W'_{0}\Hom^{-1}(K_{\Lmd}, L_{\Lmd})
$$
for
$ \Lmd = A, {\Bbb C} $ by assumption.
Put
$ h_{{\Cal D}} = 0 $.
We consider
$$
w_{\Lmd} = u'_{\Lmd} - v'_{\Lmd} - \alpha _{\Lmd}\circ h_{\Lmd,{\Bbb C}} + 
h_{{\Bbb C}}\circ \alpha _{\Lmd} \in  
G^{0}W'_{0}\Hom^{-1}(K_{\Lmd,{\Bbb C}}, L_{{\Bbb C}})
$$
for  $ \Lmd = {\Cal D}, A $.
Then
$ dw_{\Lmd} + w_{\Lmd}d = 0 $,
and
$$
w_{\Lmd} = dh'_{\Lmd} - h'_{\Lmd}d\quad  \text{with }h'_{\Lmd} \in  
G^{0}W'_{0}\Hom^{-2}(K_{\Lmd,{\Bbb C}}, L_{{\Bbb C}}),
$$
because
$ H^{-1}G^{0}W'_{0}\Hom^{\ssbull }(K_{\Lmd,{\Bbb C}}, L_{{\Bbb C}}) = 0 $ by 
the second spectral sequence of (2.5.4) with
$ k = 0 $.
So
$ u - v $ is homotopic to zero in
$ {C}_{\Cal H}^{b}(X,A)_{{\Cal D}} $,
and (2.7.2) is well-defined.

Now we have to show that (2.7.2) induces
$$
K^{b}\MHM(X,A) \rightarrow  {K}_{\Cal H}^{b}(X,A)_{{\Cal D}},
\tag 2.7.4
$$
i.e.,
$ u $ is homotopic to zero in
$ {C}_{\Cal H}^{b}(X,A)_{{\Cal D}} $,
if it corresponds to a morphism of
$ C^{b}\MHM(X,A) $ which is homotopic to zero.
Let
$ (\overline{K}_{\Lmd},W') $ be the image of
$ (K_{\Lmd};G,W') $ by (2.5.3),
i.e.,
$$
(W'_{a}/W'_{b}){\overline{K}}_{\Lmd}^{i} = {\Gr}_{G}^{i}(W'_{a}/W'_{b})K_{\Lmd}[i]
\quad \text{in } \Perv(X,\Lmd).
$$
Put
$ (\overline{K}_{{\Cal D}};F,W') = (K_{{\Cal D}};F,W') $.
(Similarly for
$ L $.)
Let
$$
\aligned
\overline{h}_{{\Cal D}} : (\overline{K}_{{\Cal D}}^{i};F,W') 
&\rightarrow  
(\overline{L}_{{\Cal D}}^{i-1};F,W'),
\\
\overline{h}_{\Lmd }: 
(\overline{K}_{\Lmd}^{i},W') &\rightarrow  (\overline{L}_{\Lmd}^{i-1},W')
\,\, (\Lmd = A, {\Bbb C})
\endaligned
$$
be homotopies which are compatible with
$ \Gr_{G}\alpha _{\Lmd} $ for
$ \Lmd = {\Cal D}, A $,
i.e.,
$$
\Gr_{G}\alpha _{\Lmd}{\circ }\overline{h}_{\Lmd,{\Bbb C}} = 
\overline{h}_{{\Bbb C}}{\circ }\Gr_{G}\alpha _{\Lmd}
\quad \text{in } \Perv(X,{\Bbb C}),
$$
where
$ \overline{h}_{{\Cal D},{\Bbb C}} = \DR_{X}(\overline{h}_{{\Cal D}}) $ and
$ \overline{h}_{A,{\Bbb C}} = \overline{h}_{A}\otimes _{A}{\Bbb C} $.
(This condition is satisfied because the homotopies come from a homotopy 
of mixed Hodge Modules.)

Let
$ \overline{u}_{\Lmd} = d\overline{h}_{\Lmd} + \overline{h}_{\Lmd}d $ for
$ \Lmd = {\Cal D}, A, {\Bbb C} $.
Then by the second spectral sequence of (2.5.4) with
$ k = 0 $, there exist morphisms in
$ C^{+}(X,\Lmd;G,W') $ :
$$
u_{\Lmd} : (K_{\Lmd};G,W') \rightarrow  (L_{\Lmd};G,W')\quad  \text{for }\Lmd = A, 
{\Bbb C},
$$
such that
$ \overline{u}_{\Lmd} $ is the image of
$ u_{\Lmd} $ by (2.5.3) (i.e.,
$ \overline{u}_{\Lmd} = \Gr_{G}u_{\Lmd}) $ for
$ \Lmd = A, {\Bbb C} $.
Put
$ u_{{\Cal D}} = \overline{u}_{{\Cal D}} $.

The image of
$ u_{\Lmd} $ in
$ H^{0}G^{-1}W'_{0}\Hom^{\ssbull }(K_{\Lmd},L_{\Lmd}) $ is zero for
$ \Lmd = A, {\Bbb C} $  due  to the spectral sequence (with
$ k = -1 $), because 
$ u_{\Lmd} $ in
$ H^{o}G^{0}W'_{0} \Hom(K_{\Lmd},L_{\Lmd}) $ corresponds by the
edge morphism (i.e. by 
$ \Gr_{G} $ or (1.5.3)) to
$ \overline{u}_{\Lmd} \in E_{1}^{0,0} $ of the second spectral sequence
of (1.5.4) with
$ k = 0 $, and
$ \overline{u}_{\Lmd} \in E_{1}^{0,0} $ in the spectral sequence with
$ k = -1 $ is the image of
$ \overline{h}_{\Lmd} $ by
$ d_{1} $.
So there exists a homotopy
$$
h_{\Lmd} \in  G^{-1}W'_{0}\Hom^{-1}(K_{\Lmd},L_{\Lmd})
$$
such that
$ u_{\Lmd} = dh_{\Lmd} + h_{\Lmd}d $ in
$ \Hom^{0}(K_{\Lmd},L_{\Lmd}) $ for
$ \Lmd = A, {\Bbb C} $.
Put
$ h_{{\Cal D}} = \overline{h}_{{\Cal D}} $.

Here we may assume
$ {\Gr}_{G}h_{\Lmd} = \overline{h}_{\Lmd} $ by replacing
$ h_{\Lmd} $ if necessary.
In fact, let
$$
\overline{v}_{\Lmd} = \Gr_{G}h_{\Lmd} - \overline{h}_{\Lmd} \in  W'_{0}
\Hom^{-1}(\overline{K}_{\Lmd},\overline{L}_{\Lmd})
$$
for
$ \Lmd = A, {\Bbb C} $.
Then
$ d\overline{v}_{\Lmd} + \overline{v}_{\Lmd}d = 0 $ in
$ \Hom^{0}(\overline{K}_{\Lmd},\overline{L}_{\Lmd})$.
Using the second spectral sequence as above with
$ k = p + q = -1 $,
we have
$$
v_{\Lmd} \in  H^{-1}G^{-1}W'_{0}\Hom^{\ssbull }(K_{\Lmd},L_{\Lmd}),
$$
such that
$ \Gr_{G}v_{\Lmd} = \overline{v}_{\Lmd} $.
So modifying
$ h_{\Lmd} $ by
$ v_{\Lmd} $ we have
$ \Gr_{G}h_{\Lmd} = \overline{h}_{\Lmd} $ for
$ \Lmd = A, {\Bbb C} $.

Now we show that 
$ u_{\Cal D} $,
$ u_{A} $,
$ u_{\Bbb C} $ are extended to a morphism of
$ C_{\Cal H}^{b}(X,A)_{\Cal D} $.
Using the spectral sequence with
$ k= 0 $ again, we have
$$
\alpha _{\Lmd}\circ u_{\Lmd,{\Bbb C}} - 
u_{{\Bbb C}}\circ \alpha _{\Lmd} = 0\quad  \text{in }
H^{0}G^{0}W'_{0}\Hom^{\ssbull }(K_{\Lmd,{\Bbb C}},L_{{\Bbb C}}),
$$
because its image in
$ E_{1}^{0,0} =  W'_{0}\Hom^{0}(\overline{K}_{\Lmd,
{\Bbb C}},\overline{L}_{{\Bbb C}}) $ by the edge morphism is zero and
$ E_{1}^{p,q} = 0 $ for
$ p < 0 $ or
$ q < 0 $.
So there exists
$$
u'_{\Lmd} \in  G^{0}W'_{0}\Hom^{-1}(K_{\Lmd,{\Bbb C}},L_{{\Bbb C}}),
$$
such that
$ \alpha _{\Lmd}\circ u_{\Lmd,{\Bbb C}} - 
u_{{\Bbb C}}\circ \alpha _{\Lmd}
 = du'_{\Lmd} + u'_{\Lmd}d $.
Then
$ u = (u_{\Cal D}, u_{A}, u_{\Bbb C}, u'_{\Cal D}, u'_{A}) $
gives a morphism of
$ C_{\Cal H}^{b}(X,A)_{\Cal D} $.

We have to show that
$ u $ is homotopic to zero.
Let
$$
w'_{\Lmd} = u'_{\Lmd} 
-\alpha _{\Lmd}\circ h_{\Lmd,{\Bbb C}} - 
h_{{\Bbb C}}\circ \alpha _{\Lmd} \in 
G^{-1}W'_{0}\Hom^{-1}(K_{\Lmd,{\Bbb C}},L_{{\Bbb C}}) $$
for
$ \Lmd = {\Cal D}, A $.
Then
$ dw'_{\Lmd} + w'_{\Lmd}d = 0 $,
and
$$
\Gr_{G}w'_{\Lmd} = - \Gr_{G}\alpha _{\Lmd}\circ \overline{h}
_{\Lmd,{\Bbb C}} + 
\overline{h}_{{\Bbb C}}\circ \Gr_{G}\alpha _{\Lmd} = 0
$$
in
$ W'_{0}\Hom^{-1}(\overline{K}_{\Lmd,{\Bbb C}},
\overline{L}_{{\Bbb C}}) $
by assumption on
$ \overline{h}_{\Lmd} $.
This implies
$$
w'_{\Lmd} = 0\quad\text{in }H^{-1}G^{-1}W'_{0}\Hom^{\ssbull }(K_{\Lmd,
{\Bbb C}},L_{{\Bbb C}}),
$$
by the second spectral sequence as above with
$ k = -1,  p + q = - 1 $,
and we have
$$
h'_{\Lmd} \in  G^{-1}W'_{0}\Hom^{-2}(K_{\Lmd,{\Bbb C}},L_{{\Bbb C}}),
$$
such that
$ w'_{\Lmd} = dh'_{\Lmd} - h'_{\Lmd}d $.
Then
$ h = (h_{\Cal D}, h_{A}, h_{\Bbb C}, h'_{\Cal D}, h'_{A}) $ satisfies
$ u = dh + hd $.

So we get (2.7.4).
It is clear that this induces (2.7.1).
\enddemo

\bigskip\noindent
{\bf 2.8.~Proposition.} {\sl Let  
$ f : X \rightarrow  Y $ be a proper 
morphism of algebraic varieties or a projective morphism of complex 
analytic spaces.
Then the direct image functors for filtered differential complexes, filtered 
complexes of
$ {\Cal D} $-Modules, and complexes of
$ A- $ or
$ {\Bbb C} $-Modules induce the direct image functors for mixed Hodge 
complexes
$$
f_{*} : {D}_{\Cal H}^{b}(X,A) \rightarrow  {D}_{\Cal H}^{b}(Y,A),\quad f_{*} : 
{D}_{\Cal H}^{b}(X,A)_{{\Cal D}} \rightarrow  {D}_{\Cal H}^{b}(Y,A)
_{{\Cal D}},
\tag 2.8.1
$$
in a compatible way with
$ \DR^{-1} $.
In the algebraic case, they are also compatible with the direct image of 
mixed Hodge Modules [15] via the functor
$ \varepsilon  $ in (2.7).
}
\bigskip

\demo\nofrills {Proof.\usualspace}
 The direct image is defined by using the canonical flasque 
resolution of Godement which is truncated by
$ \tau _{\le i} $ for
$ i $ large enough.
This induces the functors in (2.8.1) by (iv) of (1.4).
The compatibility of the direct image with
$ \DR^{-1} $ is clear.

So it remains to show the compatibility with
$ \varepsilon  $.
Let
$ M^{\ssbull } $ be a bounded complex of mixed Hodge Modules such that
$ H^{i}f_{*}M^{j} = 0 $ for
$ i \ne  0 $.
It is enough to construct a natural isomorphism
$$
\varepsilon f_{*}M^{\ssbull } = f_{*}\varepsilon M^{\ssbull }
\tag 2.8.2
$$
in a functorial way for
$ M^{\ssbull } $.

By definition
$ f_{*}M^{\ssbull } $ is the complex whose
$ j $-th component is
$ H^{0}f_{*}M^{j} $.
Let
$ K = f_{*}\varepsilon M^{\ssbull} \in D_{\Cal H}^{b} (X,A)_{\Cal D} $, and
$ G $ be the filtration on
$ K $ induced by
$ \sigma  $ on
$ M^{\ssbull } $ as in the proof of (2.7).
If
$ M^{j} = 0 $ for
$ j \ne 0 $, the isomorphism (2.8.2) is obtained by using the truncation
$ \tau $ on
$ K $ (see Remark (i) after (2.9)), because
$ {\Cal H}^{i}K = 0 $ for 
$ i \ne 0 $.
In general, we have to use the truncations relative to the filtration
$ G $, which are
$ \Dec G $, 
$ \Decdual G $ in [1], [4].
We can construct these in our situation.
See (2.9) below.
Then we have the weak quasi-isomorphisms
$$
(\Dec G)^{0}K/(\Decdual G)^{1}K \leftarrow  (\Dec G)^{0}K \rightarrow  K,
\tag 2.8.3
$$
Here the first term is the mapping cone of
$ (\Decdual G)^{-1}K \rightarrow (\Dec G)^{0}K $ where
$ W $ is {\it not} shifted as in (2.1).
So we get the assertion, because the first term is isomorphic to
$ \varepsilon f_{*}M^{\ssbull } $.
\enddemo

\bigskip
\noindent
{\bf 2.9.} Let
$ {C}_{\Cal H}^{b}(X,A;G)_{{\Cal D}} $ denote the category of objects of
$ {C}_{\Cal H}^{b}(X,A)_{{\Cal D}} $ with a finite filtration
$ G $.
An object of
$ {C}_{\Cal H}^{b}(X,A;G)_{{\Cal D}} $ is
$$
K = ((K_{{\Cal D}};F,W,G), (K_{A};W,G), (K_{{\Bbb C}};W,G); \alpha 
_{{\Cal D}}, \alpha _{A})
$$
such that the three filtrations
$ F, W, G $ on
$ K_{{\Cal D}} $ are compatible in the sense of [14] (see also (A.1)), and
$ (G^{i}/G^{j})K $ belongs to
$ {C}_{\Cal H}^{b}(X,A) $ for any
$ i < j $.
Here
$ (G^{i}/G^{j})K $ denotes
$$
(((G^{i}/G^{j})K_{{\Cal D}};F,W), ((G^{i}/G^{j})K_{A},W), 
((G^{i}/G^{j})K_{{\Bbb C}},W); \alpha _{{\Cal D}}, \alpha _{A}).
$$

A morphism of
$ {C}_{\Cal H}^{b}(X,A;G)_{{\Cal D}} $ is a morphism of
$ {C}_{\Cal H}^{b}(X,A)_{{\Cal D}} $ such that
$ u_{\Lmd}, u'_{\Lmd} $ preserve
$ G $ (similarly for homotopy).
A weakly acyclic object is filtered acyclic 
for
$ G $.
A weak quasi-isomorphism is a filtered quasi-isomorphism for
$ G $.
We can define
$ {K}_{\Cal H}^{b}(X,A;G)_{{\Cal D}}, {D}_{\Cal H}^{b}(X,A;G)_{{\Cal D}} $ as 
before.

If
$ K \in C_{\Cal H}^{b}(X,A;G)_{\Cal D} $ satisfies
${}^{p}{\Cal H}^{j}{\Gr}_{G}^{i}{\Gr}_{k}^{W}K_{A} = 0 $ for
$ i \ne j $, and there exists a filtration
$ W' $ on
$ K_{\Cal D} $, 
$ K_{A} $, 
$ K_{\Bbb C} $, such that
$ W $ is the convolution of 
$ W' $ and
$ G $ (see the proof of (2.7)), then we have
$ M^{\ssbull} \in C^{b}\MHW(X,A) $ such that
$$
M^{i} = ({\Cal H}^{i}{\Gr}_{G}^{i}(K_{\Cal D};F,W'), 
{}^{p}{\Cal H}^{i}{\Gr}_{G}^{i}(K_{A},W'))
$$
where the isomorphism
$ \DR_{X}({\Cal H}^{i}{\Gr}_{G}^{i}K_{\Cal D}) =
({}^{p}{\Cal H}^{i}{\Gr}_{G}^{i}K_{A})\otimes_{A}{\Bbb C} $ is
induced by 
$ \alpha _{\Cal D} $,
$ \alpha _{A} $.
Furthermore,
$ \varepsilon M^{\ssbull} $
is isomorphic in
$ K_{\Cal H}^{b}(X,A) $ to
$ K $ with
$ G $ forgotten.

Using Remarks (ii), (iii) below, we have for
$ K \in  {K}_{\Cal H}^{b}(X,A;G)_{{\Cal D}} $
$$
(\Dec G)^{m}K \in {K}_{\Cal H}^{b}(X,A;G)_{{\Cal D}} 
$$
with a morphism
$ (\Dec G)^{m}K \rightarrow  K $ in
$ {K}_{\Cal H}^{b}(X,A;G)_{{\Cal D}} $, where
$ (\Dec G)^{m}K $ is
$$
(((\Dec G)^{m}K_{{\Cal D}};F,W,G), ((\Dec G)^{m}K_{A};W,G), 
((\Dec G)^{m}K_{{\Bbb C}};W,G)), 
$$
and
$ (\Dec G)^{m}K_{{\Cal D}} $,
$ (\Dec G)^{m}K_{{\Lmd}} $ for
$ \Lmd = A, {\Bbb C} $
are as in Remarks (ii), (iii) below.
(Here we assume that
$ {\Gr}_{k}^{W}{\Gr}_{G}^{j}K_{\Lmd} $ are injective for
$ \Lmd = A, {\Bbb C} $ as in the proof of (2.7).)
By Remark (ii),
$ (\Dec G)^{m}K $ is essentially unique in
$ {K}_{\Cal H}^{b}(X,A;G)_{{\Cal D}} $.
In fact, if
$ K_{a} \rightarrow  K \,\, (a = 1, 2) $ satisfies the condition of
$ (\Dec G)^{m}K $,
we have
$ K_{3} \rightarrow  K $ satisfying the condition of
$ (\Dec G)^{m}K $ together with morphisms
$ K_{3} \rightarrow  K_{a} \,\, (a  = 1, 2) $ compatible with the morphisms 
to
$ K $ in
$ {K}_{\Cal H}^{b}(X,A;G)_{{\Cal D}} $ (and the independence of
$ K_{3} $ can be proved).
We have also the functoriality of
$ (\Dec G)^{m}K $.
In fact, for a morphism
$ K \rightarrow  L $ in
$ {K}_{\Cal H}^{b}(X,A;G)_{{\Cal D}} $,
we have a unique morphism
$ (\Dec G)^{m}K \rightarrow  (\Dec G)^{m}L $ in
$ {K}_{\Cal H}^{b}(X,A;G)_{{\Cal D}} $ giving a commutative diagram.
We have a similar assertion for
$ \Decdual G $.

\noindent\bigskip
{\it Remarks.} (i)
In the case the filtration
$ G $ is trivial (i.e.,
$ {\Gr}_{G}^{i}K = 0 $ for
$ i \ne  0) $,
we get the truncation
$ \tau _{\le m} $ on the mixed 
$ A $-Hodge
$ {\Cal D} $-complexes.
The proof is by induction on the length of
$ W $.
See Remark (iv) below.
In the case $ X = pt $, this can be applied to the mixed Hodge complexes
in the sense of [4].

(ii) Let
$ {D}_{c}^{b}(X;A;G,W) $ be the full subcategory of
$ D^{b}(X,A;G,W) $ (see the proof of (2.7)) consisting of objects
$ (K;G,W) $ such that
$ {\Gr}_{G}^{j}{\Gr}_{W}^{k}K \in  {D}_{c}^{b}(X;A) $.
Then there exists the
$ t $-structure
$ (D^{\le a}, D^{\ge b}) $ on
$ {D}_{c}^{b}(X;A;G,W) $ such that
$ (K;G,W) \in  D^{\le a} $ (resp.
 $ D^{\ge b}) $ if and only if
$$
{}^{p}{\Cal H}^{i}{\Gr}_{G}^{j}W^{k}K = 0\quad  \text{for }i > j + a\,\,\, 
\text{(resp. }i < j + b).
$$
In fact, if
$ X $ is a point, we have the desired
$ t $-structure using
$ \Dec G $ in [4].
The argument is similar (up to a shift of
$ t $-structure) if
$ {\Cal H}^{i}{\Gr}_{G}^{j}{\Gr}_{W}^{k}K $ are local systems.
In general, let
$ Z $ be a pure dimensional smooth closed subvariety of
$ X $ such that
$ {\Cal H}^{i}{\Gr}_{G}^{j}{\Gr}_{W}^{k}K|_{Z} $ are local systems.
Let
$ U = X \backslash  Z $ with natural inclusions
$ i : Z \rightarrow  X, j : U \rightarrow  X $.
Assume we have the desired
$ t $-structure on
$ U $ and
$ Z $ by induction on strata.
Then the condition of
$ (K;G,W) \in  D^{\le a} $ on
$ X $ is equivalent to
$$
\aligned
{}^{p}{\Cal H}^{i}i^{*}{\Gr}_{G}^{j}W^{k}K &= 0\,\,\, \text{for }i > j + a,
\\ 
{}^{p}{\Cal H}^{i}j^{*}{\Gr}_{G}^{j}W^{k}K &= 0\,\,\, \text{for }i > j + a.
\endaligned
$$
This is equivalent to
$ i^{*}(K;G,W) \in  D^{\le a}, j^{*}(K;G,W) \in  D^{\le a} $,
and similarly for
$ D^{\ge b} $.
So the assertion follows from the gluing argument of
$ t $-structures in [2].
The associated truncation
$ \tilde{\tau}_{\le m} $ is denoted by
$ ({}^{p}\Dec G)^{-m} $.
We have the dual
$ t $-structure
$ ({}^{d}D^{\le a}, {}^{d}D^{\ge b}) $ where
$ {}^{p}{\Cal H}^{i}{\Gr}_{G}^{j}W^{k}K $ is replaced by
$ {}^{p}{\Cal H}^{i}{\Gr}_{G}^{j}(K/W^{k}K) $ in the condition.
The associated truncation
$ \tilde{\tau}_{\le m} $ is denoted by
$ ({}^{p}\Decdual G)^{-m} $.

(iii) Let
$ (M,F) $ be a complex of filtered
$ {\Cal D} $-Modules with decreasing finite filtrations
$ W, G $ such that
$ F, W, G $ are compatible on
$ {M}_{U\rightarrow V}^{j} $ for any
$ \{ U \rightarrow  V\} \in  LE(X) $ (see (1.3)) in the sense of [14] (see also 
(A.1)).
Let
$ G' = \Dec G $.
See [4] (and also (A.3)).
Assume the three filtrations
$ F, W, G' $ on
$ {\Gr}_{G}^{i}{M}_{U\rightarrow V}^{j} $ are compatible.
(This condition is satisfied in the case of mixed
$ A $-Hodge
$ {\Cal D} $-complexes by Remark (iv) below.)
Then the four filtrations
$ F, W, G, G' $ on
$ {M}_{U\rightarrow V}^{j} $ are compatible.
In fact, it is enough to verify that  natural morphisms
$$
F^{p}W^{k}G^{\prime q}G^{i}{M}_{U\rightarrow V}^{j} \rightarrow  
F^{p}W^{k}G^{\prime q}{\Gr}_{G}^{i}{M}_{U\rightarrow V}^{j}
$$
are surjective by [14, 1.2.12].
But this follows from the compatibility of the filtrations
$ F, G, W $ on
$ {M}_{U\rightarrow V}^{j} $,
because the latter implies the surjectivity of
$$
F^{p}W^{k}G^{i}{M}_{U\rightarrow V}^{j} \rightarrow  
F^{p}W^{k}{\Gr}_{G}^{i}{M}_{U\rightarrow V}^{j}
$$
by loc.~cit, and
$ G^{\prime q}G^{i}{M}_{U\rightarrow V}^{j} $ is the inverse image of
$ G^{\prime q}{\Gr}_{G}^{i}{M}_{U\rightarrow V}^{j} $ by the projection
$ G^{i}{M}_{U\rightarrow V}^{j} \rightarrow  
{\Gr}_{G}^{i}{M}_{U\rightarrow V}^{j} $.

Similar argument applies to the dual notion
$ \Decdual G $.
See (A.3).

(iv) Let
$ {\Cal C}, {\Cal A} $ be as in (A.1) where (A.2.2) is 
assumed.
Let
$ (K,W) $ be a complex of
$ {\Cal C} $ with an increasing finite filtration
$ W $.
Then for
$ k \in  {\Bbb Z} $,
 $ W_{k}K $ is strict if
$ W_{k-1}K, {\Gr}_{k}^{W}K $ are strict and the morphism
$ H^{j}{\Gr}_{k}^{W}K \rightarrow  H^{j+1}W_{k-1}K $ associated with the 
distinguished triangle
$$
\rightarrow  W_{k-1}K \rightarrow  W_{k}K \rightarrow  {\Gr}_{k}^{W}K 
\rightarrow \quad  \text{in }{\Cal A}
$$
is a strict morphism in
$ {\Cal C} $.

Assume
$ W_{k}K $ and
$ {\Gr}_{k}^{W}K $ are strict and
$ H^{j}{\Gr}_{k}^{W}K \rightarrow  H^{j+1}W_{k-1}K $ are strict morphisms 
of
$ {\Cal C} $ for any
$ k $.
Then the induced filtration
$ W $ on the canonical truncation
$ \tau _{\le m}K $ in
$ {\Cal A} $ gives a filtration in
$ {\Cal C} $ (i.e.,
$ {\Gr}_{k}^{W}\tau _{\le m}K^{m} \in  {\Cal C} $,
where
$ \tau _{\le m}K^{m} = \Ker \,d^{m}) $.
This follows from the snake lemma applied to
$$
\CD
0 @>>>    W_{k-1}K^{j}   @>>>    W_{k}K^{j}   @>>>   
{\Gr}_{k}^{W}K^{j}   @>>>    0 \\
@. @VVV @VVV @VVV \\
0 @>>>    W_{k-1}K^{j+1}   @>>>    W_{k}K^{j+1}   @>>>   
{\Gr}_{k}^{W}K^{j+1} @>>> 0
\endCD
$$
because
$ \Ker \,{\Gr}_{k}^{W}d^{j} \rightarrow  \Coker \, W_{k-1}d^{j} $ is 
factorized by
$ H^{j}{\Gr}_{k}^{W}K \rightarrow  H^{j+1}W_{k-1}K $,
and
$$
\Ker \,{\Gr}_{k}^{W}d^{j} \rightarrow  H^{j}{\Gr}_{k}^{W}K,\quad 
H^{j+1}W_{k-1}K \rightarrow  \Coker \, W_{k-1}d^{j} 
$$
are strictly surjective or injective morphisms.

We have furthermore
$$
\aligned
H^{j}W_{k}\tau _{\le m}K 
&= H^{j}W_{k}K\,\,\, \text{for }j \le  m,\,\,\, \text{and 0 } 
\text{otherwise.} \\
H^{m}{\Gr}_{k}^{W}\tau _{\le m}K 
&= \Coker(H^{m}W_{k-1}K \rightarrow  H^{m}W_{k}K) \\
&= \Ker(H^{m}{\Gr}_{k}^{W}K \rightarrow  H^{m+1}W_{k-1}K)
\endaligned
$$
using the morphism of long exact sequences
$$
\CD
H^{m}W_{k-1}\tau _{\le m}K @>>>   H^{m}W_{k}\tau _{\le m}K    
@>>>  H^{m}{\Gr}_{k}^{W}\tau _{\le m}K @>>>      0 \\
@\vert @\vert @VVV @VVV \\
H^{m}W_{k-1}K @>>>  H^{m}W_{k}K @>>> H^{m}{\Gr}_{k}^{W}K 
@>>> H^{m+1}W_{k-1}K.
\endCD
$$
We have the dual assertion using
$$
\rightarrow  {\Gr}_{k}^{W}K \rightarrow  K/W_{k-1}K \rightarrow  
K/W_{k}K \rightarrow 
$$
instead of
$ \rightarrow  W_{k-1}K \rightarrow  W_{k}K \rightarrow  {\Gr}_{k}^{W}K 
\rightarrow  $.

\bigskip\noindent
{\bf 2.10.~Proposition.} {\sl
Assume
$ X = pt $.
Let
$ \MHS(A)^{p} $ be the category of (graded) polarizable mixed
$ A $-Hodge structures, which is naturally identified with
$ \MHM(pt,A) $.
See [15].
Then a natural functor
$$
\varepsilon : D^{b}\MHS(A)^{p} \rightarrow  {D}_{\Cal H}^{b}(pt,A)\,\,
(= {D}_{\Cal H}^{b}(pt,A)_{\Cal D})
$$
in (2.7) (which is trivial in this case) is an equivalence of categories.
}
\bigskip

\demo\nofrills {Proof.\usualspace}
For a polarizable mixed
$ A $-Hodge structure
$ H $ and a weak quasi-isomorphism 
$ K \to L $ in
$ C_{\Cal H}^{b}(pt,A) $, we show the isomorphism
$$
\Hom_{K^{b}} (H,K) \simto \Hom_{K^{b}} (H,L)
\tag 2.10.1
$$
where
$ K^{b} $ means 
$ K_{\Cal H}^{b}(pt,A) $, and
$ \MHS(A)^{p} $ is naturally identified with a subcategory of
$ C_{\Cal H}^{b}(pt,A) $. 
Let
$ A^{H} $ be the Hodge structure of type
$ (0,0) $ whose underlying
$ A $-module is
$ A $.
Then the assertion is reduced to the case
$ H = A^{H} $ by replacing
$ K $ with
$ H^{*}\otimes K $,
where
$ H^{*} $ is the dual mixed Hodge structure of
$ H $, and
$ L := H^{*}\otimes K $ is defined by
$ L_{\Lmd} = H_{\Lmd}^{*}\otimes K_{\Lmd} $ for
$ \Lmd = F, A, {\Bbb C} $.
(Here
$ H_{F} = H_{\Bbb C} $.)
Since
$ K_{\Cal H}^{b}(pt,A) $ is a triangulated category, it is enough 
to show that if
$ K $ is weakly acyclic, we have
$$
\Hom_{K^{b}} (A^{H},K) = 0.
\tag 2.10.2
$$

Let
$ u = (u_{F}, u_{A}, u_{{\Bbb C}}, u'_{F}, u'_{A}) : A^{H} \rightarrow  K $ be a 
morphism of
$ {C}_{\Cal H}^{b}(X,A) $,
and put
$$
m_{\Lmd} = u_{\Lmd}1 \,\, (\Lmd = F, A, {\Bbb C}),
\quad m'_{\Lmd} = u'_{\Lmd}1 \,\, (\Lmd = F, A).
$$
Then
$ m = (m_{F}, m_{A}, m_{{\Bbb C}}, m'_{F}, m'_{A}) $ belongs to
$$
 F_{0}W_{0}{K}_{F}^{0} \oplus  W_{0}{K}_{A}^{0} \oplus  W_{0}{K}
_{\Bbb C}^{0} \oplus  W_{0}{K}_{\Bbb C}^{-1} \oplus  W_{0}{K}_{\Bbb C}^{-1},
$$
and satisfies
$$
dm_{\Lmd} = 0 \,\, (\Lmd = F, A, {\Bbb C}),\quad
\alpha _{\Lmd}m_{\Lmd} - m_{{\Bbb C}} = dm'_{\Lmd} \,\,(\Lmd = F, A),
$$
where
$ \alpha _{A} $ denotes also its composition with the natural inclusion
$ K_{A} \rightarrow  K_{A,{\Bbb C}} $.

Let
$ h = (h_{F}, h_{A}, h_{{\Bbb C}}, h'_{F}, h'_{A}) : A^{H} \rightarrow  K[-1] $ 
be a homotopy, and put
$$
n_{\Lmd} = h_{\Lmd}1 \,\, (\Lmd = F, A, {\Bbb C}),
\quad n'_{\Lmd} = h'_{\Lmd}1 \,\, (\Lmd = F, A).
$$
Let
$ W' = \Dec W $,
$ W'' = \Decdual W $.
See (A.3).
Then
$ n = (n_{F}, n_{A}, n_{{\Bbb C}}, n'_{F}, n'_{A}) $ belongs to
$$
F_{0}W'_{0}{K}_{F}^{-1} \oplus  W'_{0}{K}_{A}^{-1} \oplus  W'_{0}{K}
_{\Bbb C}^{-1} \oplus  W_{1}{K}_{\Bbb C}^{-2} \oplus  W_{1}{K}_{\Bbb C}^{-2},
$$
and satisfies
$$
\alpha _{\Lmd}n_{\Lmd} - n_{{\Bbb C}} + dn'_{\Lmd} \in  W_{0}{K}_{\Bbb C}^{-1} 
\,\, (\Lmd = F, A).
$$
By (2.1.2), the condition
$ u = dh + hd $ corresponds to the relation
$$
m_{\Lmd} = dn_{\Lmd} \,\, (\Lmd = F, A, {\Bbb C}),
\quad m'_{\Lmd} = dn'_{\Lmd} + \alpha _{\Lmd}n_{\Lmd} - n_{{\Bbb C}} \,\, 
(\Lmd = F, A).
$$
Let
$ Z $ (resp.
 $ R) $ denote the group consisting of
$ m $ (resp.
 $ n) $ satisfying the above conditions.
Then we have the morphism
$ \phi :  R \rightarrow  Z $ by the above relation so that
$ \Coker \, \phi $ is isomorphic to
$ \Hom_{K^{b}}(A^{H}, K) $.
We have furthermore
$$
\Hom_{K^{b}}(A^{H}, K) = \Coker \, \phi  = \Coker \, \phi '.
\tag 2.10.3
$$
where the morphism
$ \phi' : R' \to Z' $ is defined by
$$
(\tilde{n}_{F}, \tilde{n}_{A}, \tilde{n}') \rightarrow  (d\tilde{n}_{F}, 
d\tilde{n}_{A}, d\tilde{n}' + \alpha _{F}\tilde{n}_{F} - \alpha 
_{A}\tilde{n}_{A}),
$$
with
$$
\aligned
Z' =
&\bigl\{(\tilde{m}_{F}, \tilde{m}_{A}, \tilde{m}') 
\in  F_{0}W_{0}{K}_{F}^{0} \oplus  W_{0}{K}_{A}^{0} \oplus  
W_{0}{K}_{\Bbb C}^{-1} :
\\
&\qquad\qquad d\tilde{m}_{\Lmd} = 0 \,\,
(\Lmd = F, A),\quad \alpha _{F}\tilde{m}_{F} - \alpha _{A}\tilde{m}_{A} - 
d\tilde{m}' = 0\bigr \},
\\
R' = 
&\bigl\{(\tilde{n}_{F}, \tilde{n}_{A}, \tilde{n}') 
\in  F_{0}W'_{0}{K}_{F}^{-1} \oplus  W'_{0}{K}_{A}^{-1} \oplus  
W_{1}{K}_{\Bbb C}^{-2} : 
\\
&\qquad\qquad
d\tilde{n}' + \alpha _{F}\tilde{n}_{F} - \alpha _{A}\tilde{n}_{A} \in  
W_{0}{K}_{\Bbb C}^{-1} \bigr\}.
\endaligned
$$
Note that the conditions imply
$$
\alpha _{F}\tilde{m}_{F} - \alpha _{A}\tilde{m}_{A}
\in W''_{-1}{K}_{\Bbb C}^{0} , \quad
\alpha _{F}\tilde{n}_{F} - \alpha _{A}\tilde{n}_{A}
\in W''_{-1}{K}_{\Bbb C}^{-1} .
$$
For the proof of (2.10.3) we have the surjective morphisms
$ q_{Z} : Z \rightarrow  Z', q_{R} : R \rightarrow  R' $ defined by
$$
m \rightarrow  (m_{F}, m_{A}, m'_{F} - m'_{A})\,\,\, \text{(same for }n).
$$
so that
$ q_{Z}\phi  = \phi 'q_{R} $
and
$ \phi  : \Ker \,q_{R} \rightarrow  \Ker \,q_{Z} $ is surjective.

Now we have to show
$ \Coker\,\phi' = 0 $ if
$ K $ is weakly acyclic.
Take
$ (\tilde{m}_{F}, \tilde{m}_{A}, \tilde{m}') \in Z' $.
By definition of
$ W' $,
$ W'' $ we have
$$
W'_{k}K_{\Lmd}^{i}/W''_{k-1}K_{\Lmd}^{i} = H^{i}\Gr_{k-i}^{W}K_{\Lmd}.
$$
It is the
$ E_{1} $-term of the weight spectral sequence (2.3.1), 
and the 
$ E_{1} $-complex is acyclic by the acyclicity of
$ K $ and the
$ E_{2} $-degeneration of the spectral sequence.
By definition of
$ Z' $, the images of
$ \tilde{m}_{F} $ and
$ \tilde{m}_{A} $ in
$ H^{0}\Gr_{0}^{W}K_{\Bbb C} $ coincide, and
$ d_{1} $ of the spectral sequence is a morphism of
Hodge structures.
So we may assume
$$
\tilde{m}_{F} \in F^{0}W''_{-1}K_{F}^{0} , \quad
\tilde{m}_{A} \in W''_{-1}K_{A}^{0} ,
$$
by replacing
$ \tilde{m} $ with
$ \tilde{m} - \phi '(\tilde{n}) $ for some
$ \tilde{n} = \{\tilde{n}_{F}, \tilde{n}_{A}, \tilde{n}_{\Bbb C}\} 
\in R' $ such that
$$
\tilde{n}_{F} \in F^{0}W'_{0}K_{F}^{-1} , \quad
\tilde{n}_{A} \in W'_{0}K_{A}^{-1} . 
$$
In fact, we have
$ \tilde{n}_{F} \in F_{0}W'_{0}K_{F}^{-1} $,
$ \tilde{n}_{A} \in W'_{0}K_{A}^{-1} $ such that
$ \tilde{m}_{\Lmd} - d\tilde{n}_{\Lmd} \in W''_{-1} $ for
$ \Lmd = F, A $ and the images of
$ \tilde{n}_{F} $,
$ \tilde{n}_{A} $ in
$ H^{-1}\Gr_{1}^{W}K_{\Bbb C} $ coincide, i.e.,
$$
\alpha_{F}\tilde{n}_{F} - \alpha_{A}\tilde{n}_{A} 
\in W''_{-1}H_{\Bbb C}^{-1} =
d(W_{1}K_{\Bbb C}^{-2}) + W_{0}K_{\Bbb C}^{-1}.
$$
Then we can choose
$ \tilde{n}' $ such that the condition of
$ R' $ is satisfied.

We may assume furthermore
$ \tilde{m}_{\Lmd} = 0 $ for
$ \Lmd = F, A $ by replacing
$ \tilde{m} $ with
$  \tilde{m} - \phi '(\tilde{n}) $, because
$ F^{0}W''_{0}K_{F} $,
$ W''_{0}K_{A} $ are acyclic, and
$ d (F^{0}W''_{-1}K_{F}^{-1}) = d(F^{0}W_{0}K_{F}^{-1}) $, etc.
Then the assertion follows from the acyclicity of
$ W''_{-1}K_{\Bbb C} $.
This proves (2.10.1--2).

So we get the isomorphism
$$
\Hom_{K^{b}}(H,K) \simto \Hom_{D^{b}}(H,K)
\tag 2.10.4
$$
for
$ H \in \MHS(A)^{p} $,
$ K \in {C}_{\Cal H}^{b}(pt,A) $, where
$ D^{b} $ means
$ {D}_{\Cal H}^{b}(pt,A) $.
Then using the truncation
$ \tau $ in (2.9), it is enough to show
$$
\aligned
&\Hom_{K^{b}}(H', H[i]) = 0 \quad\text{for } i > 1, \\
&\Hom_{\MHS(A)^{p}}(H', H) = \Hom_{K^{b}}(H', H), \\
&\Ext^{1}_{\MHS(A)^{p}}(H', H) = \Hom_{K^{b}}(H', H[1]),
\endaligned
\tag 2.10.5
$$
for
$ H, H' \in \MHS(A)^{p} $.
Here we may assume
$ H' = A^{H} $ as before.
Then the first two assertions are clear by (2.10.3), 
and the last  is reduced to the following :
 
For
$ H = ((H_{\Bbb C}; F,W), (H_{A},W)) \in \MHS(A)^{p} $, we have
$$
\Ext^{1}_{\MHS(A)^{p}}(A^{H}, H) = 
\Coker((F^{0}W_{0}H_{\Bbb C} \oplus  W_{0}H_{A})^{(0)} 
\rightarrow  W_{-1}H_{\Bbb C}),
$$
where
$ (F^{0}W_{0}H_{\Bbb C} \oplus  W_{0}H_{A})^{(0)} = 
\Ker((F^{0}W_{0}H_{\Bbb C} \oplus  
W_{0}H_{A}) \rightarrow  {\Gr}_{0}^{W}H_{\Bbb C}) $.

This can be verified, for example, by using
$$
\Ext^{1}_{\MHS(A)^{p}}(A^{H}, W^{-1}H) = 
\Coker((F^{0}W_{-1}H_{\Bbb C} \oplus  W_{-1}H_{A})
\rightarrow  W_{-1}H_{\Bbb C}),
$$
(see [3], [13])) together with the long exact 
sequence associated with
$
0 \rightarrow  W_{-1}H \rightarrow  W_{0}H \rightarrow  {\Gr}_{0}^{W}H 
\rightarrow  0.
$
See also [1, 2.2].

We can show that the last isomorphism of (2.10.5) is actually
induced by the functor
$ \varepsilon $, by using the quasi-isomorphism
$ H[1] \to C(\tilde{H} \to A^{H}) $, where
$ 0 \to H \to \tilde{H} \to A^{H} \to 0 $ is an exact sequence of
$ \MHS(A)^{p} $.
\enddemo

\bigskip\noindent
{\it Remark.}
The reader may notice the difference between the proofs of
(2.10) and [1, 3.4].
In the latter,
$ \Hom_{K^{b}}(H,K) $ depends on the resolution of
$ K $, and the effaceability of
$ \Gamma^{1} $ in loc.~cit is shown by taking a good resolution of
$ K $.
The mixed Hodge complexes on
$ pt $ in this paper are closer to those in [4] and also to
$ \tilde{p} $-Hodge complexes in [1].
It may be easier to prove Lemma 3.11 in [1] by showing the
effaceability of
$ \Gamma^{1} $ as in the proof of 3.4, than showing that
the constructed functor is a right quasi-inverse.

\bigskip\bigskip
\centerline{\bf Appendix to \S 2}

\bigskip\noindent
We review here some facts from the theories of compatible filtrations [14]
and the realization functor [2].

\bigskip
\noindent
{\bf A.1.} Let
$ {\Cal A} $ be an abelian category.
We say that
$ E_{1}, \dots , E_{n} $ are compatible subobjects of
$ E \in  {\Cal A} $ if there exists an
$ n $-ple  complex of short exact sequences
$ \{K^{\nu }\} $ (i.e.,
 $ K^{\nu } = 0 $ if
$ |\nu _{i}| > 1 $ for some
$ i $,
and
$ K^{\nu -1_{i}} \rightarrow  K^{\nu } \rightarrow  K^{\nu +1_{i}} $ are 
exact for any
$ \nu  \in  {\Bbb Z}^{n}) $ such that
$$
K^{0} = E,\quad K^{-1_{i}} = E_{i},
$$
where the
$ j $-th component of  
$ 1_{i} \in  {\Bbb Z}^{n} $ is
$ \delta _{i,j} $.
We say that filtrations
$ F_{1}, \dots , F_{n} $ of
$ E \in  {\Cal A} $ are compatible filtrations, if
$ {F}_{1}^{{p}_{1}}E, \dots , {F}_{n}^{{p}_{n}}E $ are compatible subobjects 
of
$ E $ for any
$ p_{1}, \dots , p_{n} \in  {\Bbb Z} $.
Let
$ {F}_{n}^{c}({\Cal A}) $ denote the additive category of objects of
$ {\Cal A} $ with compatible
$ n $ filtrations.
We say that
$$
0 \rightarrow  (M'; F_{1}, \dots , F_{n}) \rightarrow  (M; F_{1}, \dots , 
F_{n}) \rightarrow  (M''; F_{1}, \dots , F_{n}) \rightarrow  0
$$
is a short exact sequence of
$ {F}_{n}^{c}({\Cal A}) $ if
$$
0 \rightarrow  {F}_{1}^{{p}_{1}} \dots {F}_{n}^{{p}_{n}}M' \rightarrow  
{F}_{1}^{{p}_{1}} \dots {F}_{n}^{{p}_{n}}M \rightarrow  {F}_{1}^{{p}_{1}} 
\dots {F}_{n}^{{p}_{n}}M'' \rightarrow  0\,\,\,
$$
are exact in
$ {\Cal A} $ for any
$ -\infty  \le  p_{i} < +\infty  $.
Then
$ {F}_{n}^{c}({\Cal A}) $ is an exact category with those exact sequences.
See (A.2) below.

For
$ (M; F_{1}, \dots , F_{n}) \in  {F}_{n}^{c}({\Cal A}) $ and
$ p_{i} \in  {\Bbb Z} $,
we see that there exists an
$ n $-ple complex of short exact sequences
$ \{K^{\nu }\} $ such that
$$
{F}_{1}^{{p}_{1}} \dots {F}_{n}^{{p}_{n}}M = K^{0},\quad {F}_{1}^{{p}_{1}} 
\dots {F}_{i}^{{p}_{i}+1} \dots {F}_{n}^{{p}_{n}}M = K^{-1_{i}}.
$$
This implies the commutativity of
$ {\Gr}_{{F}_{i}}^{{p}_{i}} \,\, (1 \le  i \le  n) $,
because the restriction to
$ \nu _{i} = 1 $ corresponds to the passage to
$ {\Gr}_{{F}_{i}}^{{p}_{i}} $.

\bigskip
\noindent
{\bf A.2.} Let
$ {\Cal C} $ be an exact category.
We choose and fix an abelian category
$ {\Cal A} $ such that
$ {\Cal C} $ is an additive full subcategory of
$ {\Cal A} $ which is stable by extensions in
$ {\Cal A} $,
and the short exact sequences of
$ {\Cal C} $ are the short exact sequences of
$ {\Cal A} $ whose components belong to
$ {\Cal C} $.
A morphism
$ u $ of
$ {\Cal C} $ is called strict if it has a factorization
$ u''u' $ such that
$ u' $ is a strict epimorphism and
$ u'' $ is a strict monomorphism, or equivalently, if
$ \Ker \,\, u, \Im \, u, \Coker \, u $ in
$ {\Cal A} $ belong to
$ {\Cal C} $.

A filtration
$ F $ of
$ E \in  {\Cal C} $ is a filtration of
$ E $ in
$ {\Cal A} $ such that
$$
F^{p}E/F^{q}E \in  {\Cal C}\quad  \text{for }-\infty  \le  p < q \le  +\infty .
$$
(See [14, 1.3.1].)
Let
$ F({\Cal C}) $ denote the category of filtered objects of
$ {\Cal C} $,
where
$$
0 \rightarrow  (E',F) \rightarrow  (E,F) \rightarrow  (E'',F) \rightarrow  0
$$
is called exact if it is filtered exact in
$ {\Cal A} $,
or equivalently, if
$$
0 \rightarrow  F^{p}E' \rightarrow  F^{p}E \rightarrow  F^{p}E'' \rightarrow  
0\,\,\,
$$
are exact for any
$ -\infty  \le  p < +\infty  $.
Then
$ F({\Cal C}) $ is also an exact category.
In fact,
$ F(C) $ is a full subcategory of the abelian category
$ {\tilde F}(\Cal A) $ whose objects are
$ \{ M^{p} \}_{p \in {\Bbb Z} \cup \{ - \infty \} } \in
\prod {\Cal A} $ with morphisms
$ M^{p} \rightarrow M^{q} $
$ (\infty > p > q \ge - \infty) $ compatible with compositions.
Let
$$
F^{n}({\Cal C}) = F(F^{n-1}({\Cal C}))  \,\, (n > 0),
\quad F^{0}({\Cal C}) = {\Cal C}.
$$
(Similarly for
$ {\tilde F}^{n}(\Cal A)$.)
Then
$ F^{n}(\Cal A) $
is a full subcategory of the abelian category
$ {\tilde F}^{n}(\Cal A)$ satisfying the above conditions.

If
$ {\Cal C} $ is an abelian category
$ {\Cal A} $,
we have an equivalence of categories
$$
F^{n}({\Cal A}) = {F}_{n}^{c}({\Cal A})
\tag A.2.1
$$
in the notation of Remark above.

We say that a complex
$ K $ of
$ {\Cal C} $ is cohomologically strict relatively to
$ {\Cal A} $ if
$ H^{j}K $ in
$ {\Cal A} $ belong to
$ {\Cal C} $,
and strict if
$ d^{j} : K^{j} \rightarrow  K^{j+1} $ and
$ \Im \, d^{j-1} \rightarrow  \Ker \,\, d^{j} $ are strict.

We assume for any complex
$ K $ of
$ {\Cal C}  $:
$$
\text{$ K $ is strict if and only if $ K $ is cohomologically strict
relatively to
$ {\Cal A} $.}
\tag A.2.2
$$
This condition is satisfied if
$ {\Cal C} $ is
$ {F}_{n}^{c}({\Cal A}) $ as above.
See [14, 1.2.4].

\bigskip
\noindent
{\it Remark.}
In the case
$ n = 2 $, for example, this notion of strictness (i.e., the bistrictness)
is different from the separate strictness
(i.e., the strictness for each filtration).
If we define
$ \Ker $,
$ \Coker $,
$ \Im $,
$\Coim $ in 
$ F_{2}^{c}({\Cal A}) = F_{2}({\Cal A}) $ as usual, 
and define the strictness of a morphism by
$ \Im = \Coim $, then we get the separate strictness.
For the correct definition of  bistrictness, we cannot
avoid the problem of three filtrations, because the third filtration
is defined by
$ \Ker $ or
$ \Im $.
(This becomes necessary, for example,
when we define the truncation
$ \tau $ on a mixed Hodge complex in the sense of [1, 3.2].)

\bigskip
\noindent
{\bf A.3.} With the above notation and assumptions, let
$ K $ be a complex of
$ {\Cal C} $ with a decreasing filtration
$ G $.
Assume
$$
{\Gr}_{G}^{i}d : {\Gr}_{G}^{i}K^{j} \rightarrow  {\Gr}_{G}^{i}K^{j+1}\,\,\, 
\text{are strict for any }i, j.
\tag A.3.1
$$
Then
$ \Dec G $ is defined by
$$
(\Dec G)^{i}K^{j} = \Ker(d : G^{i+j}K^{j} \rightarrow  
\Gr_{G}^{i+j}K^{j+1}).
\tag A.3.2
$$
(See [4] in the case
$ {\Cal C} = {\Cal A} $.)
More precisely, it is defined in
$ {\Cal A} $,
and belongs to
$ {\Cal C} $,
because
$$
\aligned
G^{i+j+1}K^{j} \subset  (\Dec G)^{i}K^{j} &\subset  G^{i+j}K^{j}\quad  
\text{with } 
\\
(\Dec G)^{i}K^{j}/G^{i+j+1}K^{j} &= \Ker \,{\Gr}_{G}^{i+j}d,
\\
G^{i+j}K^{j}/(\Dec G)^{i}K^{j} &=  \text{Coim  } {\Gr}_{G}^{i+j}d.
\endaligned
$$
(See also [14, 1.3.7].)
Similarly
$ \Decdual G $ is defined by
$$
(\Decdual G)^{i}K^{j} = G^{i+j}K^{j} + \Im(d : G^{i+j-1}K^{j-1} \rightarrow  
K^{j}),
\tag A.3.3
$$
and
$ (\Dec G)^{i}K \rightarrow  (\Decdual G)^{i}K $ is a quasi-isomorphism in
$ {\Cal A} $.
(This does not hold if
$ {\Gr}_{G}^{i}d $ are not strict.
See (2.9).)
Note that
$ \Dec $ (i.e.
the passage of
$ G $ to
$ \Dec G) $ commutes with an exact functor from
$ {\Cal C} $.
For example, if
$ {\Cal C} $ is the category of filtered objects
$ (M,F) $ of an abelian category,
$ \Dec $ commutes with the functor
$ (M,F) \rightarrow  F^{p}M/F^{q}M $.

\bigskip\noindent
{\it Remarks.} (i)
Assume
$ G $ is a finite filtration, and the complexes
$ {\Gr}_{G}^{i}K $ and
$ (E_{r},d_{r}) \,\, (r\ge 1) $ of the spectral sequence associated with
$ G $ (which is defined in
$ {\Cal A}) $ are strict so that the
$ {E}_{r}^{p,q} $ belong to
$ {\Cal C} $ inductively.
Then
$ {E}_{\infty }^{p,q} \in  {\Cal C} $ so that we have a spectral sequence in
$ {\Cal C} $,
and
$ K $ is also strict by (A.2.2).
(This is a generalization of the lemma on two filtrations in [4].)

(ii) If furthermore the spectral sequence degenerates at
$ E_{2} $,
then
$ (K, \Dec G) $ is cohomologically strict, because it is strict as a filtered 
complex of
$ {\Cal A} $ by the
$ E_{2} $-degeneration of the spectral sequence [4] and we have
$$
{\Gr}_{\Dec G}^{i}H^{j}K = {E}_{\infty }^{i+j,-i} \in  {\Cal C}
\tag A.3.4
$$
so that
$ (H^{j}K,\Dec G) \in  F({\Cal C}) $.
(This can be applied to
$ (K_{{\Cal D}};F,W) $ in (2.3) where
$ G = W, {\Cal C} = MF({\Cal D}_{X}) $ and
$ {\Cal A} $ is the category of graded
$ \bigoplus _{p} F_{p}{\Cal D}_{X} $-Modules.
Here
$ X $ is assumed smooth by taking a local embedding, because the 
assertion is local.)

\bigskip
\noindent
{\bf A.4.} Let
$ {\Cal A} $ be an abelian category having enough injectives.
We denote by
$ C^{*}({\Cal A};G,W) $ the category of complexes of objects of
$ {\Cal A} $ which satisfy the bounded condition corresponding to
$ * = +, -, b $,
and have two finite filtrations
$ G, W $.
We define
$ K^{*}({\Cal A};G,W) $,
$ D^{*}({\Cal A};G,W) $ as in [2], [21].
(Similarly for
$ C^{*}({\Cal A};W) $,
$ K^{*}({\Cal A};W) $,
$ D^{*}({\Cal A};W) $, etc.)
Let
$ {\Cal D} $ be a triangulated full subcategory of
$ D^{b}({\Cal A}) $ which has a
$ t $-structure whose heart is denoted by
$ {\Cal C} $.
See [2].
Let
$ D^{b}({\Cal A};G_{\text{b\^ete}},W) $ denote the full subcategory of
$ D^{b}({\Cal A};G,W) $ defined by the condition :
$$
{\Gr}_{G}^{j}{\Gr}_{W}^{k}K[j] \in  {\Cal C}
\tag A.4.1
$$
for
$ (K; G, W) \in  D^{b}({\Cal A};G,W) $.
This implies
$$
{\Gr}_{G}^{j}(W^{a}/W^{b})K[j] \in  {\Cal C},
$$
and we get a filtration
$ W $ on
$ \overline{K}^{i} := {\Gr}_{G}^{i}K[i] \in  {\Cal C}\quad  $ such that
$$
(W^{a}/W^{b})\overline{K}^{i} = {\Gr}_{G}^{i}(W^{a}/W^{b})K[i]\quad  \text{in }
{\Cal C}.
\tag A.4.2
$$
Using the short exact sequence
$ 0 \rightarrow  {\Gr}_{G}^{i+1} \rightarrow  G^{i}/G^{i+2} \rightarrow  
{\Gr}_{G}^{i} \rightarrow  0 $,
we have the morphism
$ d : (\overline{K}^{i}, W) \rightarrow  (\overline{K}^{i+1}, W) $ in
$ \Cal C $ such that
$ d^{2} = 0 $.
See [2, 3.1.7].
(Note that, if the filtration $ G $ on
$ K^{i} $ splits in a compatible way with 
$ W $ (by replacing
$ (K; G, W) $ with an injective resolution,  see (A.7) below), then
we get a representative of
$ d : (\overline{K}^{i}, W) \rightarrow  (\overline{K}^{i+1}, W) $
as in the construction of
$ \varphi $ in [21, 2.4].)

Thus we get a functor
$$
D^{b}({\Cal A};G_{\text{b\^ete}},W) \rightarrow  C^{b}({\Cal C};W)
\tag A.4.3
$$
such that
$ W^{a}/W^{b} $ of the image of
$ (K; G, W) $ by (A.4.3) is naturally isomorphic to the image of
$ ((W^{a}/W^{b})K, G) \in  D^{b}({\Cal A};G_{\text{b\^ete}}) $ by the functor
$$
D^{b}({\Cal A};G_{\text{b\^ete}}) \rightarrow  C^{b}{\Cal C}
\tag A.4.4
$$
in [2].
By the same argument as in [2, 3.1.8], 
(A.4.3) is an equivalence of categories.
(See (A.5) below.)
Taking the composition of the quasi-inverse of (A.4.3) with the forgetful 
functor which forget the filtration
$ G,  $ we get
$$
C^{b}({\Cal C};W) \rightarrow  D^{b}({\Cal A};W),
\tag A.4.5
$$
and this induces
$$
D^{b}({\Cal C};W) \rightarrow  D^{b}({\Cal A};W).
$$
The argument is the same as in [2, 3.1.10].
In fact, (A.4.5) factors through
$ K^{b}({\Cal C};W) $ (using the morphism of the spectral sequences in (A.5)), 
and preserves mapping cones.
So it is enough to show that if the image of
$ {\Gr}_{W}^{k}(K,G) $ by (A.4.4) is acyclic for
$ (K; G, W) \in  D^{b}({\Cal A};G_{\text{b\^ete}},W) $,
then
$ {\Gr}_{W}^{k}K $ is acyclic.
But this is proved in loc.~cit.

\bigskip\noindent
{\bf A.5.~Proposition.} {\sl
With the above notation, 
the functor (A.4.3) is an equivalence of categories.
(See  [2, 3.1.8] for the case
$ W $ is trivial.)
}
\bigskip

\demo\nofrills {Proof.\usualspace}
Let
$ (K; G, W) \in  C^{b}({\Cal A};G,W) $,
$ (L; G, W) \in  C^{+}({\Cal A};G,W) $
such that they are isomorphic to objects of
$ D^{b}({\Cal A};G_{\text{b\^ete}},W) $ in
$ D^{+}({\Cal A};G,W) $ and the components of
$ {\Gr}_{G}^{j}{\Gr}_{W}^{k}L $ are injective.
Then we have
$$
({\Gr}_{G}^{p}\Hom^{\ssbull }(K, L),W) = \prod _{j-i=p} (\Hom^{\ssbull }
({\Gr}_{G}^{i}K, {\Gr}_{G}^{j}L),W),
\tag A.5.1
$$
(see (A.7) below), where
$ \Hom^{\ssbull }(K,L) $ is a bifiltered complex whose
$ i^{th} $ component is
$ \prod _{m-n=i} \Hom(K^{n},L^{m}) $, and has the filtrations
$ G, W $ in the usual way.
We have furthermore
$$
\aligned
{\Ext}_{{D}^{+}({\Cal A};W)}^{j}((K,W), (L,W[m])) 
&= H^{j}W^{m}\Hom^{\ssbull }(K, L), \\
{\Ext}_{{D}^{+}({\Cal A};G,W)}^{j}((K; G,W), (L; G[k],W[m])) 
&= H^{j}G^{k}W^{m}\Hom^{\ssbull }(K, L),
\endaligned
\tag A.5.2
$$
where
$ {\Ext}_{{D}^{+}({\Cal A};W)}^{j}((K,W), (L,W[m])) = \Hom_{D^{+}
({\Cal A};W)}((K,W), (L[j],W[m]))  $, etc.
So we get spectral sequences
$$
\aligned
{E}_{1}^{p,q} = \prod _{j-i=p} {\Ext}_{{D}^{+}({\Cal A};W)}^{q}
&({\Gr}_{G}^{i}(K,W)[i], {\Gr}_{G}^{j}(L,W)[j]) \\
\Rightarrow  \,
&{\Ext}_{{D}^{+}({\Cal A};W)}^{p+q}((K,W), (L,W)), \\
{E}_{1}^{p,q} =  \text{(same as above if }
&p \ge  k,\,\,\, \text{and }0\,\,\, \text{otherwise)  } \\
\Rightarrow  \,
& {\Ext}_{{D}^{+}({\Cal A};G,W)}^{p+q}((K; G,W), (L; G[k],W)),
\endaligned
\tag A.5.3
$$
where  $ {E}_{1}^{p,q} = 0 $ for
$ q < 0 $.
Combining this with (A.6) below, we get the full faithfulness of 
(A.4.3).
(Here
$ k = 0 $ in the above spectral sequence.)
We can verify the essential surjectivity of (A.4.3) by an argument similar 
to the second solution in [2, 3.1.8].
\enddemo

\bigskip\noindent
{\bf A.6.~Lemma.} {\sl
With the notation and assumption of (A.4), let
$ (K,W), (L,W) \in  D^{b}({\Cal A};W) $ such that
$ {\Gr}_{W}^{k}K, {\Gr}_{W}^{k}L \in  {\Cal C} $ for any
$ k $.
Then
$ K, L \in  {\Cal C} $, and
$ W $ determines a filtration
$ \widetilde{W} $ on
$ K, L $ in
$ {\Cal C} $,
so that a natural morphism
$$
\Hom_{D^{b}({\Cal A};W)}((K,W),(L,W)) \rightarrow  \Hom
_{{\Cal C}}((K,\widetilde{W}),(L,\widetilde{W}))
\tag A.6.1
$$
is an isomorphism.
}
\bigskip

\demo\nofrills {Proof.\usualspace}
Let
$$
k = \min\{i \in  {\Bbb Z} : {\Gr}_{W}^{i}K \ne  0\}.
$$
Then we may assume
$ W^{k}L = L $ (replacing
$ L $ with
$ W^{k}L $ if necessary), because any filtered quasi-isomorphism
$ (K',W) \rightarrow  (K,W) $ may be replaced with its composition with
$ (W^{k}K',W) \rightarrow  (K',W) $.
We prove the bijectivity of (A.6.1) by induction on the length of the 
filtration
$ W $ on
$ K $.
Here we may assume that
$ {\Gr}_{W}^{k}L^{j} $ are injective objects of
$ {\Cal A} $.

Assume first
$ W $ on
$ K $ is trivial, i.e.,
$ {\Gr}_{W}^{i}K = 0 $ for  $ i \ne  k $.
In this case we have
$$
\Hom_{D^{b}({\Cal A};W)}((K,W),(L,W)) = \Hom_{D^{b}{\Cal A}}(K,L),
$$
because
$ D^{b}({\Cal A};W), D^{b}{\Cal A} $ can be replaced with
$ K^{b}({\Cal A};W), K^{b}{\Cal A} $ by assumption on
$ (L,W) $.
We have a similar isomorphism for
$ (K,\widetilde{W}),(L,\widetilde{W}) $,
and the assertion follows.

In general, let
$ K' = W^{k+1}K, K'' = {\Gr}_{W}^{k}K $ so that 
we have an exact sequence
$$
\aligned
0 &\rightarrow  \Hom_{D^{b}({\Cal A};W)}((K'',W),(L,W)) \rightarrow  
\Hom_{D^{b}({\Cal A};W)}((K,W),(L,W)) \\
&\rightarrow  \Hom_{D^{b}({\Cal A};W)}((K',W),(L,W)) \rightarrow  
{\Ext}_{{D}^{b}({\Cal A};W)}^{1}((K'',W),(L,W)),
\endaligned
$$
and the last term is isomorphic to
$ {\Ext}_{{D}^{b}{\Cal A}}^{1}(K'',L) $ by the same argument as above
(using the injectivity of
$ (L, W) $).

On the other hand, we have an exact sequence
$$
\aligned
0 &\rightarrow  \Hom_{{\Cal C}}((K'',\widetilde{W}),(L,\widetilde{W})) 
\rightarrow  \Hom_{{\Cal C}}((K,\widetilde{W}),(L,\widetilde{W})) \\
&\rightarrow  \Hom_{{\Cal C}}((K',\widetilde{W}),(L,\widetilde{W})) 
\rightarrow  
{\Ext}_{{D}^{b}({\Cal C};\widetilde{W})}^{1}
((K'',\widetilde{W}),(L,\widetilde{W})),
\endaligned
$$
and the last term is isomorphic to
$ {\Ext}_{C}^{1}(K'',L) $.
Since we have a natural morphism between the exact sequences, the 
assertion follows from the inductive hypothesis.
\enddemo

\bigskip\noindent
{\bf A.7.~Proposition.} {\sl
Let
$ {F}_{n}^{c}({\Cal A}) $ be as in (A.1).
Assume the filtrations
$ F_{1}, \dots , F_{n} $ are finite.
Let
$$
E = (M; F_{1}, \dots , F_{n}) \in  {F}_{n}^{c}({\Cal A})
$$
such that
$ {\Gr}_{{F}_{1}}^{{p}_{1}} \dots {\Gr}_{{F}_{n}}^{{p}_{n}}M $ are injective.
Then it is an injective object of
$ {F}_{n}^{c}({\Cal A}) $ (i.e.
 $ \Hom(E_{2},E) \rightarrow  \Hom(E_{1},E) $ is surjective for any strict 
monomorphism
$ E_{1} \rightarrow  E_{2} $ of
$ {F}_{n}^{c}({\Cal A})) $,
 $ M $ has a decomposition
$ \bigoplus _{\nu \in {\Bbb Z}^{n}} M_{\nu } $ such that
$ {F}_{1}^{{p}_{1}} \dots {F}_{n}^{{p}_{n}}M $ is the direct sum of
$ M_{\nu } $ for
$ \nu _{i} \ge  p_{i} $,
and
$$
({\Gr}_{{F}_{n}}^{p}\Hom(M',M), F_{i}  \,\, (i \ne  n)) = \prod _{j-i=p} 
(\Hom({\Gr}_{{F}_{n}}^{i}M', {\Gr}_{{F}_{n}}^{j}M), F_{i}  \,\, (i \ne  n))
$$
for
$ E' = (M'; F_{1}, \dots , F_{n}) \in  {F}_{n}^{c}({\Cal A}) $,
where
$ \Hom(M,M') $ has the filtration
$ F_{1}, \dots , F_{n} $ as usual.
}
\bigskip

\demo\nofrills {Proof.\usualspace}
The assertions follow by induction on
$ n $.
If the assertions hold for
$ n - 1 $,
we get the decomposition by induction on the length of
$ F_{n} $ using a splitting of the exact sequence
$$
0 \rightarrow  ({F}_{n}^{p+1}M, F_{i}  \,\, (i \ne  n)) \rightarrow  
({F}_{n}^{p}M, F_{i}  \,\, (i \ne  n)) \rightarrow  ({\Gr}_{{F}_{n}}^{p}M, F_{i}  
\,\, (i \ne  n)) \rightarrow  0.
$$
Then the third assertion follows, and then the first.
\enddemo

\bigskip\noindent
{\it Remark.}
The first assertion implies the isomorphism
$$
\Hom_{K^{+}} (E, E') = \Hom_{D^{+}} (E, E')
\tag A.7.1
$$
for complexes  
$ E = (K, F_{1}, \dots, F_{n}) $, 
$ E' = (K', F_{1}, \dots, F_{n}) $ of
$ F_{n}^{c}(\Cal A) $ such that
$ E $  is bounded,
$ E' $ is bounded below, and
$ \Gr_{F_{1}}^{p_{1}} \dots \Gr_{F_{n}}^{p_{n}}K^{\prime i} $ are injective.
(In fact, the left-hand side is zero if
$ E $ is strictly acyclic.)
Here
$ \Hom_{K^{+}} (E, E') $, 
$ \Hom_{D^{+}} (E, E') $ mean the group of morphisms in
$ K^{+}F_{n}^{c}(\Cal A) $,
$ D^{+}F_{n}^{c}(\Cal A) $ respectively.

\bigskip\bigskip\centerline{{\bf 3. Geometric Complexes}}

\bigskip\noindent
In this section, we define the mixed 
$ A $-Hodge complex
$ {\Cal C}^{H}_{X\langle D \rangle}A $ on a complex algebraic variety with a 
closed subvariety
$ D $

\bigskip
\noindent
{\bf 3.1.} Let
$ k $ be a field, and
$ {\Cal V}_{k} $ the category of varieties over
$ k $.
Let
$ A $ be a commutative ring.
We denote by
$ {\Cal V}^{A}_{k} $ the
$ A $-linear category consisting of varieties 
$ X $ over
$ k $ such that the group of morphisms
$ \Hom^{A}_{k}(X,Y) $ is the free
$ A $-module
$ \bigoplus_{f_{j,i}} A[f_{j,i}] $ with  basis the morphisms
$ f_{j,i} $ of varieties over
$ k $ from a connected component
$ X_{i} $ of
$ X $ to a connected component
$ Y_{j} $ of
$ Y $.
The composition is defined by
$ [g_{k,j'}]\circ [f_{j,i}] = [g_{k,j}f_{j,i}] $ if
$ j = j' $ and zero otherwise.
We have a natural functor
$ \iota^{A} : {\Cal V}_{k} \to {\Cal V}^{A}_{k} $ such that
$ \iota^{A}(X) = X $ and
$ \iota^{A}(f) = \sum [f_{j(i),i}] $ where the
$ f_{j(i),i} : X_{i} \to Y_{j(i)} $ denote the restrictions of
$ f : X \to Y $ to the connected components of
$ X $,
$ Y $.

Let
$ {\Cal B} $ be a category such that its objects form a finite set,
$ \Hom(a,b) $ are finite for
$ a, b \in {\Cal B} $, and
$ \Hom(b,b) = \{id\} $ for
$ b \in {\Cal B} $.
We assume that the objects of 
$ {\Cal B} $ permit an ordering such that
$ \Hom(a,b) = \emptyset $ unless
$ a \geqq b $.

A
$ k_{{\Cal B}} $-variety
$ X_{{\Cal B}} = \{ X_{b} \} $ is by definition a covariant functor
$ \gamma_{X} $ from
$ {\Cal B} $ to
$ {\Cal V}_{k} $, where
$ X_{b} = \gamma_{X}(b) $.
A proper
$ k_{{\Cal B}} $-variety over
$ X_{{\Cal B}} $ is a 
$ k_{{\Cal B}} $-variety 
$ Y_{{\Cal B}} $ with a functorial morphism
$ \pi : Y_{{\Cal B}} \rightarrow X_{{\Cal B}} $ such that
$ \pi_{b} : Y_{b} \to X_{b} $ is proper for any
$ b \in {\Cal B} $.
An
$ A $-linear morphism of proper
$ k_{{\Cal B}} $-varieties
$ Y_{{\Cal B}} \rightarrow Y'_{{\Cal B}} $ over
$ X_{{\Cal B}} $ is a functorial morphism
$ \iota^{A}\gamma_{Y} \to \iota^{A}\gamma_{Y'} $,
where we assume that the functorial morphisms
$ Y_{b} \to Y'_{b} $ in
$ {\Cal V}^{A}_{k} $ are given by 
$ A $-linear combinations of
$ [f_{b,i,j}] $ such that the
$ f_{b,i,j} $ commute with the restrictions of
$ \pi_{b} $ to the connected components of
$ Y_{b} $,
$ Y'_{b} $.
We denote by
$ {\Cal V}^{A}_{k}(X_{\Cal B}) $ the 
$ A $-linear category of proper
$ k_{{\Cal B}} $-varieties over
$ X_{{\Cal B}} $ with
$ A $-linear functorial morphisms.
The direct sum is given by the disjoint union, 
and the zero object by the empty set.

Let
$ C({\Cal V}^{A}_{k}(X_{\Cal B})) $ denote the category of complexes of
$ {\Cal V}^{A}_{k}(X_{\Cal B}) $,
and
$ K({\Cal V}^{A}_{k}(X_{\Cal B})) $ the category whose groups of morphisms are
divided  by the subgroups of morphisms which are homotopic to zero.
We have also
$ C^{*}({\Cal V}^{A}_{k}(X_{\Cal B})) $, 
$ K^{*}({\Cal V}^{A}_{k}(X_{\Cal B})) $ for
$ * = +, - , b $ as in [21].
Then
$ K({\Cal V}^{A}_{k}(X_{\Cal B})) $ and
$ K^{*}({\Cal V}^{A}_{k}(X_{\Cal B})) $ are triangulated categories by loc.~cit.

Let
$ D_{b} $ be closed subvarieties of
$ X_{b} $ such that the
$ X_{b} \backslash D_{b} $ are stable by 
$ \gamma_{X}(u) $ for  morphisms 
$ u $ of
$ {\Cal B} $.
We say that
$ Y_{{\Cal B}} $ is nonsingular over
$ X\langle D \rangle_{{\Cal B}} $ 
if the
$ Y_{b} $ are nonsingular and the
$ (\pi_{b} ^{-1}(D_{b}))_{\red} $ are locally either divisors 
with normal crossings or connected 
components, where
$ \pi_{b}  : Y_{b} \rightarrow  X_{b} $ is 
the morphism given by
$ \pi : Y_{\Cal B} \to X_{\Cal B} $.
We denote by
$ {\Cal V}^{A}_{k}(X\langle D \rangle_{\Cal B})^{ns} $ the full subcategory of
$ {\Cal V}^{A}_{k}(X_{\Cal B}) $ consisting of such objects.
Then we can define
$ C^{*}({\Cal V}^{A}_{k}(X\langle D \rangle_{\Cal B})^{ns}) $,
$ K^{*}({\Cal V}^{A}_{k}(X\langle D \rangle_{\Cal B})^{ns}) $ as above.

Assume 
$ k = {\Bbb C} $.
Let
$ C^{+}(X_{b}^{\an},A) $ denote the category of complexes of
$ A $-Modules which are bounded below on the associated analytic space
$ X_{b}^{\an} $.
For
$ b \in {\Cal B} $, we have a contravariant additive functor
$$
\bold{R}{\pi }_{b,*}^{A} : {\Cal V}^{A}_{k}(X_{\Cal B}) \rightarrow 
C^{+}(X_{b}^{\an},A) $$
by associating
$ (\pi_{b}) _{*}C^{\ssbull }(A_{Y_{b}}) $ to
$ Y_{b} $ where
$ \pi_{b}  : Y_{b} \rightarrow  X_{b} $ is as above, and
$ C^{\ssbull }(A_{Y_{b}}) $ is the canonical flasque 
resolution
of  Godement of the constant sheaf
$ A_{Y_{b}} $ on
$ Y_{b}^{\an} $.
It induces the contravariant additive functors
$$
\aligned
\bold{R}{\pi }_{b,*}^{A} 
&: C^{+}({\Cal V}^{A}_{k}(X_{\Cal B})) \rightarrow  C^{+}(X_{b}^{\an},A),\\
\bold{R}{\pi }_{b,*}^{A} 
&: K^{+}({\Cal V}^{A}_{k}(X_{\Cal B})) \rightarrow  K^{+}(X_{b}^{\an},A).
\endaligned
$$
The image of
$ Y_{{\Cal B}}^{\ssbull } \in  K^{+}({\Cal V}^{A}_{k}(X_{\Cal B})) $ by  
$ \bold{R}{\pi }_{b,*}^{A} $ is 
denoted by
$ \bold{R}(\pi_{b}) _{*}A_{Y_{b}^{\ssbull}} $.
We say that
$ Y_{{\Cal B}}^{\ssbull } \in  K^{+}({\Cal V}^{A}_{k}(X_{\Cal B})) $ is
$ A $-acyclic over
$ X\langle D \rangle_{{\Cal B}} $,
if the
$ \bold{R}(\pi_{b}) _{*}A_{Y_{b}^{\ssbull}} $ 
are acyclic over
$ X_{b}^{\an} \backslash  D_{b}^{\an} $ for
$ b \in {\Cal B} $,
and a morphism
$ Y_{{\Cal B}}^{\ssbull } \rightarrow  Y^{\prime \ssbull }_{{\Cal B}} $ of
$ {\Cal V}^{A}_{k}(X_{\Cal B}) $ is an
$ A $-quasi-isomorphism over
$ X\langle D \rangle_{{\Cal B}} $,
if its mapping cone is
$ A $-acyclic over
$ X\langle D \rangle_{{\Cal B}} $.
Then inverting
$ A $-quasi-isomorphisms over
$ X\langle D \rangle_{{\Cal B}} $ of
$ K^{*}({\Cal V}^{A}_{k}(X_{\Cal B})) $ and
$ K^{*}({\Cal V}^{A}_{k}(X\langle D \rangle_{\Cal B})^{ns}) $,
we get
$ D^{*}({\Cal V}^{A}_{k}(X\langle D \rangle_{\Cal B}))$,
$ D^{*}({\Cal V}^{A}_{k}(X\langle D \rangle_{\Cal B})^{ns}) $.

If
$ |{\Cal B} | = 1 $ (i.e., if 
$ {\Cal B} $ consists of one object), then
$ X_{{\Cal B}} $,
$ X\langle D \rangle_{{\Cal B}} $ will be denoted by
$ X $,
$ X\langle D \rangle $.

\bigskip

\noindent
{\it Remark.} Constructions similar to
$ C({\Cal V}^{A}_{k}(X_{\Cal B})) $ seem to be well-known in
the theory of motives,  where
$ [f_{j,i}] $ is replaced by a cycle on
$ Y \times  Y' $.
In the Hodge setting, however, we cannot replace
$ [f_{j,i}] $ with a cycle, because we would have to take further
a simplicial resolution of the cycle (if it is singular),
and the push-down (i.e. the Gysin) morphism of
mixed Hodge complexes is not easy  to
construct at the level of complexes in general.

\bigskip\noindent
{\bf 3.2.~Proposition.} {\sl For any
$ Y_{{\Cal B}}^{\ssbull } \in  C^{b}({\Cal V}^{A}_{k}(X_{\Cal B})) $,
we have an
$ A $-quasi-isomorphism
$ Y^{\prime \ssbull }_{{\Cal B}} \rightarrow  
Y^{\ssbull }_{{\Cal B}} $  over
$ X\langle D \rangle_{{\Cal B}} $ such that
$ Y^{\prime \ssbull }_{{\Cal B}} $ is nonsingular over
$ X\langle D \rangle_{{\Cal B}} $.
So we have an equivalence of categories
$ D^{b}({\Cal V}^{A}_{k}(X\langle D \rangle_{\Cal B})^{ns}) 
\simto  D^{b}({\Cal V}^{A}_{k}(X\langle D \rangle_{\Cal B})) $.
}
\bigskip

\demo\nofrills {Proof.\usualspace}
By the same argument as in [9, I, 2.6], we have
$ Y^{\prime \ssbull }_{{\Cal B}} \rightarrow  Y_{{\Cal B}}^{\ssbull } $ in
$ C^{b}({\Cal V}^{A}_{k}(X_{\Cal B})) $ such that
$ Y_{{\Cal B}}^{\prime \ssbull}  $ is nonsingular over
$ X\langle D \rangle_{{\Cal B}} $, and it is given by a 
desingularization
$ Y^{\prime \ssbull }_{{\Cal B}} \rightarrow  Y_{{\Cal B}}^{\ssbull } $ in a generalized sense
(i.e., it is proper surjective, the
$ Y^{\prime k}_{b} $ are nonsingular, 
and there exists dense open subvarieties of
$ Y_{b}^{k} $ which are isomorphic to their inverse images in
$ Y^{\prime k}_{b} $).
Then the
$ Y_{b}^{k} $ have closed subvarieties
$ Z_{b}^{k} $ such that 
$ Y_{b}^{k} \backslash Z_{b}^{k} $ are 
isomorphic to their inverse images in
$ Y^{\prime k}_{b} $, the
$ Z_{b}^{k} $ are stable by the morphisms of
$ {\Cal B} $ and also by the morphisms
$ f_{j,i} $ appearing in the differential of
$ Y^{\ssbull}_{{\Cal B}} $, and
$ \max \dim Z_{b}^{k} < \max \dim Y_{b}^{k} $.
Let
$ Z^{\prime k}_{b} $ be the inverse image of
$ Z^{k}_{b} $ in
$ Y^{\prime k}_{b} $.  Then the associated single complex of
$$
Z^{\prime \ssbull}_{{\Cal B}} \rightarrow Z^{\ssbull}_{{\Cal B}}
 \oplus Y^{\prime \ssbull}_{{\Cal B}} \rightarrow Y^{\ssbull}_{{\Cal B}}
\tag 3.2.1
$$
is
$ A $-acyclic, and we can proceed by induction on 
$ \dim Y^{\ssbull} $ as in loc.~cit.
The second assertion is clear by [21].
\enddemo

\bigskip
\noindent
{\bf 3.3.} Let
$  X $ be a complex algebraic variety, and
$ D $ a closed subvariety of
$ X $.
Assume
$ A $ is a subfield of
$ {\Bbb R} $.
Applying (3.2) to the case 
$ |{\Cal B} | = 1 $, there exists an
$ A $-quasi-isomorphism
$ Y^{\ssbull } \rightarrow  X $ over
$ X\langle D \rangle $ such that
$ Y^{\ssbull } \in  C^{b}({\Cal V}^{A}_{k}(X\langle D \rangle_{\Cal B})^{ns}) $.
We define an
$ A $-mixed Hodge complex
$$
{\Cal C}^{H}_{X\langle D \rangle}A = ((\tilde{\Omega}^{\ssbull} _{X\langle D 
\rangle};F,W), (\tilde{A}_{X\langle D \rangle},W), 
(\tilde{\Bbb C}_{X\langle D \rangle},W); \alpha _{F}, \alpha 
_{A})
$$
by
$$
\aligned
((\tilde{\Omega}^{\ssbull} _{X\langle D \rangle};F,W) 
&= \pi _{*}{C}_{\le r}^{\ssbull}(\Omega _{Y^{\ssbull}}(\log \pi ^{-1}D);F,W), 
\\
(\tilde{A}_{X\langle D \rangle},W) 
&= \pi _{*}(j_{*}j^{*}{C}_{\le r}^{\ssbull}(A_{Y^{\ssbull,\an}}),W), \\
(\tilde{\Bbb C}_{X\langle D \rangle},W) 
&= \pi _{*}(j_{*}j^{*}{C}_{\le r}^{\ssbull}(\Omega _{Y^{\ssbull,\an}}),W),
\endaligned
$$
for
$ r $ large enough.
Here
$ \pi  : Y^{\ssbull } \rightarrow  X $ and
$ j : \pi ^{-1}(X \backslash  D) \rightarrow  Y^{\ssbull } $ denote natural 
morphisms, and
$ {C}_{\le r}^{\ssbull} $ is the truncation
$ \tau _{\le r} $ of the canonical flasque resolution of Godement
(with respect to the Zariski or the classical topologies) and is applied
to each component of
$ \Omega_{Y^{\ssbull}}(\log \pi^{-1}D) $, etc.
(More precisely,
$ \tilde{\Omega}^{\ssbull} _{X\langle D \rangle} $ is the associated 
single complex of a triple complex, and the differential with respect to 
the index of
$ Y^{\ssbull } $ is induced by the differential of
$ Y^{\ssbull } $.)
The filtration
$ F $ on
$ \Omega _{Y^{i}}(\log \pi ^{-1}D) $ is the filtration ``b\^ete''
$ \sigma  $ of [4] (where
$ F_{p} = F^{-p}) $, and
 $ W_{k} $ on
$ \Omega _{Y^{i}}(\log \pi ^{-1}D), j_{*}j^{*}{C}_{\le 
r}^{\ssbull}(A_{Y^{i,\an}}) $ is defined by
$$
W_{k+i}\Omega _{Y^{i}}(\log \pi ^{-1}D),\quad \tau _{\le 
k+i}j_{*}j^{*}{C}_{\le r}^{\ssbull}(A_{Y^{i,\an}}),
$$
(similarly for
$ j_{*}j^{*}{C}_{\le r}^{\ssbull}(\Omega _{Y^{i,\an}})) $.
The morphisms
$ \alpha _{F}, \alpha _{A} $ are induced by  natural morphisms
$$
({C}_{\le r}^{\ssbull}(\Omega _{Y^{\ssbull}}(\log \pi ^{-1}D)))^{\an}
\rightarrow   
j_{*}j^{*}{C}_{\le r}^{\ssbull}(\Omega_{Y^{\ssbull,\an}}),\quad
A_{Y^{\ssbull,\an}}  \rightarrow  {\Bbb C}_{Y^{\ssbull,\an}} 
\rightarrow 
\Omega _{Y^{\ssbull,\an}}. 
$$ 

The condition of mixed Hodge complex is satisfied by using the stability 
by the direct images under proper morphisms (2.8).
We have the independence of the choice of
$ Y^{\ssbull } $ by (3.2).
The underlying
$ A $-complex
$ \tilde{A}_{X\langle D \rangle} $ is isomorphic to
$ \bold{R}j_{*}j^{*}A_{X^{\an}} $.

Applying (3.2) to the case 
$ {\Cal B} $ is an ordered set
$ \{0,1\} $, we see that
$ {\Cal C}^{H}_{X\langle D \rangle}A $ is functorial for 
$ (X, D) $, i.e., we have a natural morphism
$$
{\Cal C}^{H}_{X\langle D \rangle}A \rightarrow 
g_{*}{\Cal C}^{H}_{X'\langle D' \rangle}A 
\tag 3.3.1
$$
for a morphism 
$ g : X' \rightarrow X $ such that 
$ g^{-1}D \subset D' $, where
$ D' $ is a closed subvariety of
$ X' $.
Furthermore we get its (shifted) mapping cone
which is denoted by
$ {\Cal C}^{H}_{X\langle D \rangle \leftarrow X'\langle D' \rangle}A $.

In fact, we have a
$ k_{\Cal B} $-variety
$ X_{{\Cal B}} = \{ X_{b} \} $ such that
$ X_{0} = X $,
$ X_{1} = X' $,
where
$ {\Cal B} = \{0,1\} $ and
$ \Hom(1,0) $ consists of one element corresponding to
$ g $.
By (3.2) we have an
$ A $-quasi-isomorphism 
$ Y_{{\Cal B}}^{\ssbull} \rightarrow  X_{{\Cal B}} $ over
$ X\langle D \rangle_{{\Cal B}} $ such that the
$ Y_{{\Cal B}}^{i} $ are nonsingular over
$ X\langle D \rangle_{{\Cal B}} $.
So the morphism (3.3.1) and its shifted mapping cone are 
defined canonically by using
$ Y_{{\Cal B}}^{\ssbull} $ (where
the degree of the first term of the shifted mapping cone is zero).
By (3.2),
$ {\Cal C}^{H}_{X\langle D \rangle \leftarrow X'\langle D' \rangle}A $
is functorial for
$ (X,D) $,
$ (X',D') $.

If
$ X' $ is a closed subvariety 
$ Z $ of
$ X $ such that 
$ D' = Z \cap D $, then
$ {\Cal C}^{H}_{X\langle D \rangle\leftarrow X'\langle D' \rangle}A $,
$ \tilde{\Omega}^{\ssbull} _{X\langle D \rangle
\leftarrow X'\langle D' \rangle} $ 
are denoted by 
$ {\Cal C}^{H}_{X\langle D \rangle,Z}A $,
$ \tilde{\Omega}^{\ssbull} _{X\langle D\rangle,Z} $.
They are denoted by
$ {\Cal C}^{H}_{X\leftarrow X'}A $,
$ \tilde{\Omega}^{\ssbull} _{X\leftarrow X'} $ ,
if
$ D $,
$ D' $ are empty.

\bigskip

\noindent
{\it Remark.} By definition
$ (\tilde{\Omega}^{\ssbull} _{X\langle D \rangle},F) $ is essentially 
same as the filtered complex of  Du  Bois in [6].
See also [9].

\bigskip\noindent
{\bf 3.4.} {\it Example.}
Assume
$ X $ is smooth,
$ D $,
$ Z $ are divisors on
$ X $ such that
$ D \cup Z $ is a divisor with normal crossings,
$ \dim D \cap Z < \dim X $, and
the irreducible components
$ Z_{j} $ of
$ Z $ are smooth.
Let
$ X^{(0)} = X $, and
$ X^{(i)} $ be the disjoint union of 
$ \bigcap_{1 \le k \le i} Z_{j_{k}} $ for
$ j_{1} < \dots < j_{i} $.
We denote the natural inclusions by
$ a_{i} : X^{(i)} \to X $.
Put
$ D^{(i)} = a^{-1}_{i}(D) $.
Then the complex
$ X \leftarrow Z $ is
$ A $-quasi-isomorphic to the complex
$ X^{(0)} \leftarrow X^{(1)} \leftarrow X^{(2)} \leftarrow 
\cdots $ over
$ X \langle D \rangle $.
So
$ \tilde{\Omega}^{\ssbull}_{X\langle D \rangle,Z} $ is
represented by
$$
\Omega^{\ssbull}_{X^{(0)}}(\log D^{(0)}) \to
(a_{1})_{*}\Omega^{\ssbull}_{X^{(1)}}(\log D^{(1)}) \to
(a_{2})_{*}\Omega^{\ssbull}_{X^{(2)}}(\log D^{(2)}) \to
\cdots,
$$
where
$ F, W $ on
$ \Omega^{\ssbull}_{X^{(i)}}(\log D^{(i)}) $ are given by
$ \sigma $ and $ W[-i] $ in [4] (where 
$ (W[-i])_{k} = W_{i+k} $).
Let
$$
(j^{\Cal D}_{*}j^{*}j^{\prime \Cal D}_{!}j^{\prime *}{\Cal O}_{X},F,W)
$$
denote the underlying filtered
$ {\Cal D}_{X} $-Module of the mixed Hodge Module
$ j_{*}j^{*}j'_{!}j^{\prime *}A_{X}^{H} $, where
$ j : X\backslash D \to X $ and
$ j' : X\backslash Z \to X $ denote natural inclusions.
Then we have a natural isomorphism in
$ D^{b}_{\hol}F(X,{\Cal D};W) $ :
$$
 (j^{\Cal D}_{*}j^{*}j^{\prime \Cal D}_{!}j^{\prime *}{\Cal O}_{X},F,W) 
= \DR_{X}^{-1}
(\tilde{\Omega}^{\ssbull}_{X\langle D\rangle,Z},F,W)[\dim X] .
\tag 3.4.1
$$
(Here the filtration
$ W $ on the right-hand side is shifted by
$ \dim X $ due to the shift of complex by
$ \dim X $, and
$ W' = W $ on the left-hand side, because
$ j^{\Cal D}_{*}j^{*}j^{\prime \Cal D}_{!}j^{\prime *}{\Cal O}_{X} $ is a
$ {\Cal D}_{X} $-Module.)

For the proof of (3.4.1), we may assume
$ X $ is a complex manifold by GAGA.
Then we have locally
$ X = X_{1}\times X_{2} $ such that
$ D = D_{1}\times X_{2} $,
$ Z = X_{1}\times D_{2} $, where
$ D_{1} $,
$ D_{2} $ are divisors with normal crossings on
$ X_{1} $,
$ X_{2} $.
We have a bifiltered isomorphism
$$
(j^{\Cal D}_{*}j^{*}j^{\prime \Cal D}_{!}j^{\prime *}{\Cal O}_{X},F,W) =
 ( (j_{1})^{\Cal D}_{*}j_{1}^{*}{\Cal O}_{X_{1}},F,W) \boxtimes
( (j_{2})^{\Cal D}_{!}j_{2}^{*}{\Cal O}_{X_{2}},F,W) ,
$$
using the functorial morphisms
$ id \to j_{*}j^{*} $,
$ j'_{!}j^{\prime *} \to id $.
We have a similar decomposition on the right-hand side
of (3.4.1).
So the assertion is reduced to the case either
$ D $ or
$ Z $ is empty.
Then the assertion follows from (v) of (1.4) forgetting 
the filtration
$ W $.
But the compatibility with
$ W $ is easy, because the weight filtration 
$ W $ of
$ j'_{!}j^{\prime *}{\Bbb C}_{X}[\dim X] $ in the abelian category of 
perverse sheaves is uniquely characterized
by the property :
$$
{\Gr}_{\dim X - k}^{W} ( j'_{!}j^{\prime *}{\Bbb C}_{X}[\dim X]) =
(a_{k})_{*} {\Bbb C}_{X^{(k)}}[\dim X - k] .
$$

\bigskip\bigskip\centerline{{\bf 4. Smooth Affine Stratification}}

\bigskip\noindent
We prove (0.2) by using a good representative of the mixed Hodge
Module
$ A_{X}^{H} $ associated with a smooth affine stratification of
$ X $.

\bigskip\noindent
{\bf 4.1.~Proposition.} {\sl Let  
$ X $ be a complex algebraic variety
of dimension 
$ n $, and
$ {\Cal S} = \{S_{i}\} $ a stratification of
$ X $ such that each stratum
$ S_{i} $ is smooth and affine.
Let
$ j_{i} : S_{i} \rightarrow  X $ denote  natural inclusions.
Then we have a complex
$ {C}_{\Cal S}^{\ssbull}A^{H} $ in
$ C^{b}\MHM(X,A) $ such that
$$
{C}_{\Cal S}^{k}A^{H} = \bigoplus _{\dim S_{i} = k} 
\,(j_{i})_{!}{A}_{{S}_{i}}^{H}[k],
\tag 4.1.1
$$
and the image 
of
$ {C}_{\Cal S}^{\ssbull}A^{H} $ in
$ D^{b}\MHM(X,A) $ is naturally isomorphic to
$ {A}_{X}^{H} $.
}
\bigskip

\demo\nofrills {Proof.\usualspace}
We first note that (4.1.1) implies for any
$ M \in \MHM(X,A) $, a natural morphism
$  C_{\Cal S}^{n}A^{H}  \rightarrow
C_{\Cal S}^{\ssbull}A^{H}[n] $ induces an isomorphism
$$
\Hom (C_{\Cal S}^{\ssbull} A^{H}[n], M) =
\Ker(\Hom (C_{\Cal S}^{n} A^{H}, M) \rightarrow
\Hom (C_{\Cal S}^{n-1} A^{H}, M)).
$$
Here 
$ \Hom $ on the left-hand (resp. right-hand) side is taken in
$ D^{b}\MHM(X,A) $
(resp. 
$ \MHM(X,A) $), and the morphism in the right-hand side is induced by 
the differential of
$ C_{\Cal S}^{\ssbull}A^{H} $.

We proceed by induction on 
$ n = \dim X $.
 Let
$ Y = \bigcup _{\dim S_{i} < n} S_{i} $,
and
$ U = X \backslash  Y $ with natural inclusions
$ i : Y \rightarrow  X, j : U \rightarrow  X $.
Then we have a distinguished triangle in
$ D^{b}\MHM(X,A) $
$$
\rightarrow  i_{*}{A}_{Y}^{H}[n-1] \rightarrow  j_{!}{A}_{U}^{H}[n] 
\rightarrow  {A}_{X}^{H}[n] \rightarrow ,
\tag 4.1.2
$$
and the assertion follows from the inductive hypothesis.
Note that the isomorphism
$ {A}_{X}^{H} = {C}_{\Cal S}^{\ssbull}A^{H} $ is uniquely determined by the 
restriction to a dense open subvariety of
$ X $,
because
$$
\End({A}_{X}^{H}) = \Hom({A}_{pt}^{H}, (a_{X})_{*}{A}_{X}^{H}) = A,
$$
if
$ X $ is connected (where
$ a_{X} : X \rightarrow  pt $ is a natural morphism).
\enddemo

\bigskip\noindent
{\bf 4.2.~Theorem.} {\sl Let  
$ X $ be a complex algebraic variety, and
$ D $,
$ D' $  closed subvarieties with natural inclusions
$ j : X \backslash  D \rightarrow  X $,
$ j' : X \backslash  D' \rightarrow  X $.
Then we have a natural isomorphism
$$
{\DR}_{X}^{-1}({\Cal C}^{H}_{X\langle D \rangle , D'}A) = \varepsilon 
(j_{*}j^{*}j'_{!}j^{\prime *}{A}_{X}^{H})\quad
\text{in }{D}_{\Cal H}^{b}(X,A)_{{\Cal D}}.
\tag 4.2.1
$$
}

\demo\nofrills {Proof.\usualspace}
By (2.8) the assertion is reduced to the case where
$ D, D' $ are locally principal divisors (using the blow-up of
$ X $ along
$ D $ and
$ D' $) so that
$ j_{*}, j'_{!} $ preserve perverse sheaves.
Let
$ \{S_{i}\} $ be as in (4.1).
Here we may assume that
$ D' $ is a union of strata.
Let
$$
X^{j} = \bigcup _{\dim S_{i} \le  j}  S_{i} \bigcup D' 
$$
for
$ -1 \le j \le n $.  By the same argument as (3.2) we have 
$ A $-quasi-isomorphisms
$ Y^{j,\ssbull } \rightarrow  X^{j} $ over
$ X^{j}\langle X^{j}\cap D\rangle $ together 
with morphisms
$ Y^{j-1,\ssbull } \rightarrow  Y^{j,\ssbull } $  such that
$ Y^{j,\ssbull }  $ is nonsingular over
$ X\langle D \rangle $,
and we have a 
commutative diagram in
$ C^{b}({\Cal V}(X\langle D \rangle,A))  $:
$$
\CD
Y^{j,\ssbull } @<<<    Y^{j-1,\ssbull } \\
@VVV @VVV \\
X^{j}   @<<<    X^{j-1}.
\endCD
$$

Using
$ Y^{j,\ssbull } $ as in (3.3), we define
$$ 
K^{j} = {\DR}_{X}^{-1}(i_{j})_{*}{\Cal C}^{H}_{X^{j}\langle D \cap X^{j}\rangle}A 
\quad\text{in }
 {C}_{\Cal H}^{b}(X,A)_{{\Cal D}}
$$
so that we have natural morphisms
$ r_{j} : K^{j} \rightarrow  K^{j-1} $ in
$ {C}_{\Cal H}^{b}(X,A)_{{\Cal D}} $.
Here
$ i_{j} : X^{j} \rightarrow X $
denotes a natural inclusion.
We consider a morphism of
$ {C}_{\Cal H}^{b}(X,A)_{{\Cal D}}  $:
$$
\bigoplus _{0\le j\le n} K^{j} \rightarrow  \bigoplus _{-1\le j<n} K^{j}
\tag 4.2.2
$$
whose restriction to
$ K^{j} $ is
$ r_{n} $ for
$ j = n $, and
$ r_{j}\oplus (-id) $ for
$ 0 \le j < n $,
where
$ n = \dim X $.
Let
$ K $ be its shifted mapping cone (i.e., the degree of the first term is 
zero).
Then
$ K $ is isomorphic to the shifted mapping cone of
$ K^{n} \rightarrow K^{-1}$ so that
$$
K = \DR_{X}^{-1}({\Cal C}^{H}_{X\langle D \rangle ,D'}A).
$$

We have a filtration
$ G $ on
$ K $ such that
$ G^{i} $ is the shifted mapping cone of
$$
\bigoplus _{i\le j\le n} K^{j} \rightarrow  \bigoplus _{i-1\le j<n} K^{j}
$$
for
$ i \ge 0 $.
Let
$ {\Cal S}' = \{ S_{i} : S_{i} \subset X \backslash D' \} $
so that
$ {\Cal S}' $ is a stratification of
$ X \backslash D' $.
Then
$ {\Gr}_{G}^{i}K $ is the shifted mapping cone of
$ r_{i} : K^{i} \rightarrow  K^{i-1} $ which corresponds to
$ j_{*}j^{*}j'_{!}{C}_{\Cal S'}^{i}A^{H}[-i] $.
So
$ G $ is similar to the filtration
$ G $ in the proof of (2.8), and we have weak quasi-isomorphisms like 
(2.8.3).  Furthermore
$$
{\Gr}_{G}^{i}((\Dec G)^{0}K/(\Decdual G)^{1}K))[i]
$$
is isomorphic to
$ \varepsilon (j_{*}j^{*}j'_{!}{C}_{\Cal S'}^{i}A^{H})  $.
See(2.9).
In fact, the last assertion is reduced to the case 
where the closure 
$ \overline{S}_{i} $ of
$ S_{i} $ is smooth and
$ \overline{S}_{i}\backslash S_{i} $ is a divisor with normal crossings
by using (2.8). 
Then it follows from (3.4).
This finishes the proof of (4.2).
\enddemo

\bigskip\noindent
{\bf 4.3.~Corollary.} {\sl
With the notation of (4.2), the two mixed Hodge structures on
$ H^{i}(X \backslash D, D'\backslash D' \cap D) $
by [4] and [15] coincide.
}

 \bigskip\bigskip\centerline{{\bf 5. Du Bois Singularity}}

\bigskip\noindent
We study Du Bois singularity, and give some application of (0.3).
 
\bigskip
\noindent
{\bf 5.1.} Let
$ X $ be a reduced complex algebraic variety.
Let
$ \pi' : X' \to X $ be the normalization of
$ X $.
Put
$ X'' = (X' \times_{X} X')_{\red} $ with the natural projections
$ \pi'' : X'' \to X $ and
$ p_{a} : X'' \to X' $
($ a = 1, 2 $).
Let
$$
{\Cal O}_{X}^{wn} = \Ker(p_{1}^{*} - p_{2}^{*} : 
\pi'_{*}{\Cal O}_{X'} \to \pi''_{*}{\Cal O}_{X''} ).
\tag 5.1.1
$$
Then
$ {\Cal O}^{wn}_{X} $ is a coherent sheaf of algebras,
and we define
$ X^{wn} = \hbox{\rm Spec}_{X}{\Cal O}^{wn}_{X} $.
The natural morphism
$ X^{wn} \to X $ induces a bijection of the underlying sets,
because
$ {\Cal O}^{wn}_{X} $ is identified with the sheaf of continuous
functions on
$ X({\Bbb C}) $ whose pull-backs to
$ X' $ come from
$ {\Cal O}_{X'} $.
We call
$ {\Cal O}_{X}^{wn} $,
$ X^{wn} $  the {\it weak normalization} of
$ {\Cal O}_{X} $,
$ X $.
We say that
$ X $ is {\it weakly normal} if
$  {\Cal O}^{wn}_{X} =  {\Cal O}_{X} $.

\bigskip\noindent
{\it Remarks.} (i)
For
$ x \in X $, the analytic-local irreducible components of
$ X $ at 
$ x $ corresponds to
$ \pi^{\prime -1}(x) $.
In particular,
the weak normalization coincides with the normalization if
$ X $ is analytic-locally irreducible at every point.

(ii) Let
$ X_{1}, \dots, X_{r} $ be reduced closed subvarieties of
$ X $.
Put
$ X_{i,j} = (X_{i} \cap X_{j})_{\red} $ for
$ i < j $.
We say that
$ X_{1}, \dots, X_{r} $ are {\it relatively vertical},
if
$$
{\Cal O}_{X}\to \bigoplus_{i}  {\Cal O}_{X_{i}}  \to \bigoplus_{i<j}
 {\Cal O}_{X_{i,j}} 
\tag 5.1.2
$$
is exact, where the last morphism is given by the Cech morphism.
We see that if
$ X_{1}, \dots, X_{r-1} $ are relatively vertical, then
$ X_{1}, \dots, X_{r} $ are relatively vertical if and only if
$ \bigcup_{1 \le i < r} X_{i} $ and
$ X_{r} $ are relatively vertical.

(iii) Let
$ X_{0} = \bigcup_{1 \le i \le r} X_{i} $.
Then we see :

(a) If
$ X_{0} $ is weakly normal, 
$ X_{1}, \dots, X_{r} $ are relatively vertical.

\noindent
Assume
$ r = 2 $. 
Then :

(b) If
$ X_{0} $,
$ X_{1,2} $ are weakly normal,
so are
$ X_{1} $,
$ X_{2} $.

(c) If
$ X_{1} $,
$ X_{2} $,
$ X_{1,2} $ are weakly normal and 
$ X_{1} $,
$ X_{2} $ are relatively vertical, then
$ X_{0} $ is weakly normal.

(iv) In the isolated singularity case,
(iii) implies that a variety with isolated singularities is weakly normal, 
if and only if so are any  \'etale-local irreducible components and 
these are relatively vertical.
(Note that the normalization and the weak normalization are compatible
with \'etale base changes.)
In the one-dimensional case, this means that
$ X $ is weakly normal if and only if it is \'etale-locally isomorphic to
$ \bigcup _{j} \{x_{i} = 0 \,(i \ne j) \} \subset {\Bbb C}^{n} $.

\bigskip\noindent
{\bf 5.2.~Proposition.}  {\sl Let
$ X $ be a reduced complex algebraic variety, and
$ (M,F) $  the underlying filtered
$ {\Cal D} $-Module of the mixed Hodge Module
$ A^{H}_{X} $ in [15].
Then with the notation of (1.3), we have a natural isomorphism
$$
{\Cal H}^{0}{\Gr}_{0}^{F}\DR_{X}(M) = {\Cal O}_{X}^{wn}.
$$
}

\demo\nofrills {Proof.\usualspace}
Let
$ a_{\ssbull} : X_{\ssbull} \to X $ be a simplicial resolution (or
a cubic resolution).
By (0.3) it is enough to show
$$
{\Cal H}^{0} {\bold R}(a_{\ssbull})_{*}{\Cal O}_{X_{\ssbull}}
 = {\Cal O}_{X}^{wn}.
\tag 5.2.1
$$
We have the spectral sequence
$$
E_{1}^{p,q} = R^{q}(a_{p})_{*}{\Cal O}_{X_{p}} 
\Rightarrow 
{\Cal H}^{p+q} {\bold R}(a_{\ssbull})_{*}{\Cal O}_{X_{\ssbull}},
$$
which implies
$ {\Cal H}^{0} {\bold R}(a_{\ssbull})_{*}{\Cal O}_{X_{\ssbull}}
= \Ker((a_{0})_{*}{\Cal O}_{X_{0}} \to (a_{1})_{*}{\Cal O}_{X_{1}}) $.
So we get a natural morphism
$$
{\Cal H}^{0} {\bold R}(a_{\ssbull})_{*}{\Cal O}_{X_{\ssbull}}
\to \pi'_{*}{\Cal O}_{X'}
\tag 5.2.2
$$
in the notation of (3.1), because
$ \pi'_{*}{\Cal O}_{X'} = \pi_{*}{\Cal O}_{\tilde{X}} $ if
$ \pi : \tilde{X} \to X $ is a desingularization.
We see that (5.2.2) is injective and the image is 
$ {\Cal O}^{wn}_{X} $ by restricting
$ a_{\ssbull} : X_{\ssbull} \to X $ over a point of
$ X $, because 
$ a_{\ssbull} : X_{\ssbull} \to X $ gives a resolution of a constant sheaf.
So we get the assertion.
\enddemo

\bigskip\noindent
{\it Remarks.} (i)
By (0.3) and (5.2)
$ X $ is Du Bois (see the introduction for the definition), if and only if 
$ X $ is weakly normal and
$ {\Cal H}^{i} {\bold R}(a_{\ssbull})_{*}{\Cal O}_{X_{\ssbull}} = 0 $
(or, equivalently,
$ {\Cal H}^{i}{\Gr}_{0}^{F}\DR_{X}(M) = 0 $) for
$ i > 0 $.
Note that
$ X $ has (at most) rational singularities, if and only if
$ X $ is normal and
$ R^{i}\pi_{*}{\Cal O}_{\tilde{X}} = 0 $ for
$ i > 0 $, where
$ \pi : \tilde{X} \to X $ is a desingularization.

(ii) Assume
$ X $ is a hypersurface.
Let
$ f $ be a local reduced defining equation of
$ X $,
and
$ b_{f}(s) $ the
$ b $-function (i.e.
the Bernstein polynomial) of
$ f $.
Then
$ X $ has rational (resp.
Du Bois) singularity if and only if the roots of
$ b_{f}(s)/(s+1) $ are less than (resp.
less than or equal to)
$ -1 $.
See [17, (0.4)]  and [18, (0.5)].

(iii) With the notation of Remark (iii) after (5.1) assume
$ r = 2 $.
Then we can easily verify :

(a) If
$ X_{0} $,
$ X_{1,2} $ are Du Bois,
so are
$ X_{1} $,
$ X_{2} $.

(b) If
$ X_{1} $,
$ X_{2} $,
$ X_{1,2} $ are Du Bois and 
$ X_{1} $,
$ X_{2} $ are relatively vertical, then
$ X_{0} $ is Du Bois.

In the isolated singularity case,
these imply that a variety with isolated singularities is Du Bois,
if and only if so are any  \'etale-local irreducible components and 
these are relatively vertical.
In the one-dimensional case, this means that
$ X $ is Du Bois if and only if 
$ X $ is weakly normal.
See Remark (iv) after (5.1).

\bigskip\noindent
{\bf 5.3.~Proposition.} {\sl Let
$ ({\Bbb D}M,F) $ denote the underlying filtered
$ {\Cal D} $-Module of
$ {\Bbb D}A_{X}^{H} $ (the dual of
$ A_{X}^{H} $), and
$ K_{X} $ the dualizing complex for 
$ {\Cal O} $-Modules.
Then we have unique morphisms in
$ D^{b}(X,{\Cal O}_{X} ) $ :
$$
{\Cal O}_{X} \to {\Gr}_{0}^{F}\DR_{X}(M), \quad
{\Gr}_{0}^{F}\DR_{X}({\Bbb D}M) \to K_{X},
\tag 5.3.1
$$
whose restrictions to the smooth part of
$ X $ are natural isomorphisms.
}
\bigskip

\demo\nofrills {Proof.\usualspace}
 The assertion on the first morphism is clear by (5.2).
For the second, we have by the duality (see [14, \S 2]) :
$$
{\Bbb D}({\Gr}_{0}^{F}\DR_{X}(M)) =
{\Gr}_{0}^{F}\DR_{X}({\Bbb D}M)
\quad\text{in } D^{b}({\Cal O}_{X}),
\tag 5.3.2
$$
if
$ X $ is embeddable into a smooth variety.
Let
$ X = \bigcup_{1 \le i \le r} U_{i} $ be an open covering of
$ X $ such that
$ U_{i} $ are embeddable into smooth varieties.
Put
$ U_{I} = \bigcap_{i \in I} U_{i} $ for
$ I \subset \{1, \dots , r\} $.
Let
$ M' = {\Gr}_{0}^{F}\DR_{X}({\Bbb D}M) $.
Then we have a spectral sequence
$$
E_{1}^{p,q} = \bigoplus _{|I|=p+1} \Hom_{D^{+}(U_{I})}
(M'|_{U_{I}},K_{U_{I}}[q]) \Rightarrow
\Hom_{D^{+}(X)}(M',K_{X}[p+q]),
$$
taking an injective resolution of
$ K_{X} $ and using the Cech complex associated with
$ \{ U_{i} \} $.
Since
$ E_{1}^{p,q} = 0 $ for
$ p < 0 $ by the duality together with (0.3) and (5.3.2),
the second morphism of (5.3.1) is locally defined.
So the assertion follows from  (5.3.2).

\bigskip\noindent
{\t Remarks.} (i) One of the morphisms of (5.3.1) is an isomorphism
if and only if 
$ X $ is Du Bois.

(ii) If 
$ X $ is Du Bois and
$ {\Bbb Q}_{X}[\dim X] $ is a perverse sheaf , then
$ X $ is Cohen-Macaulay.
This is because the filtration
$ F $ on
$ ({\Bbb D}M)_{U\to V} $ is strict.
See (iii) of (1.4).

\bigskip\noindent
{\bf 5.4.~Theorem.} {\sl A rational singularity is Du Bois
}
\bigskip

\demo\nofrills {Proof.\usualspace}
Assume
$ X $ has (at most) rational singularities.
Let
$ \pi : \tilde{X} \to X $ be a desingularization.
By (0.3),
we have a natural morphism
$$
{\Gr}_{0}^{F}\DR_{X}(M) \rightarrow  \bold{R}\pi _{*}
{\Cal O}_{\tilde{X}} 
$$
whose composition with the first morphism
of (5.3.1) is a natural morphism
$ {\Cal O}_{X} \to \bold{R}\pi _{*}{\Cal O}_{\tilde{X}} $.
So
$ {\Cal O}_{X} $ is a direct factor of
$ {\Gr}_{0}^{F}\DR_{X}(M) $ because
$ X $ has only rational singularities.
Assume the singularity of
$ X $ is not Du Bois.
Then
$ {\Gr}_{0}^{F}\DR_{X}(M) $ has a nontrivial 
direct factor in
$ D^{b}({\Cal O}_{X}) $ whose (cohomological) support is contained in 
$ \Sing X $.

Since the assertion is local, we may assume that
$ X $ is a closed subvariety of a smooth variety
$ Y $.
Let
$ ({\Bbb D}M,F ) $ be as in (5.3).
To simplify the notation, we denote also by
$ (M,F) $,
$ ({\Bbb D}M,F ) $ the filtered
$ {\Cal D}_{Y} $-Modules representing them for  
$ \{ X \to Y \} $.
Then
$ ({\Bbb D}M, F) $ is really the dual of
$ (M,F) $ as filtered
$ {\Cal D}_{Y} $-Module, and (5.3.2) becomes
$$
{\Bbb D}{\Gr}_{0}^{F}\DR_{Y}(M) =
{\Gr}_{0}^{F}\DR_{Y}({\Bbb D}M) = F_{0}({\Bbb D}M)
\quad\text{in } D^{b}({\Cal O}_{Y}),
\tag 5.4.1
$$
where the last isomorphism follows from
$ {\Gr}^{F}_{p}({\Bbb D}M) = {\Gr}^{F}_{p}\DR_{Y}({\Bbb D}M) = 0 $ 
for
$ p < 0 $.
(The last vanishing follows from (0.3) and (5.3.2).)

So the assumption implies that
$ F_{0}({\Bbb D}M) $ has a nontrivial direct factor
$ N $ in
$ D^{b}({\Cal O}_{Y}) $ whose (cohomological) support is contained in 
$ \Sing X $.
Then we have a morphism of
$ {\Cal D}_{Y} $-Modules
$$
N\otimes _{{\Cal O}_{Y}}{\Cal D}_{Y} \rightarrow  {\Bbb D}M
\tag 5.4.2
$$
such that
$ {\Cal H}^{i}N\otimes _{{\Cal O}_{Y}}{\Cal D}_{Y} \rightarrow  
{\Cal H}^{i}({\Bbb D}M) $ is nontrivial if
$ {\Cal H}^{i}N \ne  0 $,
because
$ {\Cal H}^{i}F_{0}({\Bbb D}M) \rightarrow  {\Cal H}^{i}({\Bbb D}M) $
 is injective by the 
strictness of the Hodge filtration.
See (iii) of (1.4).

Let
$ Z = \bigcup _{i}  \Supp  {\Cal H}^{i}N $.
We may assume that
$ Z $ is smooth with pure dimension
$ m < \dim X $,
and also $ \Supp {\Cal H}^{i}N = Z $ or
$ \emptyset  $ by deleting a proper closed subvariety of
$ Z $ from
$ X $.
Let
$$
M' = \bold{R}\Gamma _{Z}{\Bbb D}M.
$$
Then (5.4.2) is the composition of  $ N\otimes _{{\Cal O}_{Y}}{\Cal D}_{Y} 
\rightarrow  M' \rightarrow  {\Bbb D}M  $.
Furthermore
$ M' $ underlies the dual of
$ {A}_{Z}^{H} $ so that
$ {\Cal H}^{i}M' = 0 $ for
$ i \ne  - m $.
So we get
$$
{\Cal H}^{-m}({\Bbb D}M) \ne  0
$$
by the property after (5.4.2).
Then, passing to the corresponding complex of
$ {\Bbb C}_{Y^{\an}} $-Modules via the functor
$ \DR_{Y} $ (and using the duality), we get
$$
{}^{p}{\Cal H}^{m}{\Bbb C}_{X^{\an}} \ne  0.
$$
But this contradicts to Remark below by restricting to a general smooth 
subvariety transversal to
$ Z $ (where the shift of index by
$ m $ comes from
$ {}^{p}{\Cal H}^{m}L = L[m] $ for a local system
$ L $ on
$ Z^{\an}) $.
\enddemo

\bigskip
\noindent
{\it Remark.} For a connected complex algebraic variety of
$ \dim X \ge  1 $,
we have
$$
{}^{p}{\Cal H}^{0}{\Bbb C}_{X^{\an}} = 0.
$$
In fact, let
$ D $ be a closed subvariety of
$ X $ such that
$ U = X \backslash  D $ is smooth affine and pure dimensional.
Let
$ j : U \rightarrow  X $ denote a natural inclusion.
Then we have a distinguished triangle in
$ {D}_{c}^{b}(X^{\an},{\Bbb C}) $
$$
\rightarrow  j_{!}{\Bbb C}_{U^{\an}} \rightarrow  {\Bbb C}_{X^{\an}} 
\rightarrow  {\Bbb C}_{Z^{\an}} \rightarrow 
$$
with
$ {}^{p}{\Cal H}^{i}j_{!}{\Bbb C}_{U} = 0 $ for
$ i \ne  \dim U $.
So the assertion is reduced to the case
$ \dim X = 1 $ by induction.
Then
$ {\Bbb C}_{X^{\an}}[1] \in   \Perv(X^{\an},{\Bbb C}) $ by the exact 
sequence
$$
0 \rightarrow  {\Bbb C}_{X^{\an}} \rightarrow  \pi _{*}
{\Bbb C}_{\tilde{X}^{\an}} \rightarrow  E \rightarrow  0,
$$
where
$ \pi  : \tilde{X} \rightarrow  X $ is the normalization, and 
$ \Supp E \subset  \Sing X $.
(The assertion follows also from
$ H_{\{0\}}^{i}{\Bbb C}_{X^{\an}} = 0 $ for
$ i \le 0 $.)

\bigskip\bigskip
\centerline{{\bf Table of notations}}

\bigskip

\item{1.1}
$ \Diff_{X}^{n}(L,L') $, 
$ \Hom_{\Diff}((L,F),(L',F)) $,
$ MF({\Cal D}_{X}) $,

\item{}
$ MF({\Cal O}_{X}, \Diff) $,
$ C^{*}F({\Cal O}_{X}, \Diff) $,
$ K^{*}F({\Cal O}_{X}, \Diff) $,
$ D^{*}F({\Cal O}_{X}, \Diff) $,

\item{}
$ MF({\Cal O}_{X}, \Diff;W) $,
$ C^{*}F({\Cal O}_{X}, \Diff;W) $,
$ K^{*}F({\Cal O}_{X}, \Diff;W) $,

\item{}
$ D^{*}F({\Cal O}_{X}, \Diff;W) $,
$ D^{*}_{\coh}F({\Cal O}_{X}, \Diff;W) $,
$ D^{*}_{\coh}F({\Cal O}_{X}, \Diff) $,

\item{}
$ \DR_{X}^{-1} : MF({\Cal O}_{X}, \Diff) \to MF({\Cal D}_{X}) $,

\item{1.2}
$ \bold{R}f_{*} : D^{+}F({\Cal O}_{X}, \Diff)  \rightarrow  D^{+}F
({\Cal O}_{Y}, \Diff)  $,

\item{1.3}
$ LE(X) $,
$ MF({\Cal D}_{V})_{U} $,
$ MF(X,{\Cal D}) $,
$ MF(X,{\Cal D};W) $,

\item{}
$ C^{*}F(X,{\Cal D}) $, 
$ C^{*}F(X,{\Cal D};W) $,
$ D^{*}F(X,{\Cal D}) $,
$ D^{*}F(X,{\Cal D};W) $,

\item{}
$ M(X,{\Cal D}) $,
$ D^{*}(X,{\Cal D}) $,
$ D^{*}(X,{\Cal D};W) $,

\item{}
$ D^{*}_{\coh} F(X,{\Cal D};W) $,
$ D^{*}_{\hol} F(X,{\Cal D};W) $,
$ D^{*}_{\coh} F(X,{\Cal D}) $,
$ D^{*}_{\hol} F(X,{\Cal D}) $,

\item{}
$ M_{\hol} (X,{\Cal D}) $,
$ D^{*}_{\hol} (X,{\Cal D}) $,
$ D^{*}_{\hol} (X,{\Cal D};W) $,

\item{}
$  D^{*}_{c} (X,{\Bbb C}) $,
$ \Perv (X,{\Bbb C}) $,

\item{}
$ \DR_{X} :D^{*}_{\hol} (X,{\Cal D}) \rightarrow 
 D^{*}_{c} (X,{\Bbb C}) $,
$ \DR_{X} : M_{\hol} (X,{\Cal D}) \rightarrow 
\Perv (X,{\Bbb C}) $,

\item{}
$ {\Gr}_{p}^{F}\DR_{X} :D^{*}_{\coh} F(X,{\Cal D}) \rightarrow 
 D^{*}_{\coh} (X,{\Cal O}_{X}) $,

\item{1.4}
$ MF_{\hol}(X,{\Cal D},A;W) $,
$ MF_{\hol}(X,{\Cal D},A) $,

\item{}
$ \MH(X,A,n) $,
$ \MHW(X,A) $,
$ \MHM(X,A) $,

\item{2.1}
$ {C}_{\hol}^{b}F({\Cal O}_{X},\Diff,A;W) $,
$ {C}_{\hol}^{b}F(X,{\Cal D},A;W)  $,
$ {C}_{\hol}^{b}F({\Cal O}_{X},\Diff,A) $,

\item{}
$ {C}_{\hol}^{b}F(X,{\Cal D},A) $,
$ \Perv(X,{\Bbb C};W) $,

\item{}
$ \DR_{X}^{-1} : {C}_{\hol}^{b}F({\Cal O}_{X},\Diff,A;W)  
\rightarrow  {C}_{\hol}^{b}F(X,{\Cal D},A;W) $,

\item{2.2}
$ C_{\Cal H}^{b} (X,A,n)_{\Cal D} $,
$ C_{\Cal H}^{b} (X,A)_{\Cal D} $,
$ C_{\Cal H}^{b} (X,A) $,

\item{}
$ K_{\Cal H}^{b} (X,A,n)_{\Cal D} $,
$ K_{\Cal H}^{b} (X,A)_{\Cal D} $,
$ K_{\Cal H}^{b} (X,A) $,

\item{}
$ {\Cal D}_{\Cal H}^{b}(X,A,n)_{\Cal D} $,
$ K_{\hol}^{b}F(X,{\Cal D},A) $,
$ D_{\hol}^{b}F(X,{\Cal D},A) $,

\item{2.4}
$ {D}_{\Cal H}^{b}(X,A) $,
$ {D}_{\Cal H}^{b}(X,A)_{{\Cal D}} $,

\item{}
$ {\DR}_{X}^{-1} : {D}_{\Cal H}^{b}(X,A) \rightarrow  {D}
_{\Cal H}^{b}(X,A)_{{\Cal D}} $,

\item{2.5}
$ D^{b}\MHM(X,A) $,
$ C^{b}\MHM(X,A) $,
$ C_{\MHM}^{b}(X,A) $,

\item{}
$ C_{c}^{*}(X,\Lmd;G,W') $,
$ K_{c}^{*}(X,\Lmd;G,W') $,
$ D_{c}^{*}(X,\Lmd;G,W') $,

\item{}
$ K_{c}^{+,b}(X,\Lmd;G,W')_{\inj} $,
$ MF_{\hol}(X,{\Cal D};W') $,
$ {K}_{\hol}^{+,b}(X,{\Cal D},A;G,W') $,

\item{}
$ D_{c}^{b}(X,\Lmd;G_{\text{b\^ete}},W') $,
$ C^{b}\Perv(X, \Lmd;W') $,

\item{2.7}
$ \varepsilon  : D^{b}\MHM(X,A) \rightarrow  {D}_{\Cal H}^{b}(X,A)_{{\Cal D}} $,

\item{2.8}
$ f_{*} : {D}_{\Cal H}^{b}(X,A) \rightarrow  {D}_{\Cal H}^{b}(Y,A) $,
$ f_{*} : {D}_{\Cal H}^{b}(X,A)_{{\Cal D}} \rightarrow  {D}_{\Cal H}^{b}(Y,A)
_{{\Cal D}} $,

\item{2.9}
$ {C}_{\Cal H}^{b}(X,A;G)_{{\Cal D}} $,
$ \Dec G $,
$ \Decdual G $,

\item{A.1}
$ {F}_{n}^{c}({\Cal A}) $,

\item{A.2}
$ F({\Cal C}) $,
$ {\tilde F}(\Cal A) $,

\item{A.3}
$ \Dec G $,
$ \Decdual G $,

\item{A.4}
$ C^{*}({\Cal A};G,W) $,
$ K^{*}({\Cal A};G,W) $,
$ D^{*}({\Cal A};G,W) $,

\item{}
$ C^{*}({\Cal A};W) $,
$ K^{*}({\Cal A};W) $,
$ D^{*}({\Cal A};W) $,
$ D^{b}({\Cal A};G_{\text{b\^ete}},W) $,

\item{3.1}
$ {\Cal V}_{k} $,
$ {\Cal V}^{A}_{k} $,
$ X_{\Cal B} $,
$ {\Cal V}^{A}_{k}(X_{\Cal B}) $,
$ C^{*}({\Cal V}^{A}_{k}(X_{\Cal B})) $, 
$ K^{*}({\Cal V}^{A}_{k}(X_{\Cal B})) $,

\item{}
$ {\Cal V}^{A}_{k}(X\langle D \rangle_{\Cal B})^{ns} $,
$ C^{*}({\Cal V}^{A}_{k}(X\langle D \rangle_{\Cal B})^{ns}) $,
$ K^{*}({\Cal V}^{A}_{k}(X\langle D \rangle_{\Cal B})^{ns}) $,

\item{}
$ D^{*}({\Cal V}^{A}_{k}(X\langle D \rangle_{\Cal B}))$,
$ D^{*}({\Cal V}^{A}_{k}(X\langle D \rangle_{\Cal B})^{ns}) $,

\item{}
$ \iota^{A} : {\Cal V}_{k} \to {\Cal V}^{A}_{k} $,

\item{}
$ \bold{R}{\pi }_{b,*}^{A} : K^{+}({\Cal V}^{A}_{k}(X_{\Cal B})) \rightarrow 
K^{+}(X_{b}^{\an},A) $,

\item{3.3}
$ {\Cal C}^{H}_{X\langle D \rangle}A $,
$ \tilde{\Omega}^{\ssbull} _{X\langle D \rangle} $,
$ \tilde{A}_{X\langle D \rangle} $,
$ \tilde{\Bbb C}_{X\langle D \rangle} $,

\item{4.1}
$ {C}_{\Cal S}^{\ssbull}A^{H} $,

\item{5.1}
$ {\Cal O}^{wn}_{X} $,
$ X^{wn} $,

\item{5.3}
$ K_{X} $,

\bigskip\bigskip
\centerline{{\bf References}}

\bigskip

\item{[1]}
Beilinson, A, Notes on absolute Hodge cohomology, Contemporary Math. 55 
(1986) 35--68.

\item{[2]}
Beilinson, A., Bernstein, J. and Deligne, P., Faisceaux pervers, Ast\'erisque, 
vol. 100, Soc. Math. France, Paris, 1982.

\item{[3]}
Calrson, J.,
Extensions of mixed Hodge structures, in Journ\'ees de 
G\'eom\'etrie Alg\'ebrique d'Angers 1979, 
Sijthoff-Noordhoff Alphen a/d Rijn, 1980, pp. 107--128.

\item{[4]}
Deligne, P., Th\'eorie de Hodge I, Actes Congr\`es Intern. Math. I (1970), 
425--430; II, Publ. Math. IHES, 40 (1971), 5--58; III, ibid. 44 (1974), 
5--77.

\item{[5]}
\SameAuthor, D\'ecompositions dans la cat\'egorie d\'eriv\'ee, Proc. Symp. 
Pure Math., 55 (1994), Part 1, 115--128.

\item{[6]}
Du Bois, Ph, Complexe de de Rham filtr\'e d'une vari\'et\'e singuli\`ere,  
Bull. Soc. Math. France, 109 (1981),
41--81.

\item{[7]}
Grothendieck, A and Dieudonn\'e, J., El\'ements de G\'eom\'etrie
Alg\'ebrique,  Publ. Math. IHES 32 (1967).

\item{[8]}
Grauert, H. and Riemenschneider, O., Verschwindungss\"atze f\"ur analytische 
Kohomologiegruppen auf
Komplexen R\"aumen,  Inv. Math. 11 (1970), 263--292.

\item{[9]}
Guill\'en, F., Navarro Aznar, V., Pascual-Gainza, P. and Puerta, F., 
Hyperr\'esolutions cubiques et descente cohomologique, Lect. Notes in 
Math., Springer, Berlin, vol. 1335, 1988.

\item{[10]}
Koll\'ar, J., Higher direct images of dualizing sheaves, I, II, Ann. 
of Math. 123 (1986), 11--42; 124 (1986), 171--202.

\item{[11]}
\SameAuthor, Shafarevich maps and automorphic forms, Princeton Univ. Press,
1995.

\item{[12]}
Kov\'acs, S., Rational, log canonical, Du Bois singularities : 
On the conjectures of Koll\'ar and Steenbrink, preprint.

\item{[13]}
Morgan, J.,
The algebraic topology of smooth algebraic varieties,
Publ. IHES 48 (1978), 137--204.

\item{[14]}
Saito, M., Modules de Hodge polarisables, 
Publ. RIMS, Kyoto Univ., 24 (1988), 849--995.

\item{[15]}
\SameAuthor, Mixed Hodge Modules, 
Publ. RIMS Kyoto Univ. 26 (1990), 221--333.

\item{[16]}
\SameAuthor, On Koll\'ar's conjecture,
Proc. Symp. Pure Math. 52 (1991), Part 2, 509--517.

\item{[17]}
\SameAuthor, On $b$-function, spectrum and rational singularity, 
Math. Ann. 295 (1993), 51--74.

\item{[18]}
\SameAuthor, On the Hodge filtration of Hodge Modules, 
preprint RIMS-1078, May 1996.

\item{[19]}
\SameAuthor, Induced $ D $-Modules and differential complexes, Bull. Soc.  Math. France, 117, (1989), 361--387.

\item{[20]}
Steenbrink, J., Mixed Hodge structures associated with isolated 
singularities, Proc. Symp. Pure Math. 40 (1983) Part 2, 513--536.

\item{[21]}
Verdier, J.L., Cat\'egories d\'eriv\'ees, in SGA 4 1/2, Lect. Notes in Math., 
Springer, Berlin, vol. 569, 1977 , pp. 262--311.

\bye